\documentclass[11pt,reqno]{amsart}

\usepackage[mathscr]{eucal}
\usepackage{mathrsfs}
\usepackage{amsfonts}
\usepackage{amsmath}
\usepackage{amsthm}
\usepackage{amssymb}
\usepackage{latexsym}
\usepackage[all]{xy}
\usepackage{pstricks,pst-node}

\theoremstyle{plain}
\newtheorem{defn}[equation]{Definition}
\newtheorem{cor}[equation]{Corollary}
\newtheorem{lem}[equation]{Lemma}
\newtheorem{prop}[equation]{Proposition}
\newtheorem{thm}[equation]{Theorem}
\newtheorem{conj}[equation]{Conjecture}
\newtheorem{problem}[equation]{Problem}

\theoremstyle{remark}
\newtheorem{claim}[equation]{Claim}
\newtheorem{rem}[equation]{Remark}
\newtheorem{examp}[equation]{Example}

\newtheorem{notation}[equation]{Notation}
\numberwithin{equation}{subsection}

\makeatletter
\renewcommand{\subsection}{\@startsection{subsection}{2}{0pt}{-3ex
plus -1ex minus -0.2ex}{-2mm plus -0pt minus
-2pt}{\normalfont\bfseries}} \makeatother

\newcommand{\erem}{\hphantom{.}\hfill$\lozenge$\end{rem}}

\newcommand{\Lmod}[1]{#1\text{-}{\mathsf{Mod}}}
\newcommand{\Rmod}[1]{\mathsf{Mod}\text{-}#1}
\newcommand{\bimod}[1]{#1\text{-}{\sf{Bimod}}}
\newcommand{\Lmof}[1]{#1\text{-}{\mathsf{mod}}}
\newcommand{\Rmof}[1]{\mathsf{mod}\text{-}#1}
\newcommand{\bimof}[1]{#1\text{-}{\sf{bimod}}}

\newcommand{\fpa}[1]{{\frac{\pa}{\pa x_{#1}}}}
\newcommand{\paf}[1]{{\frac{\pa\Phi}{\pa x_{#1}}}}
\newcommand{\hdot}{{\:\raisebox{2pt}{\text{\circle*{1.5}}}}}
\newcommand{\idot}{{\:\raisebox{2pt}{\text{\circle*{1.5}}}}}

\DeclareMathOperator{\Tor}{\mathrm{Tor}}
\DeclareMathOperator{\Iso}{\mathrm{Isom}}
\DeclareMathOperator{\Ext}{\mathrm{Ext}}
\DeclareMathOperator{\RHom}{\mathrm{RHom}}
\DeclareMathOperator{\coh}{\mathrm{Coh}}

\DeclareMathOperator{\Lie}{\mathrm{Lie}}
\DeclareMathOperator{\Tr}{\mathrm{Tr}}
\DeclareMathOperator{\Rep}{\mathrm{Rep}}

\DeclareMathOperator{\Irrep}{\mathsf{Irrep}}
\DeclareMathOperator{\Ad}{\mathrm{Ad}}
\DeclareMathOperator{\act}{\mathtt{act}}
\DeclareMathOperator{\La}{{\mathsf{\Lambda}}}
\DeclareMathOperator{\cone}{{\mathsf{Cone}}}

\DeclareMathOperator{\im}{\mathtt{image}}

\DeclareMathOperator{\sym}{{\mathrm{Sym}}}

\DeclareMathOperator{\loc}{{{\scr L}\!{\text{\it oc}}_\dd}}
\DeclareMathOperator{\conn}{{\mathrm{Conn}}_\dd}
\renewcommand{\div}{{\mathrm{div}^{\,}}}

\def\map{\longrightarrow}
\newcommand{\dis}{\displaystyle}
\newcommand{\step}[1]{\noindent\vskip 2pt{{\text{\sc{Step }}#1}\en}}

\newcommand{\bex}{\vskip 7pt $\hphantom{x}\hskip 40mm\dis }
\newcommand{\eex}{\hfill$\lozenge$\break}
\newcommand{\beq}{\begin{equation}\label}
\newcommand{\eeq}{\end{equation}}

\DeclareMathOperator{\Spec}{\mathrm{Spec}}
\DeclareMathOperator{\pr}{pr}
\newcommand{\iso}{{\,\stackrel{_\sim}{\to}\,}}
\newcommand{\cd}{\!\cdot\!}
\newcommand{\qis}{\stackrel{\oper{qis}}\too }

\DeclareMathOperator{\GL}{{\mathrm{GL}}}
\DeclareMathOperator{\PGL}{{\mathrm{PGL}}}

\def\ccirc{{{}_{\,{}^{^\circ}}}}

\newcommand{\g}[1]{\mathfrak{#1}}
\newcommand{\scr}[1]{\mathscr{#1}}

\renewcommand{\b}[1]{\mathbf{#1}}

\newcommand{\Coh}{\mathrm{Coh}} 

\DeclareMathOperator{\End}{\mathrm{End}}

\newcommand{\op}{\operatorname}

\DeclareMathOperator{\Hom}{\mathrm{Hom}}

\DeclareMathOperator{\Der}{\mathrm{Der}}
\DeclareMathOperator{\dder}{\mathbb{D}\mathrm{er}}

\newcommand{\y}{^\diamond }

\newcommand{\id}{\mathrm{Id}}
\newcommand{\Id}{\mathrm{Id}}

\newcommand{\Mat}{\mathrm{Mat}}
\newcommand{\KK}{{\mathbb K}}
\newcommand{\CK}{{\mathcal K}}
\newcommand{\xxx}{x_1,\ldots,x_n }
\newcommand{\KA}{{\mathcal K}_A }
\newcommand{\perf}{{\mathsf{perf}}}
\newcommand{\fpr}{{\mathsf{fproj}}}

\newcommand{\fin}{{\mathsf{finj}}}
\newcommand{\dash}{\text{-}}
\newcommand{\crit}{{\mathsf{crit}}}
\newcommand{\obs}{{\mathsf{Obstruct}}}
\newcommand{\BD}{{\mathbb D}}
\newcommand{\BS}{{\mathbb S}}
\newcommand{\bth}{{\boldsymbol\xi}}

\newcommand{\LOM}{{{\mathbb L}\Om^1_\idot(F,\Phi)}}
\newcommand{\nat}{^\flat}

\newcommand{\De}{{\boldsymbol{\delta}}}
\newcommand{\br}{_{\operatorname{cyc}}}

\newcommand{\CZ}{{\mathcal Z}}

\newcommand{\TTT}{{\mathbb T}}
\newcommand{\out}{\stackrel{_{\oper{out}}}\otimes }
\newcommand{\ino}{\stackrel{_{\oper{in}}}\otimes }
\newcommand{\bw}{\De}
\newcommand{\VV}{{\mathfrak V}}

\newcommand{\si}{\sigma}

\newcommand{\G}{\Gamma}

\renewcommand{\d}{{\mathfrak d}}
\newcommand{\Ph}{\Phi}
\newcommand{\spad}{{^{_{_\heartsuit}}}}
\newcommand{\eps}{\varepsilon}

\newcommand{\arr}{\overset{{\,}_\to}}
\newcommand{\bi}{\imath}
\newcommand{\bl}{{\scr L}}
\newcommand{\boplus}{\mbox{$\bigoplus$} }
\newcommand{\botimes}{\mbox{$\bigotimes$}}

\newcommand{\X}{{\mathfrak{X}}}

\newcommand{\vir}{^{\operatorname{vir}}}

\newcommand{\van}{{{\scr V}\!an}}
\newcommand{\opp}{{\operatorname{op}}}
\newcommand{\bj}{{\mathbf j}}
\newcommand{\bp}{{\mathbf p}}

\DeclareMathOperator{\Inn}{\mathrm{Inn}}
\DeclareMathOperator{\Aut}{\mathrm{Aut}}

\renewcommand{\mid}{\enspace\big|\enspace}

\newcommand{\CF}{{\mathcal F}}
\newcommand{\CE}{{\mathcal E}}
\newcommand{\DZ}{{\mathcal Z}D^b}
\newcommand{\scc}{{\scr C}}
\renewcommand{\o}{\otimes }
\newcommand{\be}{\beta}
\newcommand{\cyd}{{ $\oper{CYd}$ }}

\newcommand{\wt}{\widetilde }

\newcommand{\grep}{{\g g(\Rep_\dd }}
\DeclareMathOperator{\Sym}{{\mathrm{Sym}}}

\DeclareMathOperator{\ad}{\mathrm{ad}}

\DeclareMathOperator{\imm}{{\mathrm{Im}}^{\,}\bi^\al}

\DeclareMathOperator{\Ker}{\mathrm{Ker}}

\newcommand{\DG}{$\,\oper{DG}\;$}

\def\BB{{\mathsf{B}}}

\newcommand{\brst}{D^{^{_{\footnotesize{\mathtt{BRST}}}}\!}}

\newcommand{\CHI}{\chi}
\newcommand{\QQ}{{\overline{Q}}}
\newcommand{\QQQ}{{\widehat{Q}}}
\newcommand{\btr}{{\blacktriangle}}

\newcommand{\DPA}{\dder(F|A) }

\newcommand{\dd}{{\mathbf d}}
\newcommand{\SP}{{\mathsf{P}}}

\newcommand{\bone}{{\boldsymbol{1}}}

\newcommand{\bll}{\{\!\{}
\newcommand{\brr}{\}\!\}}

\newcommand{\mto}{\mapsto}
\newcommand{\lto}{\longmapsto}

\newcommand{\inv}{^{-1}}
\newcommand{\vi}{${\en\sf {(i)}}\;$}
\newcommand{\vii}{${\;\sf {(ii)}}\;$}
\newcommand{\viii}{${\sf {(iii)}}\;$}
\newcommand{\iv}{${\sf {(iv)}}\;$}

\newcommand{\sset}{\subset}

\newcommand{\into}{{}^{\,}\hookrightarrow^{\,}}
\newcommand{\too}{\,\longrightarrow\,}
\newcommand{\onto}{\twoheadrightarrow}
\newcommand{\tooo}{{\;{-\!\!\!-\!\!\!-\!\!\!-\!\!\!\longrightarrow}\;}}
\newcommand{\eqq}{{=\!\!\!=\!\!\!=\!\!\!=\!\!\!=\!\!\!=\!\!\!=\!\!\!=}}
\newcommand{\oper}{\operatorname}

\newcommand{\C}{\mathbb C}

\newcommand{\Z}{\mathbb Z}

\newcommand{\wb}{\overline }

\newcommand{\HH}{{\mathbb{H}}}
\newcommand{\SH}{{\mathbb D}}

\renewcommand{\th}{\theta}

\newcommand{\bmu}{{\boldsymbol{\mu}}}
\newcommand{\bt}{{\boldsymbol{\alpha}}}
\newcommand{\Om}{\Omega}

\newcommand{\pp}{{\mathbf{p}}}
\newcommand{\Th}{\Theta}
\renewcommand{\P}{\xi}

\DeclareMathOperator{\gr}{\mathrm{gr}}

\def\ip<#1,#2>{\left\langle#1,#2\right\rangle}
\def\sp<#1>{\left\langle#1\right\rangle}

\newcommand{\ve}{^\vee }

\newcommand{\llb}{{\boldsymbol(}\!^{\!}{\boldsymbol(}}
\newcommand{\rrb}{{\boldsymbol)}\!^{\!}{\boldsymbol)}}

\newcommand{\lo}{\,{\stackrel{L}\otimes}}

\def\ip<#1,#2>{\left\langle#1,#2\right\rangle}

\def\pbox{\parbox[t]{130mm}}
\def\npb{\noindent$\bullet\quad$\parbox[t]{143mm}}

\def\hp{\hphantom{x}}

\newcommand{\half}{\mbox{$\frac{1}{2}$}}

\newcommand{\bd}{{\mathbf{d}}}

\newcommand{\om}{\omega}
\newcommand{\pa}{\partial}

\newcommand{\omh}{{\mathbf{w}}}
\newcommand{\ka}{{\boldsymbol{\kappa}}}

\newcommand{\al}{{\alpha}}

\newcommand{\la}{{\lambda}}
\newcommand{\A}{{\mathfrak A}}
\newcommand{\FD}{{\mathfrak D}}
\newcommand{\ffd}{\wh{\mathfrak D}}
\newcommand{\FI}{{\mathfrak I}}
\newcommand{\en}{{\enspace}}

\newcommand{\wh}{\widehat}

\def\oo{{\mathcal O}}

\def\U{{\mathcal U}}

\def\aa{{\mathcal{A}}}

\def\gl{{\mathfrak{g}\mathfrak{l}}}

\def\k{{\mathbb C}}

\def\gld{{\mathrm{GL}_\dd}}

\def\SL{{\mathrm{SL}}}

\newcommand{\mm}{M}

\newcommand{\nn}{1,\ldots,n}

\def\T{{\mathcal{T}}}

\begin{document}

\title{\Large{\textbf{CALABI-YAU ALGEBRAS}}}

\author{{\sc Victor Ginzburg}}

\maketitle
\vskip -25mm

\begin{abstract} We introduce some new algebraic
structures arising naturally in the geometry of 
CY manifolds and mirror symmetry. We give a universal
construction of CY algebras in terms of a  noncommutative  symplectic
DG algebra  resolution.

In dimension 3, the resolution is determined
by a {\em  noncommutative  potential}. Representation varieties of the CY
algebra  are intimately related to the set of 
critical points, and to the sheaf of vanishing cycles of 
the potential.
Numerical invariants, like ranks of cyclic homology
groups, are expected to be
given by `matrix integrals' over
representation varieties.

We discuss examples of CY algebras involving   quivers,
 3-dimensional McKay
correspondence, crepant resolutions,  Sklyanin algebras,   hyperbolic
3-manifolds and Chern-Simons.
Examples related to quantum Del Pezzo surfaces are discussed in ~\cite{EtGi}. 
\end{abstract}

\bigskip

\centerline{\sf Table of Contents}
\vskip -1mm

$\hspace{20mm}$ {\footnotesize \parbox[t]{115mm}{
\hp${}_{}$\!\hp1.{ $\;\,$} {\tt Algebras and potentials}\newline
\hp2.{ $\;\,$} {\tt Representation functor, critical points and
vanishing cycles}\newline
\hp3.{ $\;\,$} {\tt Calabi-Yau condition}\newline
\hp4.{ $\;\,$} {\tt Quiver algebras and McKay correspondence in
dimension 3}\newline
\hp5.{ $\;\,$} {\tt Calabi-Yau algebras of dimension  3}\newline
\hp6.{ $\;\,$} {\tt Fundamental groups,  Chern-Simons, and  the Weil
 representation}\newline
\hp7.{ $\;\,$} {\tt Calabi-Yau algebras and Calabi-Yau manifolds}\newline
\hp8.{ $\;\,$} {\tt Noncommutative Hessian}\newline
\hp9.{ $\;\,$} {\tt Some homological algebra}
}}

{\large\section{{\textbf{Algebras and potentials}}}\label{pc}}
\subsection{Introduction} In this paper, we study some new algebraic
structures, such as  Calabi-Yau (CY) algebras,
arising naturally in the geometry of CY manifolds.
The ultimate goal of introducing CY algebras is to transplant
most of  conventional CY geometry to the setting of
{\em  noncommutative } geometry.  Some motivation for this, to be
explained in \S\ref{motiv}, comes  from
mirror symmetry. 
Besides that, it will be demonstrated by numerous concrete examples
that CY algebras do arise `in nature'. Furthermore,
CY algebras enjoy quite intriguing homological properties
which are closely related to  algebro-geometric properties of the corresponding
representation schemes, see \S\ref{rep_sec}.

In  this section \ref{pc},
we begin with  most essential, elementary algebraic  constructions,
and introduce them  in  as `ground to earth' a
way as possible. Thus, our exposition won't be
strictly  logical; the definitions and constructions
of \S\ref{pc} will not be given in full generality.
A more general and more conceptual approach, as well as motivation for
these constructions,
 will be provided in subsequent sections \ref{rep_sec}, \ref{cy_sec},
 and \ref{cy3}.

Among our most important results, 
we would like to mention Theorem \ref{BV}, Theorem \ref{main}, Theorem
 \ref{mca},
and Theorem \ref{CYconditions}.


The point of view on Calabi-Yau algebras
advocated in this paper is not the only possible one.
The reader is  referred to \cite{Bo},
\cite{CR}, \cite{IR}, and \cite{KR} for alternative
approaches.

\subsection{Acknowledgements}{\small
I am very much indebted to Maxim Kontsevich for 
generously sharing with me his unpublished ideas.
Maxim   Kontsevich has  read a preliminary version of
this paper and made a number of very interesting comments
which are incorporated in the present version.

I have 
benefited a lot from useful discussions and E-mail exchanges with  Kevin Costello, 
 Raphael Rouquier,
and  Michel Van den Bergh. Thanks are also due to Pavel Etingof, Dmitry
Kaledin,
Nikita Nekrasov,
Viktor Ostrik, Travis Schedler,
Amnon Yekutiely, and James Zhang.}

\subsection{Algebras defined by a potential.}\label{pot}
We will work with unital associative $\k$-algebras.
Given such an algebra $A$, let $[A,A]$ denote  the $\C$-vector 
subspace in $A$ spanned by the commutators and
write $A\br=A/[A,A]$ for  the commutator
quotient space.

 Let  $F=\k\langle \xxx\rangle$ be a free associative algebra
with  $n$
generators. The vector space
$F\br$ has an obvious basis labelled by
 cyclic words in the alphabet $\xxx$.

For each $j=1,\ldots,n,$ M.~Kontsevich  \cite{Ko1} introduced 
a linear map
$\frac{\pa}{\pa x_j}: F\br\to F,$
$\Phi\mto \frac{\pa \Phi}{\pa x_j},$ as follows.
Given a cyclic word $\Phi=x_{i_1}x_{i_2}\ldots x_{i_r}$,
one finds all occurrences of the variable $x_j$ in $\Phi$.
Deleting one such  occurrence of $x_j$  breaks up the cycle $\Phi$,
thus creating an ordinary (not cyclic) word. We carry out this procedure
with each occurrence of the variable $x_j$ in $\Phi$, one at a time.
The element $\frac{\pa \Phi}{\pa x_j}$ is defined to be the
sum of all the resulting words. More formally, we put
\beq{partial}
\paf{j}:=\sum_{\{s\;|\; i_s=j\}}x_{i_s+1}\,x_{i_s+2}\ldots  x_{i_r}\; x_{i_1}\,x_{i_2}\ldots
x_{i_s-1}\in \k\langle \xxx\rangle,
\eeq
and extend this definition to linear combinations of cyclic
words by $\C$-linearity.

Thus, given any element $\Phi\in F\br$, to be referred to as {\em
potential}, we have a well defined collection of elements
$\frac{\pa \Phi}{\pa x_i}\in F,\,i=1,\ldots,n.$

A key object of interest for us is an associative algebra
\beq{AF}
\A(F,\Phi):=\C\langle  x_1,\ldots,x_n\rangle\big/\llb \pa \Phi/\pa  x_i\rrb_{i=1,\ldots,n},
\eeq
a quotient of the free algebra $F$ by the two-sided ideal
generated by all $n$ partial derivatives of the potential $\Phi$.

\begin{examp}[Basic example]\label{basic} Let $F =\C\langle x,y,z\rangle,$ and let
$\Phi=xyz-yxz\in F\br,$ a difference of two cyclic words.
We compute
$$
\mbox{\Large{$\frac{\pa \Phi}{\pa x}$}}
= yz-zy,\quad
\mbox{\Large{$\frac{\pa \Phi}{\pa y}$}}= zx-xz,\quad
\mbox{\Large{$\frac{\pa \Phi}{\pa z}$}}= xy-yx.
$$
Therefore, taking the quotient modulo the ideal generated by the above elements
produces a polynomial algebra, i.e., we get
\beq{bas_iso}
\A(\C\langle x,y,z\rangle, \Phi)=\C[x,y,z].
\eeq
\end{examp}

A key feature of the above example is that we work 
with a polynomial algebra, $\C[x,y,z],$ in 
$\mathbf{3}$ variables.  The algebra $\C[x,y,z]$ is, in effect,
a basic
example of  CY algebra of dimension 3.  One of the main messages of this paper is that,
roughly speaking:
\vskip 2pt

\quad\pbox{\textbf{\em Any Calabi-Yau algebra of dimension 3 
`arising in nature' is defined by a potential, i.e. has the form
$\A(F,\Phi)$.}}

\vskip 2pt

\noindent
We refer to \S\ref{CY_defi} for the definition of CY algebras, and to \S\ref{cot} for
 a more precise version
 of the claim above. 
\begin{rem} \vi Not every algebra of the  form $\A(F,\Phi)$,
where $F=\C\langle\xxx\rangle$,
is a CY algebra of dimension 3.
There seems to be no simple characterisation of those
potentials $\Phi$ for which  $\A(F,\Phi)$ is a CY algebra of dimension
 3, cf. however Theorem \ref{CYconditions}.

\vii It should be emphasized that, for a given
algebra $A$, its presentation in the form $A=\A(F,\Phi)$
is by no means determined by the algebra itself.
A given algebra may have many different presentations involving
different free algebras $F$ and different potentials ~$\Phi$.

\viii Let $\Aut(F)$ be the group of algebra
automorphisms $f: x_1\mto f(x_1),\ldots, x_n\mto f(x_n),$ of
the  free algebra $F=\C\langle\xxx\rangle$.
The group   $\Aut(F)$, as well as
the Lie algebra $\Der(F,F)$ of all derivations $F\to F$, acts naturally on
 the commutator quotient space $F\br$. Potentials
from the same $\Aut(F)$-orbit clearly give rise to isomorphic algebras
$\A(F,\Phi)$.

A potential $\Phi\in F\br$ is said to be {\em isolated}
if the corresponding  $\Aut(F)$-orbit is `infinitesimally open' in the
sense
that we have $\Der(F,F)(\Phi)=F\br$. It is  an interesting open
problem
to describe all ($\Aut(F)$-orbits of)  isolated potentials $\Phi$ such that
 $\A(F,\Phi)$ is a CY algebra of dimension
 $3.$ It seems likely that the CY algebras  arising in the theory of 
 cluster categories, see \cite{IR}, \cite{KR},
provide examples  of such  CY algebras of dimension
 $3$ associated to an isolated potential.
\erem
 
Below, we provide a few more 
interesting examples of algebras defined by a potential. The philosophy
behind them is to search for deformations of the
polynomial algebra, like in  Example \ref{basic},
obtained by an appropriate deformation of the potential.

\begin{examp}\label{qexamp} Fix $q\in\C$ and let $\Phi=xyz-q\cdot yxz-f$
for some $f\in  F\br,$ where $F=\C\langle x,y,z\rangle,$ as before. 
Then, the corresponding algebra $\A(F,\Phi)$ is a quotient of the free
algebra $\C\langle x,y,z\rangle$ by the relations
\beq{abc}
xy-q\cd yx=\pa f/\pa z,\quad zx-q\cd xz=\pa f/\pa y,\quad
yz-q\cd zy=\pa f/\pa x.
\eeq

In particular,  put
$$\Phi_q=xyz-q\cd yxz+\mbox{\large{$\frac{1}{2}$}}(x^2+y^2+z^2).
$$
For $q=1$, we get
$\A(F,\Phi_1)={\mathcal U}(\mathfrak{s}\mathfrak{l}_2),$
 the  enveloping algebra of the Lie algebra
$\mathfrak{s}\mathfrak{l}_2$.

As another example, in \eqref{abc}, put  $f\in\C[x],$
a polynomial independent of the variables $y$ and $z$.
Further, take $q$ to be a primitive $n$-th root of unity.
Then the algebra
$\A(F,\Phi)$ has a large center, $\CZ(F,\Phi):=\CZ(\A(F,\Phi)).$ Specifically,
the elements
 $\b x:=x^n,\,\b y:=y^n,\,\b z:=z^n$ are  central.

Assume in addition
that,  in the expansion
$f=\sum_{r=1}^d a_r\cdot x^r$, the coefficient $a_r\in\C$ 
 vanishes whenever $n|r$. Then,
the center $\CZ(F,\Phi)$ contains  an extra element 
$\b u:=\sum_{r=1}^d {\mbox{\Large$\frac{r\cdot a_r}{1-q^r}$}\cdot x^r.}$
One can show that the center is, in effect,
generated by the  elements $\b x,\b y,\b z,\b u$. Moreover, we have
$\CZ(F,\Phi)\cong\C[\b x,\b y,\b z,\b u]/\llb \b x\cdot\b y\cdot\b z-\phi(\b x, 
\b u)\rrb,$ where $\phi\in\C[\b x,\b u]$ is a certain (complicated) polynomial 
in two variables which is determined by ~$f$. For more discussion,
cf. Example
\ref{qexamp1} and also \cite{EtGi}.
\end{examp}

\begin{examp}[Sklyanin algebras]\label{sklya} Keep $F =\C\langle x,y,z\rangle$.
Given a triple $a,b,c\in\C,$  we  put
$$\Phi=axyz+byxz+c(x^3+y^3+z^3).$$

The corresponding algebra $\A(F,\Phi)$ is 
called 
 {\em quadratic Sklyanin algebra}. The original definition
in
[AS], [ATV], [VdB2] was not given in this form but is  equivalent to
ours, as can be readily seen
from e.g.  \cite[formula (8)]{VdB2}.

The {\em cubic} Sklyanin algebra  (of type $A$) may be defined as an algebra
of the form
$$\A(\C\langle x,y\rangle,\Phi),\qquad\Phi=ax^2y^2+bxyxy+c(x^4+y^4).
$$

If $a\neq b$, the same  algebra has the following different presentation, see \cite{Ma},
$$\A(\C\langle x,y,z\rangle,\,\Phi_{p,q}),\qquad
\Phi_{p,q}=\half z^2+xyz-yxz+p(x^2y^2+xyxy)+q(x^4+y^4).
$$
Here, the parameters $p,q\in\C$ are related to parameters $a,b,c$ by the equations
$a=-r(p+1),\,b=r(2-p),\,c=-rq,$
for some $r\neq 0$.

Note that, for $p=q=0$, the algebra  $\A(F,\Phi_{p,q})$
degenerates to the enveloping algebra of the 3-dimensional
Heisenberg Lie algebra such that $z$ is a central element
and $[x,y]=z$.
%
\end{examp}

\subsection{A DG algebra.}\label{dg_sec} We will 
denote the  grading on a DG algebra by either upper or
lower index depending
on whether the differential in the DG algebra has degree
$+1$ or $-1$, respectively.
From now until the end of \S\ref{pc}, we let
$F=\C\langle  x_1,\ldots,x_n\rangle.$ 

We will view the algebra
$F$ as the degree zero component in  a  free {\em graded} algebra
\beq{Dalg}
\FD=\C\langle  x_1,\ldots,x_n,\th_1,
\ldots,\th_n,t\rangle,\en \deg
t=2,\en\deg \th_j=1,\en j=1,\ldots,n.
\eeq

Any potential $\Phi\in F\br$ gives rise to
a differential $\d: \FD_\idot\to\FD_{\idot-1}$
defined as follows
\beq{d_def}
\d:\ t\mto\sum_{j=1}^n [x_j,\th_j],\quad
\th_j\mto\pa \Phi/\pa x_j,\quad x_j\mto0,\en \forall j=1,\ldots,n.
\eeq
This  assignment on generators can be uniquely extended to  an odd
super-derivation 
on  ~$\FD.$ 

Using Proposition \ref{frob}(i) below, one checks that $\d^2=0$. We write
$\big(\FD_\idot(F,\Phi),\,\d\big)$ for the resulting DG algebra.
It is immediate from  formulas \eqref{d_def} that
the zeroth homology of  $\FD(F,\Phi)$ is nothing
but the algebra defined by our potential $\Phi$; thus, one has
a diagram
 \beq{HD}
\xymatrix{
\FD(F,\Phi)=\boplus_{_{r\geq 0}}\,\FD_r(F,\Phi)
\ar@{->>}[r]&\FD_0(F,\Phi)
\ar@{->>}[r]&
 H_0\big(\FD(F,\Phi),\,\d\big)=\A(F,\Phi).
}
\eeq

The geometric meaning of the DG algebra $\FD_\idot(F,\Phi)$ will
be clarified in \S\ref{cy_rep}.
Specifically,  let $G$ be a Lie group and $X$ a $G$-manifold.
It will explained that  any $G$-invariant smooth function
$\phi$ on $X$
 gives rise to a complex, called  Batalin-Vilkovisky complex,
whose cohomology are related to critical points of the function $\phi$.
We will show that  the 
above defined DG algebra $\FD_\idot(F,\Phi)$ is nothing but a  noncommutative  version
of the   Batalin-Vilkovisky complex.

One of the main results of this paper 
says that (a suitable completion of) the DG algebra $\big(\FD_\idot(F,\Phi),\,\d\big)$ is 
acyclic in positive degrees,  that is, the projection
\eqref{HD} is a quasi-isomorphism, iff
$\A(F,\Phi)$ is a CY algebra of dimension 3, cf. Theorem \ref{CYconditions}.
This result is at the origin of an especially nice homological
behavior of  CY algebras of dimension 3.

\subsection{Noncommutative calculus} To proceed further, we need
to introduce a few basic concepts of  noncommutative  geometry.

Write $\o=\o_\C$.
For any algebra $A$, the space $A\otimes A$  has two
{\em commuting} (and isomorphic) $A$-bimodule structures,
called  the {\em outer}, resp., {\em inner}, bimodule structure.
These two   bimodule structures are given by
\beq{inout} a(b'\otimes b'')c:=(ab')\otimes (b''c),\quad\text{resp.},\quad
a(b'\otimes b'')c:=(b'c)\otimes(ab''),\quad\forall
b',b'',a,c\in A.
\eeq
We will view $A\otimes A$  as a bimodule  with
respect to the outer structure, unless specified otherwise;
 thus $A\otimes A=A\out A$  is a rank one free
$A$-bimodule.

Let $\bimod{A}$ be the  category of  $A$-bimodules. There is a
contravariant duality functor
\beq{DUAL}
\Hom_{\bimod{A}}(-,A \o A):\
\bimod{A}\too\bimod{A},\qquad M\lto M^\vee.
\eeq
Here the target bimodule $A\o A$
is taken  with respect to the outer  structure.
The inner   structure  on $A\o A$ survives in
the $\Hom$-space, making it an $A$-bimodule again.

The bimodule $\Om^1A$, of {\em noncommutative differentials},
is defined as the kernel of multiplication map $\bmu: A\o A\to A$.
The  dual, $\dder A:=(\Om^1A)^\vee,$
is called the bimodule  of {\em double derivations},
a noncommutative counterpart of the space of vector fields on
a manifold.

Elements  of $\dder A$ may be identified with derivations $A\to A\o A$
as follows. First of all, there is a canonically defined
\footnote{Note that there is no {\em  canonical} ordinary derivation $A\to A$.} distinguished
double derivation (de Rham differential)
\beq{De}
\De: A\to\Om^1A\into A\o A,\quad a\mto da:=1\o a - a\o 1\in \Om^1A\sset A\o A.
\eeq
Now, the double derivation corresponding
to an element $\th\in \Hom_{\bimod{A}}(\Om^1A,A \o A)$ is 
given by the assignment $A\ni a\mto \th(da).$ 

Associated with each
element  $a'\o a''\in A\o A,$ there is an {\em inner}
double derivation $\ad(a'\o a''): u\mto ua'\o a''-a'\o a''u.$
By definition, one has $\De=\ad(1\o1).$ 

Clearly, we have $(A\out A)^\vee=A\ino A$. The duality
functor, $M\mto M^\vee,$ interchanges   $\jmath_A$, the tautological $A$-bimodule
imbedding below,
with the  $A$-bimodule map $\ad: a' \o a''\mto \ad(a'\o a'')$, as follows
\beq{morph}
\jmath_A:\ \Om^1A\into A\out A\qquad
\xymatrix{
\ar@{.>}[rr]^<>(0.5){\text{duality}}_<>(0.5){\eqref{DUAL}}&&}
\qquad\ad=\jmath_A^\vee:\ A\ino A\to \dder A.
\eeq

\begin{examp}
The bimodule $\Om^1F$ of 1-forms for the free algebra  $F=\C\langle\xxx\rangle$
is a rank $n$ free $F$-bimodule with basis $dx_1,\ldots,dx_n$.
In this basis, one can write the double derivation  \eqref{De}, for $A=F$, in the form
\beq{Ombasis}\De:\ F\too \Om^1F=\bigoplus_{j=1}^n F\cd dx_j\cd F,\qquad
f\mto df=\sum_{j=1}^n \left(\mbox{\Large$\frac{\pa f}{\pa x_j}$}\right)'
\cd dx_j\cd
\left(\mbox{\Large$\frac{\pa f}{\pa x_j}$}\right)''.
\eeq

For each $j=\nn,$
the corresponding term in the sum in the RHS of the above formula
determines a certain element $\left(\frac{\pa f}{\pa x_j}\right)'\o
\left(\frac{\pa f}{\pa x_j}\right)''\in F\o F.$
Here and elsewhere,  we use Sweedler's notation and write
$u'\o u''$ instead of $\sum_i u'_i\o u''_i$ for an element
in a tensor product.

Thus, one obtains a collection of maps
\beq{dj}
\mbox{\Large$\fpa{j}$}:\ F\to F\o F,\qquad f\mto
\left(\mbox{\Large$\frac{\pa f}{\pa x_j}$}
\right)'\o
\left(\mbox{\Large$\frac{\pa f}{\pa x_j}$}\right)'',\quad j=\nn.
\eeq

It is immediate  that each of  these maps is a double derivation;
moreover, these double derivations form a basis of $\dder F,$
dual to the basis $\{dx_j,\,j=\nn\}$ on 1-forms, i.e.,
\beq{Dbasis}
\dder F=\bigoplus_{j=1}^n F\cd \frac{\pa}{\pa x_j}\cd F,\qquad
\fpa{j}(x_i)= \begin{cases} 1\o 1&\text{if}\en j=i\\
0&\text{if}\en j\neq i.
\end{cases}
\eeq

In terms of this basis,  the distinguished double derivation
\eqref{De} reads
\beq{inner}
\De=\ad(1\o 1)=\sum_{j=1}^n \left(x_i\cd\mbox{\Large$\fpa{i}$}-\mbox{\Large$\fpa{i}$}\cd
x_i\right)=\sum_{j=1}^n \left[x_i,\,\mbox{\Large$\fpa{i}$}\right].
\eeq
Here, for any $F$-bimodule $M$ and elements $x\in F,\,m\in M$,
we write $[x,m]:=xm-mx.$
\begin{rem} 
We  use the same notation
$\fpa{j}$ both in \eqref{dj} and in \eqref{partial} since
the effect of the action of   double derivations \eqref{dj}
on  a (non cyclic) word 
is similar to that of formula  \eqref{partial}.
\erem

\noindent
{\sc Higher derivatives.}\en To any potential $\Phi\in F\br,$ one
associates its {\em Hessian},
 a noncommutative
version of  the Hessian of  a smooth
function on  $\k^n$.
The noncommutative Hessian,  $\|\frac{\pa^2\Phi}{\pa x_i\pa x_j}\|$, is an $n\times n$-matrix with
entries in $F\o F$; its $(i,j)$-th
entry is defined as the image of $\Phi$ under the composite map,
cf. \eqref{partial} and \eqref{dj},
$$\xymatrix{
F\br\ar[r]^<>(0.5){\frac{\pa}{\pa x_j}}&
F\ar[r]^<>(0.5){\frac{\pa}{\pa x_i}}&F\o F},
\qquad
\Phi\mto \mbox{\Large$\frac{\pa^2\Phi}{\pa x_i\pa x_j}=
(\frac{\pa^2\Phi}{\pa x_i\pa x_j})'\o (\frac{\pa^2\Phi}{\pa x_i\pa x_j})''$}\in F\o F.
$$

\begin{rem}\label{locus} The Hessian of a function $\phi$  on an arbitrary
manifold $X$
is well defined only on the  critical locus $\crit(\phi)\sset X$, 
 the zero scheme of the 1-form $d\phi$.
A  noncommutative  analogue of such a construction will be discussed
in \S\ref{nchess}.
\erem

Similarly, for each  $r\geq 1$,
there is a map $\fpa{i}: F^{\o r}\to F^{\o (r+1)}.$
Thus, for any $\Phi\in F\br$ and any $r$-tuple of indices,
one inductively defines the elements 
$\frac{\pa^r\Phi}{\pa x_{i_1}\pa x_{i_2}\ldots\pa x_{i_r}}\in F^{\o r}.$

The classic result on  {\em total} symmetry of the tensor of
$r$-th 
derivatives of a function gets replaced, in  noncommutative  geometry, by
{\em cyclic} symmetry. Specifically,  in $F^{\o r}$, one has
\beq{opp}
\si\left(\frac{\pa^r\Phi}{\pa x_{i_1}\pa x_{i_2}\ldots\pa
x_{i_r}}\right)
=
\frac{\pa^r\Phi}{\pa x_{\si(i_1)}\pa x_{\si(i_2)}\ldots\pa x_{\si(i_r)}},
\qquad\forall \Phi\in F\br,\;r\geq 1.
\eeq
In the LHS of this formula,  $\si$   denotes
 the map $F^{\o r}\to F^{\o r}$ given by 
the cyclic permutation of the
tensor factors while in the RHS $\si$ stands
for a cyclic permutation  of indices.
 
One can also prove  the following
\begin{prop}\label{frob} \vi {\sc(Poincar\'e lemma).}\en 
For an $n$-tuple  $\{f_i\in F\}_{i=\nn},$ we have
$$\sum_{i=1}^n [x_i,f_i]=0\quad
\Longleftrightarrow\quad
\exists \Phi\in F\br\en\oper{such\en that}\en f_i=\pa\Phi/\pa x_i,\en
i=\nn;
$$

\vii {\sc(Frobenius theorem).}\en 
Let $\|f_{i,j}\|$ be an $F\o F$-valued $n\times n$-matrix such
that
$$\si(f_{i,j})=f_{\si(i),\si(j)},
\quad\operatorname{and}\quad
\si\left(\frac{\pa f_{i,j}}{\pa x_k}\right)=
\frac{\pa f_{\si(i),\si(j)}}{\pa x_{\si(k)}}
,\quad\forall 1\leq i,j,k\leq n.
$$

Then,
there exists a potential $\Phi\in F\br,$ such that we have
$\|f_{i,j}\|=\|\frac{\pa^2\Phi}{\pa x_i\pa x_j}\|$.\qed
\end{prop}

The implication `$\Leftarrow$' in part (i) has been first
noticed by M. Kontsevich \cite[\S6]{Ko1}.
\end{examp}

\subsection{Noncommutative cotangent complex}\label{ncc}
Write $\T_X$ for the  tangent, resp. $\T^*_X$ for the
cotangent, sheaf of an algebraic variety (or scheme) $X$.
Given a scheme imbedding $Y\into X$ one gets, by restriction, two sheaves
$\T^*_{X}|_Y$ and $\T_{X}|_Y$, and also
the conormal sheaf, ${\scr N}^*_{X|Y}.$ There is a standard
 short exact sequence of sheaves on $Y$,
\beq{conormal}
\xymatrix{
0\ar[r]&{\scr N}^*_{X|Y}\ar[rr]&&\T^*_{X}|_Y\ar[rr]^<>(0.5){p_{X|Y}}&&\T^*_Y\ar[r]&0.
}
\eeq

The geometric setting above can be imitated in algebra.
A morphism of schemes corresponds to
 an algebra map $B\to A$.
Given such a map, we introduce the following
 $A$-bimodules which are algebraic counterparts of 
 $\T^*_{X}|_Y$ and $\T_{X}|_Y$, respectively,
\beq{B|A}
\Om^1(B|A):=A\o_B\Om^1B\o_B A\quad\text{resp.},\quad
\dder(B|A):=A\o_B\dder B\o_B A.
\eeq
The  map $B\to A$ induces a canonical  $A$-bimodule
map $p_{B|A}:\Om^1(B|A)\to\Om^1A.$

In the special case where  $A=B/I$, one has the following
short exact sequence
of $A$-bimodules, cf. \cite{CQ}, which is a  noncommutative  analogue of 
\eqref{conormal},
\beq{CQ}
\xymatrix{
0\ar[r]&I/I^2 \ar[rr]^<>(0.5){d_{B|A}}&&
\Om^1(B|A)\ar[rr]^<>(0.5){p_{B|A}}&& \Om^1A\ar[r]& 0.
}
\eeq
Here,  $I\sset B$ is 
 a two-sided ideal, and the  map $d_{B|A}$ is induced by  restriction to $I$ of
the de Rham differential $d: B\to\Om^1B,\, b\mto db,$ cf. \eqref{De}.

Now, fix
a potential $\Phi\in F\br,$ on $F=\C\langle\xxx\rangle$.
We have the  algebra
  $\A=\A(F,\Phi)=F/\llb\paf{i}\rrb_{i=1,\ldots,n}$, 
 the DG algebra
$\FD=\FD(F,\Phi),$ and the  
algebra projection $\FD\onto\A$, see \eqref{HD}.
Thus, we can form
the corresponding $\A$-bimodule $\Om^1(\FD|\A).$ 
\begin{defn}\label{lom} 
Define the  {\em cotangent complex} associated to $(F,\Phi)$ to be
$$
\LOM:=\Om^1(\FD|\A)\,=\,\A(F,\Phi)\,\botimes_{_{\FD(F,\Phi)}}\,\Om^1\FD(F,\Phi)
\,\botimes_{_{\FD(F,\Phi)}}\,
\A(F,\Phi).
$$
This  is a
$\operatorname{DG}$ 
$\A$-bimodule, with the grading   $\LOM=\bigoplus_{r\geq 0}
{\mathbb L}\Om^1_{r}(F,\Phi)$ and differential $\d:
{\mathbb L}\Om^1_\idot(F,\Phi)\to{\mathbb L}\Om^1_{\idot-1}(F,\Phi),$
both being  induced from those on
the DG algebra ~$\FD.$
\end{defn}

To give a more explicit  description of the  cotangent complex,
consider the composite
$$
p_{\FD|\A}:\
\xymatrix{{\mathbb L}\Om^1_0(F,\Phi)\ar@{->>}[rr]&& \Om^1\A
\ar@{^{(}->}[rr]^<>(0.5){\jmath_{_\A}}&&\A\o\A,}
$$
an   $\A$-bimodule map
induced by  the projection $\FD\onto\A$
and the tautological imbedding $\jmath_\A$, cf. ~\eqref{morph},
and also a similar map
$\dis p_{F|\A}:
\Om^1(F|\A)\onto\Om^1\A\into\A\o\A.
$

Further, it follows from \eqref{Ombasis}-\eqref{Dbasis}
  that the spaces  $\Om^1(F|\A)$ and $\dder(F|\A)$ are both free
$\A$-bimodules,
with bases  $\{dx_i\}_{i=1,\ldots,n}$ and 
$\{\fpa{i}\}_{i=1,\ldots,n},$
respectively. 
M. Van den Bergh
\cite{VdB3} introduced  an important {\em contraction 
map} $\dis\dder(F|\A)\to\Om^1(F|\A)$ defined
~by 
\beq{imap}
 \xymatrix{ a'\o \frac{\pa}{\pa x_i} \o a''\;
\ar@{|->}[r]& \; i_{a'\frac{\pa}{\pa x_i}a''} (d^2\Phi) :=
\sum\limits_{j=1}^n\,a'
\left({\Large\frac{\pa^2\Phi}{\pa x_i\pa
x_j}}\right)'\o dx_j\o
\left({\Large\frac{\pa^2\Phi}{\pa x_i\pa x_j}}\right)''a''\!.
}
\eeq

\begin{prop}\label{cotangent_thm} \vi  The $\operatorname{DG}$-module $\LOM$
is concentrated in degrees $0,1,2,$ and there is an isomorphism
between the complexes in two rows of the following diagram
\beq{diexac}
\xymatrix{0\ar[r]&
{\mathbb L}\Om^1_2(F,\Phi)\ar@{=}[d]\ar[r]^<>(0.5){\d_2}&
{\mathbb L}\Om^1_1(F,\Phi)\ar@{=}[d]\ar[rr]^<>(0.5){\d_1}&&
{\mathbb L}\Om^1_0(F,\Phi)\ar@{=}[d]\ar[r]^<>(0.5){p_{\FD|\A}}&\A\otimes
\A\ar@{=}[d]\ar[r]&0\\
0\ar[r]&\A\o \A\ar[r]^<>(0.5){\ad}&
 \dder(F|\A)\ar[rr]^<>(0.5){\eqref{imap}}&&
\Om^1(F|\A) \ar[r]^<>(0.5){p_{F|\A}}&
\A\otimes \A\ar[r]&0.}
\eeq 

\vii The nontrivial homology groups of the cotangent complex are as follows
$$
H_r\big(\LOM,\,\d\big)=
\begin{cases}
\quad\Om^1\A&\operatorname{if}\en r=0;\\
\Ext_{\bimod{A}}^1(\A,\A\o \A)&\operatorname{if}\en r=1;\\
\Hom_{\bimod{A}}(\A,\A\o \A)&\operatorname{if}\en r=2.
\end{cases}
$$

\viii Each row in \eqref{diexac} is self-dual, i.e., it goes to itself under
the duality \eqref{DUAL}.
\end{prop}

We will refer to the
complex in either row of diagram \eqref{diexac} as 
the {\em extended cotangent complex}.
The complex in the bottom row
 has been independently introduced by M. Van den
Bergh [VdB3]. 
The proof of Proposition \ref{cotangent_thm} will be given in \S\ref{cy3_pf}.

{\large\section{{\textbf{Representation functor, critical points, and 
vanishing cycles}}}\label{rep_sec}}
\subsection{Informal outline}\label{info}
Representation functor
 provides a bridge between noncommutative and the usual commutative
geometry. The main objective  of \S\ref{rep_sec} is to describe various geometric
structures which arise once one applies the
representation functor either to the DG algebra $\FD_\idot(F,\Phi)$
or to the cotangent complex $\LOM$.
It may be instructive to keep in mind
the following general setup, cf.
\S\ref{general}
for more details.

One starts with a data $(F,\al),$ where
\vskip 2pt 

\npb{$F$ is a {\em smooth} algebra, cf. Definition \ref{ww},
for instance, $F=\C\langle\xxx\rangle$;}

\npb{$\al\in (\Om^1F)\br$ is a {\em closed} cyclic 1-form, e.g. $\al=d\Phi,$ for
$\Phi\in F\br$.}
\vskip 2pt 

Associated with the data $(F,\al)$, one defines an algebra $A=\A(F,\al)$.
This is a quotient
of $F$ by an appropriate two sided ideal; in the special case where
$F=\C\langle\xxx\rangle$ and $\al=d\Phi$,
the  algebra $A$
 reduces to  $\A(F,\Phi)$, the algebra considered in the previous section.

 For any $\dd=1,2,\ldots,$ one has
the scheme  $\Rep_\dd F,$ of $\dd$-dimensional $F$-representations.
It is
a smooth manifold and the closed cyclic 1-form $\alpha$ gives rise
to  $\Tr\wh\al$, an ordinary closed 1-form  on
 $X:=\Rep_\dd F$, see \S\ref{rep}. The scheme 
 $\Rep_\dd A$ turns out to be a closed subscheme in $X$ which is
equal to the zero locus of the 1-form  $\Tr\wh\al$.

The crucial geometric feature here is that
the zero locus of a closed 1-form 
 may be identified with an intersection of two {\em Lagrangian submanifolds}
in $T^*X$, the total space of the cotangent bundle on $X$ (the first submanifold is
the graph of the  1-form and the second submanifold is the zero section
of  $T^*X$).
This applies, in particular, to the scheme
 $\Rep_\dd A$.

In the above setting, one of the goals that one would like to achieve 
 is to express various
numerical invariants of the algebra  $A=\A(F,\al)$, such as ranks
of cyclic or Hochschild homology groups of $A$, or the corresponding
Euler characteristics, in terms of certain integrals, cf. Problems
\ref{qu}-\ref{ququ} below.
The integrals in question should be similar to those
used in Witten's {\em nonabelian localization theorem} for equivariant
cohomology, ~\cite{JK}.

The problem of expressing homological invariants of an algebra $A$ 
in terms of integrals has
been already studied in \cite{EG} in the case where  the scheme $\Rep_\dd A$
was  a complete
intersection. For algebras of the form  $A=\A(F,\al)$, however,
the scheme  $\Rep_\dd A$ has {\em virtual dimension zero}.
Such a scheme can {\em never} be a complete
intersection in $\Rep_\dd F$ except for the trivial case
where it consists of isolated points. Thus, the technique 
of \cite{EG} requires a serious modification which is not
known at the moment. 

It seems very likely that the right approach to the problem is
provided by the  so-called
Batalin-Vilkovisky (BV) formalism.
In general, let $X$ be a manifold, $\phi$  a regular function on $X$,
and $L\sset T^*_\text{odd}X$ a Lagrangian submanifold in
the `odd' cotangent bundle on $X$.
 Following BV-formalism, one considers
 integrals
of the form
\beq{integ}\int_L e^{\phi\nat},
\eeq
where $\phi\nat$ is a  function on $T^*_\text{odd}X$ obtained
by a slight modification of $q^*\phi$, the pull-back
of the function $\phi$ via the bundle projection $q:T^*_\text{odd}X\to X$.

To make sense of the integral in \eqref{integ}, one needs to specify
a manifold
$L$ as well as  a volume-form
on $L$. In general, there is neither a preferred choice of
  Lagrangian submanifold $L$ nor a canonically defined volume on
it. It is known that a choice of
nowhere vanishing  volume-form
on $X$ itself provides a natural  volume-form on
any  Lagrangian
submanifold
in $T^*_\text{odd}X$. We see that in order to apply  BV-formalism,
$X$ has to be a CY manifold. Now, given a  CY manifold
$X$, there is a canonically
defined second order differential operator
 $\Delta$ on $T^*_\text{odd}X$. 
 Furthermore, one proves a version of Stokes theorem saying that the integral
in \eqref{integ} depends only on the isotopy class of
the  Lagrangian
submanifold $L\sset T^*_\text{odd}X$, provided the
 function $\phi\nat$  satisfies
the {\em quantum master equation} $\Delta \phi\nat +\{\phi\nat,\phi\nat\}=0$,
cf. \eqref{qmaster} and \eqref{master}.

Now, let $F=\C\langle\xxx\rangle$. Then, each of the manifolds
$\Rep_\dd F, \, \dd=1,2,\ldots,$ is a vector space, hence,
it  has a natural  Eucledian volume.
Any potential
 $\Phi\in F\br$ gives rise to a polynomial $\Tr\wh\Phi$
on the vector  space $\Rep_\dd F$. Further, let $\A(F,\Phi)$
be  the algebra
associated with the potential $\Phi$.
 The corresponding  representation scheme, $\Rep_\dd\A(F,\Phi)$,
may by identified with the critical set of
the  polynomial $\Tr\wh\Phi$.

Following BV-formalism, one  might expect that interesting
algebraic invariants of the  algebra $\A(F,\Phi)$
are provided by the asymptotics of integrals 
of the form, cf. \eqref{integ},
\beq{int2}
\int_{L_\dd\sset T^*_\text{odd}(\Rep_\dd F)}e^{\dd\cdot\Tr\wh\Phi\nat}\qquad\text{where}\quad\dd\to\infty.
\eeq

More generally, let $F$ be an {\em arbitrary}, not necessarily free,
{\em smooth} algebra.
In \cite{GS2}, we  analyze noncommutative structures
on the algebra $F$ that are  necessary to insure 
that each of the manifolds $\Rep_\dd F, \, \dd=1,2,\ldots,$ be a CY
manifold,
so that it makes sense to consider integrals like \eqref{int2}.
It turns out  that the
DG algebra $\FD(F,\Phi)$ plays the role of  $T^*_\text{odd}X$.
Specifically, in  \cite{GS2}, we introduce the notion of a
noncommutative  BV structure; furthermore, we show that giving
the algebra $\FD(F,\Phi)$  a noncommutative  BV structure provides
each of the manifolds $T^*_\text{odd}(\Rep_\dd F)$
with a natural BV-operator $\Delta_\dd.$

This way, using the formalism developed in \cite{GS2},
 it is not difficult to extend 
 the constructions  discussed  in the present section below
in the special case where $F=\C\langle\xxx\rangle$
to the more general setting of  an arbitrary smooth algebra
equipped with a  noncommutative  BV structure.
 
\subsection{Reminder}\label{rep}  
We write $\Mat_\dd$ for the algebra of complex $\dd\times \dd$-matrices.
Given an algebraic variety (or scheme) $X$, 
let $\oo_X$ be the structure sheaf of $X$
and put $\C[X]=\G(X, \oo_X).$ We
let  $\La^\hdot\T(X)$, resp.
 $\La^\hdot\T^*(X)$, denote 
the graded algebra of regular polyvector
fields, resp.  differential forms,  on $X$.

Let $A$ be a finitely presented $\C$-algebra.
For each integer $\dd\geq 1$, the set
$\Hom_\text{alg}(A,\Mat_\dd),$
of all algebra homomorphisms $A\to \Mat_\dd$,
has the natural structure of a (not necessarily reduced) affine
scheme of finite type over $\C$, to be denoted
$\Rep_\dd A$. 

The group $\GL_\dd$ acts on $\Mat_\dd$ by algebra automorphisms
via conjugation. 
This gives a $\gld$-action on the scheme $\Rep_\dd A$ by
base change transformations.
We write $\C[\Rep_\dd A]^\gld$ for the algebra of
$\gld$-invariant regular functions on $\Rep_\dd A$.

Any algebra homomorphism  $f:B\to A$
induces a $\gld$-equivariant  morphism of schemes
$f^*: \Rep_\dd A\to\Rep_\dd B.$ This way, we obtain a
contravariant functor 
\beq{rep_fun}
\Rep_\dd:\ \textsf{Algebras}\too\textsf{Affine $\gld$-schemes}.
\eeq

To any element $a\in A$, one associates naturally
a $\gld$-equivariant polynomial map $\wh a: \Rep_\dd A\to\Mat_\dd,\,
\rho\mto
\rho(a).$ Taking the trace of the matrix $\rho(a)\in\Mat_\dd$,
yields a $\gld$-invariant regular function $\Tr\wh a: \Rep_\dd A\to\C,
\,\rho\mto\Tr\rho(a)$. For  $a\in [A,A],$ we have
$\Tr\wh a=0$, due to symmetry of the trace.
Hence, we obtain a well defined linear map
\beq{tr}
\Tr_\dd: \ A\br=A/[A,A]\too \C[\Rep_\dd A]^\gld, \en
a\mto\Tr\wh a.
\eeq
\begin{notation}\label{OA}  Given a   {\em graded} associative 
(super) algebra
$B^\hdot$,
 the notation $B\br$ will be used to
denote the {\em super}-commutator quotient, $B/[B,B]_\text{super},$
where $[B,B]_\text{super}$ stands for the $\C$-linear span of the
elements of the form $ab- (-1)^{|a|\cdot|b|} ba,$
for all homogeneous $a,b\in B$ of degrees $|a|$ and $|b|,$
respectively. Given
a $\Z/2\Z$-graded  super-vector space $V=V^\text{even}\oplus V^\text{odd}$,
we use the notation
 $\Sym V:=(\Sym V^\text{even})\otimes (\La V^\text{odd})$.

For any  positive
integer $\dd$,  we introduce a (super)-commutative algebra
$$
\oo_\dd(B):=\Sym(B\br)\big/\llb 1_{\Sym}-\dd\cd 1_B\rrb
\;=\;\frac{\Sym(B/[B,B]_\text{super})}{\llb 1_{\Sym}-\dd\cd 1_B\rrb}.
$$
Here, we write $1_B$, resp.
 $1_{\Sym}$, for the unit element of $B$, resp. of $\Sym(B\br)$,
Thus, $\llb 1_{\Sym}-\dd\cdot 1_B\rrb$ stands for the two-sided ideal
generated by the element $1_{\Sym}-\dd\cdot 1_B\in
\Sym^0(B\br)+\Sym^1(B\br).$\eex
\end{notation}

It is clear that the linear map in \eqref{tr}
sends the unit element $1_A$ to the constant function
$\Tr(\wh1_A)=\dd$. Therefore, the map $\Tr_\dd$
can be uniquely extended, by multiplicativity, to a
({\em surjective}, cf. \cite{LP}) morphism
of commutative algebras
\beq{tralg}
\Tr_\dd: \
\oo_\dd(A)\onto\C[\Rep_\dd A]^\gld,\quad
a_1\,a_2\,\ldots\,a_r\mto
\Tr\wh a_1\cd\Tr\wh a_2\cd\ldots\cd\Tr\wh a_r.
\eeq

\begin{rem} Assume that $B$ is
a {\em graded}
finitely presented algebra, and
let $\dd\to\infty$. Then,  according to \cite{EG},  
the map \eqref{tralg} becomes 
{\em asymptoticalyy bijective}
in the sense of \cite{EG}, Definition 4.2.1.
Similar remark aplies also to the maps
\eqref{tr2}-\eqref{tr3} below, as well as to various other trace morphisms
introduced later in this section, cf.   Proposition ~\ref{DGqis}.
\end{rem}

\begin{defn}\label{OmDer}
One defines the  graded algebra of
noncommutative
differential forms, resp.,   polyvector fields, for an associative algebra $A$,
 as the 
tensor algebra
 $\Om^\hdot A:=T^\hdot_A(\Om^1A)$, resp., as $\Th^\hdot
 A:=T^\hdot_A(\dder A).$
\end{defn}

There is a  noncommutative   analogue of  de Rham differential, ~\cite{CQ} and also  \eqref{De},
$$ d=d_\text{DR}:\ \Om^\hdot A\to
\Om^{\hdot+1} A,\quad
a_0\,da_1\,da_2\ldots da_p\mto da_0\,da_1\,da_2\ldots da_p.
$$

For any algebra $A$,
the super-commutator quotient space $(\Th^\hdot A)\br$ has
a canonical {\em odd} Lie bracket $\{-,-\}$,
 a noncommutative analogue of the
Schouten bracket on polyvector fields, cf. \cite{VdB4}.
Moreover, to each potential $\Phi\in (\Th^\hdot A)\br$
one can canonically associate a derivation
$\bth_\Phi: \Th^\hdot A\to\Th^\hdot A$ (analogous to the vector field 
with Hamiltonian $\Phi$).

The construction of the trace map may be
also extended to other
objects associated with the algebra $A$. Specifically,
there is a natural trace map for differential forms 
\beq{tr2}
\oo_\dd(\Om^\hdot A)\to \La^\hdot \T^*(\Rep_\dd A)^\gld,\quad
a_0\,da_1\,da_2\ldots da_p\mto
\Tr(\wh a_0\,d\wh a_1\,d\wh a_2\ldots d\wh a_p).
\eeq
This map clearly commutes with the de Rham differentials $d$ on each side.

There is also a similar trace map for polyvector fields
\beq{tr3}
\oo_\dd(\Th^\hdot A)\too \La^\hdot
\T(\Rep_\dd A)^\gld,\quad\eta_1\,\eta_2\ldots\eta_p\mto
\Tr(\wh\eta_1\,\wh\eta_2\ldots \wh\eta_p).
\eeq
It intertwines the   noncommutative 
and ordinary Schouten brackets, respectively, cf.
~\cite{VdB4}.

\begin{examp}\label{Cx} For $A=\C[x]$, we have
$\Rep_\bd A=\Mat_\bd.$
Further, it is clear that $\Om^\hdot A\cong\C\langle x, dx\rangle,$
a free associative algebra on two generators, $x$ and
$dx=d_\text{DR}(x),$ cf. \eqref{De},
of degrees $0$ and $1$ respectively. Similarly,
we have $\Th^\hdot A\cong \C\langle x, \frac{\pa}{\pa x}\rangle,$
where $\frac{\pa}{\pa x}$ is as in \eqref{Dbasis}.
\end{examp}

\begin{notation}
We write $A*B$ for a free product (over $\C$) of two  associative algebras.
\hfill$\lozenge$
\end{notation}

Free product plays the role of coproduct in the category of
noncommutative associative algebras. This is reflected, for instance,
in the fact that,
for any pair of algebras, $A,B$, we have 
\beq{reprep}
\Rep_\bd(A*B)=\Rep_\bd A\times \Rep_\bd B.
\eeq
There are also canonical graded algebra isomorphisms
\beq{freeprod}
\Om^\hdot(A*B)\cong(\Om^\hdot A)*(\Om^\hdot B),\quad\text{resp.},\quad
\Th^\hdot(A*B)\cong (\Th^\hdot A)*(\Th^\hdot B).
\eeq

\begin{examp}[{\textbf{Special case:}} $F=\C\langle\xxx\rangle$]\label{Cxx}
Clearly, we have
$\C\langle\xxx\rangle
\cong\C[x_1]*\C[x_2]*\ldots*\C[x_n]$,
a free product of $n$ copies of the  polynomial algebra. Therefore, 
$\Rep_\dd F\cong\Mat_\dd\times\ldots\Mat_\dd$
($n$ factors). The corresponding $\gld$-action on the representation
scheme thus
becomes the  $\gld$-diagonal action on 
$\Mat_\dd\times\ldots\Mat_\dd,$ by conjugation.

In this case, using  Example \ref{Cx} and isomorphism \eqref{freeprod}, we obtain
\beq{ThF}
\Th^\hdot F=\C\langle\xxx,\th_1,\ldots,\th_n\rangle,
\quad{\textrm{where}}\quad \th_i:=\mbox{\Large$\fpa{i}$}, \en
i=1,\ldots,n.
\eeq

The corresponding  noncommutative  
Schouten bracket, usually referred
to as  {\em necklace bracket}, was introduced
first by M. Kontsevich \cite{Ko1}, cf. \cite{Gi1}, \cite{BLB1}.
In terms of the
generators from \eqref{ThF}, it is given by
the   classic formula
\beq{bracket}
\{\Phi,\Psi\}=\sum_{i=1}^n 
\left(\frac{\pa \Phi}{\pa x_i}\cd\frac{\pa \Psi}{\pa \th_i}-
\frac{\pa \Phi}{\pa \th_i}\cd\frac{\pa \Psi}{\pa x_i}\right)\;\text{mod}\;
[\Th^\hdot F,\Th^\hdot F]_{\text{super}},
\quad\forall \Phi,\Psi\in(\Th^\hdot F)\br.
\eeq

It is instructive to reinterpret \eqref{ThF} in a coordinate
free way. 
To this end, given
an $n$-dimensional $\C$-vector  space $V,$  we
let  $T_\C V=TV=\bigoplus_{r\geq 0}
T^rV$  denote its
 tensor algebra. Thus, a choice of 
basis
$\xxx$ of
$V$ provides an algebra isomorphism
$TV=\C\langle\xxx\rangle$. Then, we have
{\em canonical} graded algebra isomorphisms (independent of the choice
of bases):
\beq{TV}
\Om^\hdot(TV)=T(V\oplus V[-1]),\quad\text{resp.},\quad \Th^\hdot(TV)=T(V\oplus V^*[-1]).
\eeq
Here, the copy $V[-1]\sset \Om^1(TV)$ is a vector space with $\C$-basis
$\{dx_1,\,dx_2,\ldots,dx_n\},$ cf. \eqref{Ombasis},
resp. the copy $V^*[ -1]\sset\Th^1(TV)$ is a vector space with the dual basis
$\{\pa/\pa x_1,\,\pa/\pa x_2,\ldots,\pa/\pa x_n\},$
cf. \eqref{Dbasis}.
The isomorphisms in \eqref{TV} can either be verified directly
or can be deduced  from the special case $\dim V=1$ using \eqref{freeprod}.
\end{examp}
\subsection{Representations and potentials} 
From now on, we let  $F=\C\langle\xxx\rangle$.
The trace map \eqref{tr}
sends any potential  $\Phi\in F\br$ to a $\gld$-invariant polynomial
$$\Tr\wh\Phi\in\C[\Rep_\dd F]^\gld=\C[\Mat_\dd\times\ldots\times\Mat_\dd]^\gld.
$$ 

It is straightforward to see by going  through definitions
that the  representation scheme for the
quotient algebra $F/\llb\pa\Phi/\pa x_i\rrb_{i\in[1,n]},$
 may be identified with  $\crit(\Tr\wh\Phi)
\sset \Rep_\dd F$, the  {\em critical locus} of the
function ~$\Tr\wh\Phi$, cf. Remark \ref{locus}. 
In more detail, 
we may apply the contravariant functor
 $\Rep_\dd$, see \eqref{rep_fun}, to the tautological algebra
projection $\pr_\A: F\onto\A(F,\Phi)$. For the induced
closed  imbedding of the corresponding representation
schemes, one obtains   a diagram
\beq{crit_diag}
\xymatrix{
\Rep_\dd\A(F,\Phi)\;\ar@{=}[d]\ar@{^{(}->}[rrr]^<>(0.5){\pr_\A^*}&&&\;\Rep_\dd F\ar@{=}[d]\\
\crit(\Tr\wh\Phi)\;\ar@{^{(}->}[rrr]&&&
\;\Mat_\dd\times\ldots\times\Mat_\dd.
}
\eeq

\begin{examp}[\textbf{Yang-Mills potential}]\label{YM}
Let $\g a$ be a  Lie algebra over $\mathbb R$ with an invariant 
positive definite inner product.
Fix an orthonormal basis, $x_1,\ldots,x_n,$ of 
$\g a$.
Write $[x_i,x_j]_{\g a}=\sum_k c_{i,j,k}\cdot x_k,$ where
$[-,-]_{\g a}: \g a\times \g a\to \g a$ denotes the Lie bracket
and where
$c_{i,j,k}\in\mathbb R$ is a totally skew-symmetric tensor
of structure constants of the Lie algebra.

Associated with our data, on  the (complexified) free tensor algebra $T\g a
\cong\C\langle  x_1,\ldots,x_n\rangle$,
we introduce  two potentials. The first potential is
\beq{PSI}
\Psi= -\frac{1}{4} \sum_{1\leq i,j\leq n}\,(x_ix_j-x_jx_i)^2
+\frac{1}{3}\sum_{1\leq p,q,r\leq n}\,c_{p,q,r}\cdot
x_px_qx_r\in (T\g a)\br.
\eeq

In the special case of an {\em abelian} Lie algebra $\g a$,
the second sum in \eqref{PSI} vanishes, and
the resulting algebra $\A(F,\Psi)$, called {\em Yang-Mills algebra},
has been studied by A. Connes and  M. Dubois-Violette
\cite{CDV1}. This algebra
is related
to the Yang-Mills theory for connections with constant coefficients
in a trivial bundle on $\mathbb{R}^n$ (or on a torus),
cf. also \cite{CDV2} and \cite{MS}.

Assume next that $\g a$ is a {\em simple} Lie algebra,
and introduce a second potential
\beq{PHI}
\Phi= -\frac{1}{4} \sum_{1\leq i,j\leq n}\,
\big(x_ix_j-x_jx_i-[x_i,x_j]_{\g a}\big)^2
=\Psi+\lambda\cdot\sum_{1\leq i\leq n} x_i^2\in (T\g a)\br.
\eeq
(One can prove that the equality on the right holds, with  $\la>0$ 
being some positive
constant that depends on the Lie algebra $\g a.$)

Let $\U\g a$,
resp.  $\Sym\g a$, be the (complexified)
enveloping, resp.   symmetric, algebra of $\g a$.
Giving a  representation of the algebra $\Sym \g a$
amounts to giving an $n$-tuple of pairwise commuting
$\dd\times\dd$-matrices. 
Representations of the algebra $\U\g a$
may (and will) be identified with Lie algebra homomorphisms
$\g a\to \Mat_\dd.$
 
There are  natural algebra projections
\beq{US}
\pr_\A:\ T\g a\onto\A,\qquad
\pr_{\U}:\
T\g a\onto\U\g a,
\qquad\pr_{\Sym}:\
T\g a\onto\Sym\g a,
\eeq
where $\A$ stands for either $\A(T\g a,\Psi)$ or  $\A(T\g a,\Phi).$
These algebra  projections induce closed imbeddings of the corresponding
representation schemes.

Further, let $\Rep^\star_\dd T\g a\sset \Rep_\dd T\g a$ denote the real vector
subspace of $\star$-representations of the algebra $T\g a$, that is,
of homomorphisms $\rho:T\g a\to\Mat_\dd=\Mat_\dd(\C)$
such that  $\rho(\g a)\sset
\Mat^{\operatorname{skew}}_\dd$ (= space of  {\em skew-hermitian}
 $\dd\times\dd$-matrices).
Explicitly, we have 
$$\pr_\U^*(\Rep^\star_\dd\U\g a)
\sset\Rep^\star_\dd T\g a=\Mat^{\operatorname{skew}}_\dd
\times\ldots\times\Mat^{\operatorname{skew}}_\dd
\sset\Mat_\dd\times\ldots\times\Mat_\dd.
$$

Restricting potentials $\Phi,\Psi$ to the subspace of $\star$-representations
gives a pair of ${\mathbb{R}}$-valued  polynomials
$\phi_\dd=\Tr\wh\Phi,\,\psi_\dd=\Tr\wh\Psi:\ \Rep^\star_\dd T\g a\to{\mathbb{R}}.$
It is clear from \eqref{PHI} that the function $\phi_\dd$  takes nonnegative values
and that the subset
$\pr_\U^*(\Rep^\star_\dd\U\g a)$ is precisely the
set of points of
{\em absolute minima} of the function $\phi_\dd$.

The polynomial $\psi_\dd$ does not necessarily have nonnegative values.
We let $\Rep^\star_{\operatorname{min}}(\psi_\dd)\sset\Rep^\star_\dd T\g a$
be the
set (possibly empty) of points of
  absolute minima of the function $\psi_\dd$.

Write $[\![a,b]\!]:=ab-ba$ for the commutator in the associative
algebra $T\g a,$ not to be confused with
the Lie bracket $[-,-]_{\g a}$ in $\g a$. With this notation,
we find
\beq{ymrel}
\frac{\pa \Psi}{\pa x_j}=\sum_i\Big[\!\!\Big[x_i,\;
\big([\![x_i,x_j]\!]-\sum\nolimits_k\,c_{i,j,k}\cdot x_k\big)\Big]\!\!\Big],\quad j=1,\ldots,n.
\eeq

\begin{conj} \vi For any simple Lie algebra $\g a$ and any $\dd\geq 1,$
the  function $\psi_\dd$ is  bounded from below;

\vii Assume, in addition, that the dimension $\dd$ is such that
$\Irrep_\dd^\star\U\g a,$ the set 
of  {\em irreducible} Lie algebra representation
$\g a \to \Mat^{\operatorname{skew}}_\dd,$ is nonempty.
Then, one has
$$
\Rep^\star_{\operatorname{min}}(\psi_\dd)
=\pr_\U^*(\Irrep_\dd^\star\U\g a).
$$
\end{conj}
\end{examp}

\subsection{Cotangent complex and symmetric obstruction
theory}\label{obs_sec}
Let $\phi$ be a regular function on a smooth variety, or a smooth
stack $\scr X$. 
The substack  $\crit(\phi)$ of critical points of 
$\phi$ has
{\em virtual dimension zero}. Moreover,
according to \cite{Be},  such a  substack 
automatically comes with a {\em symmetric perfect obstruction theory} provided
by the length two complex
\beq{obs}\obs(\phi)=\big[\xymatrix{\T_{\scr X}\,\ar[r]^<>(0.5){d\ccirc i_{d\phi}}&
\,\T^*_{\scr X}}\big].
\eeq
Here, $d$ stands for the de Rham differential,
$i_{d\phi}: \T_{\scr X}\to\oo_{\scr X},\,\xi\mto \xi(\phi)$
is the contraction map, and we write $(-)^*$ for the dual
of  a vector space or of a sheaf of  vector spaces.

Let $G$ be a linear algebraic group with  Lie algebra $\g g$,
and $X$ a $G$-variety. Let $\g g_X:=\oo_X\o\g g$ be a
free
$\oo_X$-sheaf.
Any element $x\in {\g g}$ gives rise to a vector field $\arr{x}$ on
$X$, and the assignment $\phi\o x\mto\phi\cdot\arr{x}$ gives
a  
morphism  of sheaves, $\act: \g g_X \to\T_X.$
Further, equip the sheaf
$\g g_X$ with a
(non $\oo_X$-linear) Lie bracket defined  by the formula
$$[\phi\o x, \psi\o y]:=(\phi\cdot\psi)\o [x,y] 
-(\phi\cdot\arr{x} (\psi))\o y+ 
(\psi\cdot\arr{y} (\phi))\o x,\quad\forall x,y\in\g g,\,\phi,\psi\in\oo_X.
$$

The morphism $\act: \g g_X \to\T_X$ is compatible with Lie brackets;
this way,
$\g g_X$ acquires a structure of {\em Lie algebroid} on $X$, with
$\act$
being the anchor map.

Next, we consider the quotient  stack
$X/G$.  Given a  $G$-equivariant sheaf ${\mathcal F}$
on $X$, we may  view  ${\mathcal F}^G$, the subsheaf of
$G$-invariant sections, as a sheaf on $X/G$.
Further, if  $X$ is smooth then $X/G$ is a {\em smooth} stack with tangent, resp. cotangent,
bundle being
represented by the following length two complex
$$
\T_{X/G}=\big[\xymatrix{(\g g_X)^G\ar[r]^<>(0.5){\act}&
(\T_X)^G}\big]\quad\text{resp.},\quad
\T^*_{X/G}=\big[\xymatrix{(\T^*_X)^G\ar[r]^<>(0.5){\act^*}&
(\g g_X^*)^G}\big].
$$

Any $G$-invariant regular  function $\phi$ on $X$
clearly gives rise to a function, $\phi_G $, on the stack  $X/G$.
The above formulas for the tangent and cotangent
bundles on ${\scr X}=X/G$
show that the symmetric obstruction theory
\eqref{obs} associated with the function $\phi_G $ takes the form
of a 4-term complex
$$
\obs(\phi_G )=\big[\xymatrix{(\g g_X)^G\ar[r]^<>(0.5){\act}&
(\T_X)^G\ar[rr]^<>(0.5){_{d\circ i_{d\phi}}}&&(\T^*_X)^G\ar[r]^<>(0.5){\act^*}&
(\g g_X^*)^G}\big].
$$

Now, let $F=\C\langle\xxx\rangle$ and let $\Phi\in F\br$ be
a potential.  We fix  a positive integer $\dd$
and take $X=\Rep_\dd F$ and $G=\gld$. So $\g g=\Lie\gld$, and
we write $\grep F)^\gld:=\big(\C[\Rep_\dd F]\o\g g\big)^\gld$ for
the space of  global sections of the
corresponding sheaf $(\g g_X)^\gld.$

The stack $\Rep_\dd F/\gld$ is  smooth, hence
we may apply  the previous discussion to the
 $\gld$-invariant polynomial function $\phi:=\Tr\wh\Phi$.
 We conclude
that the corresponding 
symmetric obstruction theory, at the level of global sections, reads
$$
\obs(\Tr\wh\Phi_G )=\big[\!\xymatrix{\g g(\Rep F)^G\ar[r]^<>(0.5){\act}&
\T(\Rep F)^G\ar[r]^<>(0.5){_{d\circ i}}&\T^*(\Rep F)^G
\ar[r]^<>(0.5){\act^*}&
\g g^*(\Rep F)^G}\!\big],
$$
where we have used simplified notation $\Rep F=\Rep_\dd F,\,
G=\gld,$
and $\dis i=i_{d(\Tr\wh\Phi)}$.

The following result explains the significance of the
(4-term)  extended cotangent complex 
~\eqref{diexac} associated with the data $(F,\Phi)$.

\begin{prop} The trace maps \eqref{tr}-\eqref{tr3} induce
a natural morphism between the following 4-terms complexes
$$\Tr_\dd: \
\big[\xymatrix{\LOM \ar[r]^<>(0.5){p_{\FD|\A}}& \A(F,\Phi)\o\A(F,\Phi)}\big]
 \en\too\en\obs(\Tr\wh\Phi_G).
$$
\end{prop}
\subsection{BRST construction reminded} The following 
well known construction
of symplectic geometry serves as a motivation for various
constructions of the present paper.

Let $\aa=\bigoplus_{r\in\Z}\aa_r$ be a graded
Poisson super-algebra with Poisson bracket
$\{-,-\}: \aa_i\times \aa_j\to \aa_{i+j-n},$
where $n$ is some fixed integer.
Thus, $\aa_n\sset\aa$ is a sub Lie super-algebra.
Any element $h\in \aa$ gives rise to a
(super) derivation $\bth_h: \aa\to\aa,\,a\mto \{h,a\}$.
 
Let  $\g g$ be a finite dimensional  Lie algebra and
$\bt: \g g\to \aa_n$  a Lie  algebra homomorphism.
The map $\g g\to \Der\aa,\,x\mto \bth_{\bt(x)}$ gives a $\g g$-action
on $\aa$ by derivations. The  ideal $\aa\cdot\bt(\g g)\sset\aa$ 
is $\g g$-stable. Moreover,
the space $(\aa/\aa\cdot\bt(\g g))^{\g g}$ inherits
from $\aa$
a  Poisson bracket and a grading
$(\aa/\aa\cdot\bt(\g g))^{\g g}=\bigoplus_{r\in\Z}(\aa/\aa\cdot\bt(\g g))^{\g g}_r$.
The resulting graded Poisson super-algebra is
called {\em Hamiltonian reduction} of $\aa$.

\begin{examp}[{\sc{Moment maps and Hamiltonian reduction}}]\label{ham_examp}
The main sourse of examples of Hamiltonian reduction
comes
from Hamiltonian actions of Lie groups on
symplectic manifolds. Specifically, let $\mm$ be a smooth 
symplectic algebraic manifold, and let 
$G$ be  a connected algebraic group
with Lie algebra $\g g$.
Given a Hamiltonian $G$-action on $M$
with moment map $f: M\to \g g^*$,
 one is interested in  $f\inv(0)/G$, the
Hamiltonian reduction of $\mm$.
By definition, this is an affine scheme with coordinate
ring $\C[f\inv(0)]^G$, the  algebra of  $G$-invariant
polynomial functions on the zero fiber of the moment map.

In this case, we may put $\aa:=\C[M]$ and let
the map $\bt: \g g\to \aa$
be the restriction to
$\g g=(\g g^*)^*\sset \C[\g g^*]$ 
of the pull-back morphism $f^*: \C[\g g^*]\to \C[M]$.
If, moreover, the group $G$ is connected,
then, with the above notation, one has
$\C[f\inv(0)]^G=\big(\aa/\aa\cdot\bt(\g g)\big)^{\g g}$.\hfill$\lozenge$
\end{examp}

In general, given a graded Poisson super-algebra and
a homomorphism $\bt: \g g \to \aa_n$, as above,
one introduces a graded algebra
\beq{degrees}
\brst_\idot(\aa,\g g,\bt):=\La(\g g^*[1])\o \aa\o\La(\g g[-n-1])=\aa\o 
\La(\g g^*[1]\oplus\g g[-n-1]).
\eeq
The canonical
symplectic structure on $\g g^* \times\g g,$ combined with
the  Poisson structure on $\aa$, make
 $\brst_\idot(\aa,\g g,\bt) $ a graded Poisson
 super-algebra.

The BRST
construction produces a homogeneous element $ h\in\brst_{n-1}(\aa,\g g,\bt) ,$
such
that $\{ h, h\}=0$. Therefore, the derivation
$\bth_ h: \brst_\idot(\aa,\g g,\bt) \to \brst_{\idot-1}(\aa,\g g,\bt)$ gives a differential, called
the {\em BRST differential}. With this differential, one has a Poisson algebra isomorphism 
\beq{brst_iso}
 H_0\big(\brst_\idot(\aa,\g g,\bt),\, \bth_
 h\big)=\big(\aa/\aa\cdot\bt(\g g)\big)^{\g g}_0.
\eeq
\vskip 2pt

The element $ h$ is defined as follows.
The map $\bt$ may be viewed as an element
$\bt\in \g g^*\o \aa\sset \La\g g^*\o \aa $.
Further, the Lie bracket on $\g g$ 
may be viewed as an element
$\bmu\in (\La^2\g g^*)\o\g g\sset \La\g g^*\o\La\g g.$
Therefore, both $\bt$ and $\bmu$ may be identified
naturally with elements of
the algebra $\brst_\idot(\aa,\g g,\bt) =\La\g g^*\o \aa \o\La\g g.$

It is easy to verify that 
 the 3 requirements that: 1) the derivation $\bth_h$ be of degree
$-1$; 2) the vector spaces $\g g, \g g^*\sset \brst_\idot(\aa,\g g,\bt)$ be
homogeneous; and 3)  $\deg\bt=\deg \bmu,$ force
the grade
degree assignment that we've made in
\eqref{degrees}.

With this degree assignment,
one puts $ h:=\bt+\bmu\in\brst_{n-1}(\aa,\g g,\bt)$. The fact that $\bt: \g g\to\aa$ is
a Lie algebra morphism, combined with
Jacobi identity, insure that $\{ h, h\}=0.$

To prove \eqref{brst_iso}, let $\KK_\idot(\aa,\bt)=\aa\o\La(\g g[-n-1])$.
This graded algebra  comes equipped with a
standard
Koszul differential $\ka:\aa^p\o\La^q(\g g[-n-1])\to \aa^{p+n+1}\o
\La^{q-1}(\g g[-n-1])$,
given, for any $a\in\aa$ and $x_1,\ldots,x_q\in{\g g}$, by the formula
\beq{kos_diff}
 a\o x_1\wedge\ldots\wedge x_q
\mto \sum_{j=1}^q(-1)^{(n+1)j}\cd \bt(x_j)\cd
 a\o x_1\wedge\ldots\wedge x_{j-1}\wedge \wh x_j\wedge x_{j+1}\wedge
\ldots\wedge x_q.
\eeq
It is clear that
$H_0(\KK_\idot(\aa,\bt),\,\ka)=\aa/\aa\cdot\bt(\g g).$
Isomorphism  \eqref{brst_iso} can now be derived from a natural DG algebra 
isomorphism
\beq{brst_X}
\Big(\brst_\idot(\aa,\g g,\bt),\, \bth_
 h\Big)\cong \Big(C^\hdot(\g g,\,\KK_\idot(\aa,\bt)),\,
 \ka+d_\text{Lie}\Big).
\eeq
Here, the RHS  stands for the standard Chevalley-Eilenberg 
 complex $(C^\hdot(\g g,-),\,d_\text{Lie}),$ of
the Lie algebra $\g g$
with coefficients in $\KK_\idot(\aa,\bt)$, viewed
as a DG $\g g$-module.

\subsection{Batalin-Vilkovisky complex}\label{BVsec} The construction
presented below is our interpretation of
the so called Batalin-Vilkovisky (BV) approach to gauge theory.
The reader should be careful not to confuse between the notion
of  Batalin-Vilkovisky {\em complex} discussed in this section and that
of   Batalin-Vilkovisky {\em algebra}, a
different 
 notion that will be used in \S\ref{max} and \S\ref{BVHOCH}, for instance.

Let $G$ be an  algebraic group  with
Lie algebra $\g g$.
The $\ad\g g$-action on $\g g$ induces one
on the symmetric algebra  $\sym {\g g}$.
Given  a smooth affine $G$-variety $X$, 
we consider a super-commutative algebra
\beq{bvcomplex}
\big(\La\mathcal{T}(X)\o\sym {\g g}\big)^{\g g},
\eeq
 of $\g g$-diagonal invariants.
We introduce a differential
$$\ka:\
\La^p\mathcal{T}(X)\o\sym^q{\g g}\too
\La^{p+1}\mathcal{T}(X)\o\sym^{q-1}{\g g}.
$$

To this end, let $\aa:=\La^\hdot\mathcal{T}(X),$
and view this algebra as a graded algebra such that
$\La^q\mathcal{T}(X)$ is placed in degree $q$.
The differential $\ka$ is now defined by  formula \eqref{kos_diff},
where $n=1$ and where
 $\bt$ stands for the map $\bt: \g
 g\to\T(X)=\La^1\mathcal{T}(X)\into\aa,
x\mto\arr{x},$  induced by the
$G$-action on $X$. Note that since $n+1=2$,
 all terms in the sum in the RHS of \eqref{kos_diff} come
with a plus sign.

Further, given a  regular function
$\phi\in\C[X]$, we
 have a map $i_{d\phi}: \La^\hdot\T(X)$ $\to
\La^{\hdot-1}\T(X),$ the
 contraction  of polyvector fields on $X$ with the 1-form  $d\phi$.
It is immediate to see that, for a $G$-invariant
function, 
the differentials $\ka$  and $i_{d\phi}\o \Id_{\Sym}$ anti-commute.
Below,
we will use simplified notation
$i_{d\phi}=i_{d\phi}\o \Id_{\Sym}$,
and write $\KK_\idot(X,\g g,\phi)$ for
the vector space \eqref{bvcomplex} equipped with
 the differential $\ka+i_{d\phi}$.

Next, we  define a grading
$$\KK_\idot(X,\g g,\phi)=\boplus_{r\geq0}\KK_r(X,\g g,\phi)
,\qquad
\KK_r(X,\g g,\phi)=\boplus_{p+2q=r} \big(\La^p\mathcal{T}(X)\o\sym^q{\g g}\big)^{\g g}.
$$
With this grading, the differential $\ka+i_{d\phi}$ has  degree $-1$.
The resulting DG algebra
$\big(\KK_\idot(X,\g g, \phi),\,\ka+i_{d\phi}\big)$ is called
the {\em classical} Batalin-Vilkovisky complex.

\begin{prop} There is a natural algebra isomorphism
$$
H_0\big(\KK_\idot(X,\g g, \phi),\,\ka+i_{d\phi}\big)\cong
\C[\crit(\phi)]^{\g g}.
$$
\end{prop}

The  complex  $\KK_\idot(X,\g g,\phi)$ is related to the
BRST construction. To explain this, observe
that the $G$-action on $X$ induces
a natural Hamiltonian $G$-action on $M=T^*X,$
the total space of the cotangent bundle on $X$
equipped with the canonical symplectic structure.
Thus, we are in the setting  of Example \ref{ham_examp}.

Recall that we have introduced earlier the graded algebra
 $\aa:=\La^\hdot\T(X)$.
This algebra may be
viewed as $\C[T_\text{odd}^*X]$,
an {\em odd} version of the coordinate ring $\C[T^*X].$ 
It comes equipped with a
natural 
 {\em odd} Poisson bracket of degree $n=1$, the Schouten
bracket. 
Therefore, according to \eqref{degrees}, 
we may form a Poisson super-algebra
\beq{bvalg}
\brst_\idot(\La^\hdot\T(X),\g g,\bt)=\La(\g g^*[1])\o \La^\hdot\T(X)\o\Sym(\g
g[-2]).
\eeq

Following the BRST construction, we
put  $h=\bt+\bmu\in \brst_0(\La^\hdot\T(X),\g g,\bt).$
Further, let $q^*\phi$ be the pull-back of $\phi\in\C[X]^{\g g}$ via the cotangent bundle projection
$q: T^*_\text{odd}X\to X.$ Clearly,  one has
$\{q^*\phi,q^*\phi\}=0$.
Also, in the Poisson algebra 
$\brst_\idot(\La^\hdot\T(X),\g g,\bt)$,
one has $\{\bmu,q^*\phi\}=0$. Finally, using
 $\g g$-invariance of
$\phi$ one finds that   $\{\bt,q^*\phi\}=0$.
Hence, $\{h,q^*\phi\}=0$. We conclude that
the element
$\phi\nat:=q^*\phi+h$ satisfies
  $\{\phi\nat,\phi\nat\}=0$.

Thus, the 
Hamiltonian vector field corresponding to $\phi\nat$
gives a differential 
\beq{bv_dif}
\bth_{\phi\nat}: \brst_\idot(\La^\hdot\T(X),\g g,\bt)
\to\brst_{\idot-1}(\La^\hdot\T(X),\g g,\bt),\quad \phi\nat=q^*\phi+\bt+\bmu.
\eeq

To establish a link with the BRST construction,
let $\llb \g g^*\rrb\sset \brst_\idot(\La^\hdot\T(X),\g g,\bt)$ denote
the ideal generated by the vector space $\g g^*$.
This ideal
is easily seen to be
$\bth_{\phi\nat}$-stable.
Further, the  augmentation $\eps: \La^\hdot(\g g^*[1])\onto\La^0(\g g^*[1])=\C$ induces 
an algebra projection
$$\eps\o\Id:\ 
\xymatrix{
\brst_\idot\big(\La^\hdot\T(X),\g g,\bt\big)\big/\llb \g
g^*[1]\rrb\ar@{->>}[rr]&&\La^\hdot\T(X)\o\Sym(\g
g[-2]),
}
$$ 
cf. \eqref{bvalg}.
With this understood, using \eqref{brst_X} one proves
the following result.

\begin{prop}\label{brst_kos}
The  projection $\eps\o\Id$
induces a DG algebra isomorphism 
\vskip 4pt

$\dis\qquad\qquad
\Big(\big[\brst_\idot\big(\La^\hdot\T(X),\g g,\bt\big)\big/\llb \g
g^*[1]\rrb\big]^{\g g},
\ \bth_{\phi\nat}
\Big)\,\iso \,\Big(\KK_\idot(X,\g g, \phi),\;\ka+i_{d\phi}\Big).
$\qed
\end{prop}

Assume next that $X$ is a CY manifold and fix a volume, i.e., 
a nowhere vanishing
regular section  $\pi\in\La^n\T(X),$ where $n=\dim X$. 
Contraction with $\pi$ provides an isomorphism
$\La^\hdot\T^*(X)\iso\La^{n-\hdot}\T(X),\,\al\mto i_\al\pi.$
Transporting the de Rham differential $d: \La^\hdot\T^*(X)\to\La^{\hdot+1}\T^*(X)$
via that isomorphism, one obtains a differential
$\Delta_\pi: \La^\hdot\T(X)\to\La^{\hdot-1}\T(X)$,  the so-called {\em BV-differential}.

It is straightforward to verify that,
for any  $G$-invariant regular function $\phi\in \C[X]^{\g g},$
the element $\phi\nat=q^*\phi+\bt+\bmu$
satisfies the so-called quantum master equation
\beq{qmaster}
\Delta_\pi \phi\nat + \{\phi\nat,\phi\nat\}=0.
\eeq

The master equation implies that one has $(\Delta_\pi+\bth_{\phi\nat})^2$ $=0;$ 
one checks similarly that $(\Delta_\pi+\ka+i_{d\phi})^2=0.$
The complex $\big(\KK_\idot(X,\g g, \phi),\;\Delta_\pi+\ka+i_{d\phi}\big)$
is called the {\em quantum} BV complex.

We will apply the above in the following situation. 
Let $F=\C\langle\xxx\rangle$, and $\Phi\in F\br$.
Given an integer $\bd\geq 1$, we
let $X:=\Rep_\bd F=\Mat_\bd\times\ldots\Mat_\bd,$
and  equip this  vector space with
the Euclidean volume $\pi_\dd$. This way, we get the following classical, resp.,  
 quantum, BV complex associated with the function
$\phi=\Tr\wh\Ph$:
\beq{qbv}
\KK_\idot(\Rep_\dd F,\g g\g l_\dd,
\Tr\wh\Phi)\en\text{with differential}\en
\bth_{\phi\nat}=\ka+i_{d\phi},
\en\text{resp.},\en
\Delta_{\pi_\dd}+\bth_{\phi\nat}.
\eeq

\subsection{Noncommutative Batalin-Vilkovisky complex}\label{ncbv}
 In the next section, we are going to show that the
DG algebra $\FD(F,\Phi)$ introduced in \eqref{Dalg} is nothing but 
 the BV-construction
of  section \ref{BVsec} transplanted to a noncommutative setting.

To this end, we first define below  a  noncommutative  analogue of the BRST DG algebra
$\brst_\idot\big(\La^\hdot\T(X),\g g,\bt\big)$.
The role of the manifold $X$ will be played
in that analogue by an arbitrary
 associative algebra
$F$, to be assumed free later on. Thus, $\La^\hdot\T(X)$ gets
replaced by $\Th^\hdot F$. The role of the Lie algebra $\g g$
will be played by an
 associative algebra
$R$.

A replacement for the action map $\bt: \g g\to \T(X)$ 
is provided by the following

\begin{defn}\label{der_act} Let $F,R$ be two associative
algebras. A linear map $\bt: R\to \dder F,{r\mto\th_r}$
is said to be an $R$-{\em action on $F$ by double derivations}
if one has
\beq{nu}
(\th_r\o \Id_F)\ccirc \th_s-
(\Id_F\o\th_s)\ccirc\th_r=(\th_{r\cdot s})^{13},\quad
\forall r,s\in R.
\eeq
\end{defn}
\noindent
Here, for any  map $\th: F \to F\o F,\, u\mto \th'(u)\o \th''(u),$
we write $\th^{13}$ for the
map $F \to$ $F\o F\o F,\, u\mto \th'(u)\o1\o \th''(u).$

From now on, we assume in addition that the algebra $R$ has
{\em finite} dimension over $\C$.
Following the strategy of \S\ref{BVsec}, we introduce a graded
associative algebra, cf. \eqref{bvalg},
\beq{BBD}\brst_\idot(F,R,\bt):=T(R^*[1])\,*\, \Th^\hdot F\,*\,T(R[-2]).
\eeq

Ignore the grading in \eqref{BBD} for a moment. Then, using
\eqref{freeprod}, \eqref{TV},
we get
\beq{BBD2}
\brst(F,R,\bt)\cong(\Th F)\,*\,T(R\oplus R^*)\cong
(\Th F)\,*\,\Th(TR)\cong\Th(F\,*\,TR).
\eeq
The algebra on the right comes equipped with the  noncommutative  Schouten bracket.
 Therefore,
 any element $H\in\brst(F,R,\bt)\br=
\Th(F\,*\,TR)\br$ gives rise to a 
 derivation $\bth_H:\brst(F,R,\bt)\to\brst(F,R,\bt).$

Next, view the map $\bt$  as an element $\bt\in R^*\o\dder
F\sset R^*\o \Th^1 F.$ Similarly, multiplication map $R\o R\to R$
determines an element $\bmu\in R\o T^2(R^*)$.
We have natural vector space imbeddings
$R^*\o \Th^1 F\into\brst_0(F,R,\bt)$,
resp. $R\o T^2(R^*)\into \brst_0(F,R,\bt)$.
Hence, we may view the elements
$\bt$ and $\bmu$ as elements of
$\brst_0(F,R,\bt)$, and form the sum $\bt+\bmu$.

Now, let $\Phi\in F\br$ be  a  potential.
Using the obvious algebra  imbeddings
$F\into \Th^\hdot F\into\brst_\idot(F,R),$ 
we identify $\Phi$ with an element of $\brst_0(F,R,\bt)\br.$
Let, cf. \eqref{bv_dif},
\beq{phinat}
\Phi\nat:=\Phi+\bt+\bmu\in \brst_0(F,R,\bt)\br,
\eeq
and write
$\bth_{\Phi\nat}: \brst_\idot(F,R,\bt)\to
\brst_{\idot-1}(F,R,\bt),$ 
for the corresponding Hamiltonian derivation.

\begin{prop}\label{chi} Fix an algebra $F$,  a finite dimensional
algebra $R$, a linear map $\bt:R\to\dder F$, and a potential $\Phi\in
F\br$. Define $\Phi\nat$ as in \eqref{phinat}. Then, we have

\vi The  equation $(\bth_{\Phi\nat})^2=0$ holds
iff the map $\bt$ gives  an $R$-action on $F$ by double derivations,
cf. Definition \ref{der_act}.

\vii For any algebra homomorphism $\chi: R\to\C$, the map
$\bt: r\mto \chi(r)\cdot\De$  gives  an $R$-action on $F$ by double derivations.
\end{prop}
\noindent
The statement of part (i) is verified by direct computation.
In part (ii), we have used the double derivation $\De\in\dder F,$ see \eqref{De}.
One checks that the equation
\beq{De_eq}
(\th\o \Id_F)\ccirc \th-
(\Id_F\o\th)\ccirc\th=\th^{13}.
\eeq
holds for the double derivation $\th=\De$. This implies part (ii) of the proposition.\qed

\begin{rem} In the special case where $F=\C\langle x_1,x_2\rangle$, a free
algebra on 2 generators, the solutions of  equation \eqref{De_eq}
have been classified by T. Schedler using MAGMA program.
He showed that, up to changes of variables, any nonzero solution 
has the form
\vskip 3pt

$\dis
\eps_1\cd x_1\frac{\pa}{\pa x_1} -\eps_2\cd  \frac{\pa}{\pa x_1}x_1
 +\eps_3\cd x_2\frac{\pa}{\pa x_2}  -\eps_4\cd  \frac{\pa}{\pa x_2}x_2,\qquad\text{where each}
\en\eps_r\en\text{is either}\en0\en\text{or}\en1.$\eex
\end{rem}

Next, let 
$\llb R^*[1]\rrb\sset\brst_\idot(F,R,\bt)$
be the two-sided ideal generated by the vector space
$R^*[1]$. This ideal is easily seen to be $\bth_{\Phi\nat}$-stable.
Thus, the quotient $\brst_\idot(F,R,\bt)/\llb R^*[1]\rrb$
inherits a DG algebra structure.

\begin{examp}[{\textbf{Special case:}} $R=\C$]\label{R=C}
In this case,  we write
$t$ for the base element of the graded vector space $R[-2]=\C[-2]$
corresponding to $1\in R$. Let $\tau\in R^*[1]$ denote the
dual base element, so $\tau(t)=1$. Thus, $\deg t=2$ and $\deg\tau=-1$.

Further,  let $\bt: \C\to \dder F$ be the map that
sends $1\mto\De$; in the notation of  Proposition \ref{chi}(ii),
this map corresponds to the homomorphism $\chi=\Id: \C\to\C$.
For these choices of $R$ and $\bt$, we get
$\bmu=t\cd\tau^2\in R\o T^2(R^*),$
resp., $\bt=\bw\cd\tau\in  (\dder F)\o R^*.$ 

Thus, for any potential $\Phi\in F\br$, we have
\beq{RC}
\brst_\idot(F,R,\bt)=(\Th^\hdot F)*\C\langle t,\tau\rangle,
\quad\text{and}\quad\Phi\nat:=\Phi+\bw\cd\tau +t\cd\tau^2.
\eeq

Further, since $\llb R^*[1]\rrb=\llb\tau\rrb$, we obtain
\beq{RCT}
\brst_\idot(F,R,\bt)/\llb R^*[1]\rrb\cong(\Th^\hdot F)*\C[t].
\eeq
\end{examp}

\subsection{{Hamiltonian interpretation of the DG algebra
$\mbox{$\FD(F,\Phi)$}$}}\label{cy_rep}
In this section, we let $F= \C\langle\xxx\rangle$ and fix a
 potential $\Phi\in F\br.$ We would like to construct
a trace morphism for $(\FD(F,\Phi),\,\d)$, the
 DG algebra  that has been introduced in
 \S\ref{dg_sec}.

Observe first that the super-derivation $\d$ clearly descends to 
a well defined differential on $\FD(F,\Phi)\br$, the
super-commutator quotient of our DG algebra. 

Next, we consider $\Sym\big(\FD(F,\Phi)\br\big)$,
and view this symetric algebra as a   super-commu-tative graded algebra with respect to the
grading induced from the one on $\FD(F,\Phi)\br,$
not to be confused with the standard grading on the
symetric algebra.
The  differential on $\FD(F,\Phi)\br$ may be further
extended (uniquely) 
to a super-derivation $\d_{\Sym},$ on the  algebra
$\Sym\big(\FD(F,\Phi)\br\big),$ such that we have $(\d_{\Sym})^2=0.$

For any  integer $\dd$, the  super-derivation $\d_{\Sym}$
clearly annihilates 
$1_{\Sym}-\dd\cdot 1_\FD$, a  degree zero element of our
graded  algebra. 
 Thus,  the quotient algebra
$\dis\oo_\dd\big(\FD(F,\Phi)\big)=
\Sym\big(\FD(F,\Phi)\br\big)/\llb1_{\Sym}-\dd\cdot 1_\FD\rrb$
inherits both the grading and the differential.
This way, we have made the algebra
$\oo_\dd\big(\FD(F,\Phi)\big)$ a super-commutative DG algebra.
Let 
$\d_\dd$ denote the corresonding differential.

The result below says that
the DG  algebra $\FD_\idot(F,\Phi)$ introduced in
\eqref{Dalg}-\eqref{d_def} plays the role of
a  noncommutative Batalin-Vilkovisky complex.
Specifically, we have defined earlier, see \eqref{qbv},
 a complex
$\big(\KK_\idot(\Rep_\dd F, \g g\g l_\dd,\Tr\wh\Phi),\;\ka+i_{d\phi}\big),$
the (classical) BV complex.

\begin{prop}\label{DGqis} Let $F=\C\langle\xxx\rangle$ and $\Phi\in F\br$. Then, we have

\vi In the setup of Example \ref{R=C},
 there is a natural
{\DG } algebra isomorphism
$$
\Big(\FD(F,\Phi),\,\d\Big)\;\cong\;
\Big(\brst_\idot(F,R,\bt)/\llb\tau\rrb,\,\bth_{\Phi\nat}\Big),
\quad\oper{where}\quad
R:=\C\en\oper{and}\en\bt: 1\mto\De.
$$

\vii The  trace morphism in \eqref{tr3}
induces a natural morphism 
of \DG algebras
$$
\Tr_\dd:\
\Big(\oo_\dd\big(\FD(F,\Phi)\big),\,\d_\dd\Big)\;\too\;
\Big(\KK_\idot(\Rep_\dd F, \g g\g l_\dd,\Tr\wh\Phi),\,\ka+i_{d\phi}\Big),
\quad\forall\bd\geq 1.
$$
\end{prop}

\noindent
We remark that part (i) above
is a noncommutative
analogue of Proposition \ref{brst_kos}.\hfill$\lozenge$

\begin{proof}[Proof of  Proposition \ref{DGqis}] For $=\C\langle\xxx\rangle,$ from \eqref{ThF} and 
\eqref{RC}, we find
\beq{FRbt}
\brst_\idot(F,R,\bt)=
\C\langle\xxx,\th_1,\ldots,\th_n,t,\tau\rangle.
\eeq
Therefore, we obtain natural graded
algebra isomorphisms, cf. \eqref{Dalg}, \eqref{RCT},
\beq{natiso}
\brst_\idot(F,R,\bt)/\llb\tau\rrb\cong(\Th^\hdot F)*\C[t]\cong
\C\langle\xxx,\th_1,\ldots,\th_n,t\rangle\cong\FD_\idot(F,\Phi),
\eeq

To compare the differentials in the algebras on the left and on the right, we recall that
 any  potential $H\in\C\langle\xxx,\th_1,\ldots,\th_n,t,\tau\rangle\br$
gives rise to a
`Hamiltonian' derivation $\bth_H\in\Der\C\langle\xxx,\th_1,\ldots,\th_n,t,\tau\rangle$.
It is given by the explicit formula, cf. \eqref{bracket},
\beq{thH}
\bth_H=\sum_{i=1}^n \left(\frac{\pa H}{\pa x_i}\pa_{\th_i}-
\frac{\pa H}{\pa \th_i}\pa_{x_i}\right)+
\left(\frac{\pa H}{\pa t}\pa_{\tau}-
\frac{\pa H}{\pa \tau}\pa_t\right).
\eeq
The notation in this  classically looking formula  are as follows.
Let $x$ be  one of the $2n+2$ generators of the
algebra \eqref{FRbt}.
Given   an element 
$u$ of that algebra,  we write 
$$u\pa_x:\
\C\langle\xxx,\th_1,\ldots,\th_n,t,\tau\rangle\to
\C\langle\xxx,\th_1,\ldots,\th_n,t,\tau\rangle
$$
for the derivation  that sends $x\mto u$, and
that annihilates the other
$2n+1$  generators of this free algebra.

Now, let $H:=\bw\cdot\tau+t\cdot\tau^2$. This is an   {\em even}
potential;
one computes $\frac{\pa H}{\pa t}=\tau^2,$ and
$\frac{\pa H}{\pa \tau}=\De.$
From this,  using  \eqref{thH} and  \cite{CBEG}, Lemma 7.1.1(ii),
we find
\beq{find}
\bth_H=\ad\tau+ \De\pa_t-\tau^2\pa_\tau.
\eeq

Further, according to \eqref{inner},
in the basis $\th_i=\fpa{i},\,i=\nn,$  one has
$
\De=\sum_{i=1}^n [x_i,\th_i].$
Thus,  the potential
$\Phi\nat$, as given  in \eqref{RC}, reads
\beq{RC2}
\Phi\nat:=\Phi+\bw\cd\tau +t\cd\tau^2 =\Phi+
\sum_{1\leq i\leq n} [x_i,\th_i]\cd\tau+t\cd\tau^2.
\eeq
The  corresponding {\em  odd} super-derivation  equals, cf. \eqref{find},
\beq{phph}
\bth_{\Phi\nat}=
\sum_{i=1}^n \frac{\pa \Phi}{\pa x_i}\pa_{\th_i}\;+\;
\ad\tau+ \sum_{i=1}^n [x_i,\th_i]\pa_t-\tau^2\pa_\tau.
\eeq

A straightforward comparison  of formula \eqref{phph}
with 
 formula  \eqref{d_def} for  the differential  $\d$ 
 yields the required   compatibility of the  differentials
in the DG algebras
$\brst_\idot(F,R,\bt)/\llb\tau\rrb$ and $\FD(F,\Phi)$, respectively, cf. \eqref{natiso}.
Part (i) of the proposition follows.

To prove part (ii),  we use the isomorphism 
$\FD(F,\Phi)\cong
(\Th^\hdot F)*\C[t]$, cf. \eqref{natiso}. Hence, in view of Example 
\ref{Cx} and \eqref{reprep}, we have
$$\Rep_\dd \FD(F,\Phi)\cong\Rep_\dd\big((\Th^\hdot F)*\C[t]\big)\cong
\Rep_\dd(\Th^\hdot F)\;\times\;\Mat_\dd.
$$
 We conclude that the trace map \eqref{tr3}
 gives rise to a natural graded algebra morphism
\beq{tr4}
\oo_\dd\big(\FD(F,\Phi)\big)=\oo_\dd\big((\Th^\hdot F)*\C[t]\big)
\too
\big(\La^\hdot\mathcal{T}(\Rep_\dd F)\o\C[\Mat_\dd]\big)^\gld.
\eeq

Now,  view
the matrix algebra $\Mat_\dd$ as a Lie algebra with respect to
commutator,
that is, identify $\Mat_\dd$ with
 the Lie algebra  ${\g g}=\Lie \gld.$ We may further identify
$\g g$ with $\g g^*$ via the trace paring. Thus,
we obtain algebra isomorphisms $\C[\Mat_\dd]\cong
\C[{\g g}]\cong\C[\g g^*]\cong\Sym\g g.$
With these  isomorphisms, the space
on the right of \eqref{tr4}
may be identified with
$\KK(\Rep_\dd F,\g g\g l_\dd,\Tr\wh\Phi)$. 

We define the
map $
\Tr_\dd$ in the statement of the proposition to be the  map \eqref{tr4}.
A direct verification shows that the map thus defined
intertwines the differential on the space
$\oo_\dd\big(\FD(F,\Phi)\big)$ induced by $\d$ with  the differential
$i_\al+\ka$.
\end{proof}

\subsection{Relation to $A_\infty$-algebras}\label{kevin}
Let $\g g$ be a finite dimensional
Lie algebra and $\bt: \g g\to \gl(V)$ a finite dimensional
representantation. In the setting of \S\ref{BVsec}, we may put
$X:=V$, a vector space with a linear $\g g$-action.
In this case, by definition, we have
\begin{align}
\brst_\idot(V, \g g,\phi)&=\La(\g g^*[1])\o \La^\hdot\T(V)\o \Sym(\g
g[-2])\label{linfty}\\
&=\Sym\Big(\g g^*[1]\oplus V^*\oplus V[-1] \oplus \g g[-2]\Big),\nonumber
\end{align}
where $\Sym(-)$ on the right is understood as the symmetric
algebra of a super-vector space.

According to \S\ref{BVsec}, any polynomial $\phi\in\C[V]$ makes the above space a
DG algebra with differential $\bth_{\phi\nat}$. 
Observe further that giving a differential
on \eqref{linfty} is the same thing as
providing the graded vector space $\g G:=\g g \oplus V[-1]\oplus V^*[-2]\oplus
\g g^*[-3]$
with a structure of $L_\infty$-algebra. This  $L_\infty$-algebra has
a nondegenerate invariant bilinear form; it has 
been studied in \cite[\S3.2]{GG} in connection with Maurer-Cartan
equations.
\smallskip

The  noncommutative BV complex introduced
in \S\ref{ncbv} has a similar interpretation. 
Specifically, fix a  finite dimensional vector space $V,$
and put $F=T(V^*)$. Let a finite dimensional associative algebra
$R$ act on $F$ by double derivations, cf. Definition
\ref{der_act}. 

From \eqref{TV} and \eqref{BBD}, we find
$$\brst_\idot(F,R,\bt)=T(R^*[1])*\Th^\hdot F*T(R[-2])=
T\big(R^*[1]\oplus V^*\oplus V[-1]\oplus R[-2]\big).
$$

Given a potential $\Phi\in F\br$, the construction
of \S\ref{ncbv} yields a DG algebra structure on  $\brst_\idot(F,R,\bt)$
with differential
$\bth_{\Phi\nat}$. The latter structure can be viewed
as a structure of $A_\infty$-algebra on the vector space
$\VV:=R \oplus V[-1]\oplus V^*[-2]\oplus
R^*[-3]$.

\begin{examp} Consider a special case where $R=\C$ and use the notation introduced
in \S\ref{cy_rep}.
Thus, we have
$\VV:=\C\cd \tau\;\boplus\; V\;\boplus\; V^*\,\boplus \,\C\cd t,$
where $\tau,\, V,\,V^*,\,t$ are placed in degrees
$0,1,2,3,$ respectively. 

For any  $\Phi\in F\br$,
the potential $\Phi\nat$ given by \eqref{RC2} 
satisfies
\begin{eqnarray}
&\textsl{master equation:}\quad&\{\Phi\nat,\Phi\nat\}=0;\label{master}\\
&\textsl{unit
axiom:}\quad&\,\bth_t(\Phi\nat)=\bw.
\label{unit}
\end{eqnarray}
Here,  equation \eqref{master} follows from  Propositions \ref{chi}--\ref{DGqis},
and equation \eqref{unit} is easily verified directly from formula \eqref{RC2}.

Equations \eqref{master}-\eqref{unit} say that the
potential $\Phi\nat$ gives the vector
space $\VV$ the structure of a {\em unital cyclic}
$A_\infty$-algebra with a {\em nondegenerate inner product},
that pairs $V$ with $V^*$ and $t$ with $\tau$.
The extra-terms in the potential $\Phi\nat$
account, in particular,  for the unit property:
$m_2(t, a)=m_2(a,t)=a,\, \forall a\in \VV$,
where $m_2: \VV\times\VV\to\VV$ is an associative product on $\VV$,
the degree 2 part
of the $A_\infty$-structure.

Following  Kevin Costello,
the resulting $A_\infty$-structure  $(\VV,\Phi\nat)$ can be described
as follows. Let $\VV_+:=V\;\boplus\; V^*\,\boplus \,\C\cd t$ be the
positive degree part of $\VV$.
 The augmentation 
$\VV\onto \VV/\VV_+\iso\C, \, \tau\to 1$
is an $A_\infty$-algebra morphism that gives the
 1-dimensional vector space $\C_\tau=\C$ the structure of an
 $A_\infty$ $\VV$-module.

K. Costello  observed that the isomorphism of
 Proposition \ref{DGqis} may be reinterpreted as the following result
\begin{prop}
There is a natural DG algebra isomorphism
$$\FD(F,\Phi)\cong \RHom_{\Lmod{\VV}}(\C_\tau,\C_\tau).
$$
\end{prop}
\begin{proof}
We may compute $\RHom$ on the right 
by means of the standard reduced Bar resolution 
$\VV\o T^\hdot(\VV_+)\qis\C_\tau$. Thus, we get
$$\RHom^\hdot_{\Lmod{\VV}}(\C_\tau,\C_\tau)=
\Hom^\hdot_{\Lmod{\VV}}(\VV\o T^\hdot(\VV_+),\,\C_\tau)
=T^\hdot(\VV_+^*[1]).
$$

Further, the inner product on $\VV$ provides a vector space  isomorphism
$\VV_+^*[1]\cong\VV/\C\cd\tau$.
Therefore, we obtain
$$
\RHom^\hdot_{\Lmod{\VV}}(\C_\tau,\C_\tau)
\cong T^\hdot(\VV/\C\cd\tau)=(T^\hdot\VV)/\llb\tau\rrb\cong
\Th^\hdot(TV\nat)/\llb\tau\rrb
\stackrel{^\text{Prop. \ref{DGqis}}}\eqq\FD(F,\Phi),
$$
where we have introduced the notation  $V\nat:=\C\tau\oplus V$.
It is easy to verify that the natural DG algebra structure
on the $\RHom$-space on the left corresponds, via this chain
of isomorphisms, to the DG algebra structure that
we have defined on $\FD(F,\Phi)$.
\end{proof}

\begin{rem}
In physics,
the  vector space $V\nat=\C\tau\oplus V$
 is to be  thought
of as the {\em chiral ring} of a topological string theory,
cf. \cite{La}. 
The element $\tau$
 corresponds
to the unit of the chiral ring.
According to   \cite{La}, a general  topological
{\em open} string theory is determined by
an {\em even} potential $\Phi\nat\in(\Th^\hdot(TV\nat))\br$,
that must satisfy  \eqref{master}-\eqref{unit}.
\end{rem}
\end{examp}

\subsection{Batalin-Vilkovisky algebra structure on $\oo_\bd\big(\FD(F,\Phi)\big)$}\label{max}
Let $\g G$ be an arbitrary super Lie {\em bi}algebra  with an odd Lie 
bracket $\{-,-\}: \Sym^2\g G\to \g G$ and with cobracket $\nu: \g G\to
\Sym^2\g G$. 

The super-commutative algebra $\Sym\g G$ acquires
the structure of a Batalin-Vilkovisky algebra. Specifically,
the Lie bracket on $\g G$ gives  $\Sym\g G$ the standard structure of
a DG  {\em Gerstenhaber algebra}, with Chevalley-Eilenberg
differential. 
Further, the Lie cobracket
$\nu$ gives rise to a BV-operator
$\Delta:\Sym\g G\to\Sym\g G.$ In more detail,
one defines an operator
\begin{align}
&&\Delta:\ \Sym^n\g G \too
&\Sym^{n+1}\g G \oplus
\Sym^{n-1}\g G,\quad\forall n>0,\nonumber\\
&&\hphantom{x}{}_{}\;\ {\hphantom{x}}\;
 u_1\cd \ldots \cd  u_n\longmapsto&\sum_{1\leq i\leq n} 
(-1)^{i-1} \nu(u_i)\cd  u_1 \cd \ldots \cd  \wh
u_i\cd \ldots \cd  u_n\label{bvoper}\\
&&&\hphantom{x}\quad+
\sum_{1\leq k<\ell\leq n} (-1)^{k+\ell-1} \{u_k,u_\ell\} \cd u_1\cd \ldots \cd \wh u_k\cd
\ldots \cd 
\wh u_\ell\cd \ldots \cd  u_n.\nonumber
\end{align}

Also, for $n=0$, we put $\Delta(1):=0.$
The compatibility of the Lie bracket and cobracket on a Lie bialgebra
is easily seen to be equivalent to
the equation $\Delta^2=0.$ 

We apply the general construction above to prove the following result
(cf. \cite{GS2} for a generalization and additional details).

\begin{prop}\label{bv11} Let $F=\C\langle\xxx\rangle.$ 
For any integer $\bd\geq 1$, the graded algebra
$\oo_\bd\big(\FD(F,\Phi)\big)$ has a natural structure of \DG
BV-algebra,
with BV-operator $\Delta_\bd$.
\end{prop}
\begin{proof}[Sketch of Proof] First, recall   the algebra isomorphism
$\FD(F,\Phi)\cong
(\Th^\hdot F)*\C[t]$. Thus, we may (and will) identify
the super-commutator quotient $\FD(F,\Phi)\br$
 with $\big((\Th^\hdot F)*\C[t]\big)\br$.

Now, one  uses
 the noncommutative Schouten double bracket
on $\Th^\hdot F$, introduced by Van den Bergh,
to define an odd Lie bracket $\{-,-\}$ on the vector space
$\big((\Th^\hdot F)*\C[t]\big)\br$, hence, on
 $\FD(F,\Phi)\br$.
 This way, the super-commutator quotient
$\FD(F,\Phi)\br$ acquires a natural structure of DG Lie super-algebra,
with  Lie bracket $\{-,-\}$ and differential $\d$. 

Further, the construction due to T. Schedler
\cite{Sc} gives an odd Lie {\em co}bracket 
\beq{cobra}
\nu:\ \FD(F,\Phi)\br\to \Sym^2\big(\FD(F,\Phi)\br\big),
\quad
f\mto \pr\left(\sum_{i=1}^n \frac{\pa^2 f}{\pa x_i\pa \th_i}-
\frac{\pa^2 f}{\pa \th_i\pa x_i}\right).
\eeq
In this formula, we have used two maps
 $\frac{\pa^2 f}{\pa x_i\pa \th_i}, \,
\frac{\pa^2 f}{\pa \th_i\pa x_i}:\
\FD(F,\Phi)\br \to \FD(F,\Phi)\o \FD(F,\Phi)$,
and also the map `$\pr$' given by the following
composition of natural projections
$$\pr:\
\FD(F,\Phi)\o \FD(F,\Phi)\onto
\FD(F,\Phi)\br\o \FD(F,\Phi)\br\onto
\Sym^2\big(\FD(F,\Phi)\br\big).
$$

It has been proven in \cite{Sc} that the bracket $\{-,-\}$ and
cobracket $\nu$
make  $\FD(F,\Phi)\br$ a Lie super bi-algebra.
Therefore, by the general construction explained earlier, 
the graded algebra  $\Sym\big(\FD(F,\Phi)\br\big)$
acquires the structure of a  BV-algebra, with  BV-operator
$
\Delta:\ \Sym\big(\FD(F,\Phi)\br\big)\too
\Sym\big(\FD(F,\Phi)\br\big)
$, cf. \eqref{bvoper}.
It is immediate from definitions that, for any $\bd\in \C$,
we have $\{1_{\Sym}-\bd\cdot  1_\FD,-\}=0$ and
$
\Delta(1_{\Sym}-\bd\cdot 1_\FD)=0.$ Thus, the
BV-operator
(as well as  the bracket)
descends to an operator $\Delta_\bd$ on $\oo_\bd\big(\FD(F,\Phi)\big)$.

Finally, on  $\oo_\bd\big(\FD(F,\Phi)\big)$,
we also have the differential $\d_\bd$. It
 is easy to verify that $\d_\bd\ccirc\Delta_\bd+\Delta_\bd\ccirc\d_\bd=0$.
We conclude that $\big(\oo_\bd\big(\FD(F,\Phi)\big),\,\Delta_\bd,\d_\bd\big)$
is a DG BV-algebra.

The reader is referred to \cite{GS2} for a much more elaborate construction.
\end{proof}

Let $G$ be a connected algebraic group
with Lie algebra $\g g$ and $X$  a smooth $G$-variety.
Recall that, associated with a regular
$G$-invariant function $\phi$ on $X$,
 is a (perverse) sheaf $\van(\phi)$, of vanishing cycles.
This is a $G$-equivariant constructible complex
whose cohomology sheaves are known to be
 set-theoretically supported on 
 $\crit(\phi)$. Let
$\chi^G(\van(\phi))\in \C[\g g]^{\Ad G}$ denote the
corresponding $G$-equivariant Euler characteristic.

Now,  fix $\bd\geq 1$, 
let $X:=\Rep_\bd F$
and $G=\GL_\bd$. Given a potential $\Phi\in F\br,$ we
may apply the above to the $\GL_\bd$-invariant polynomial
$\phi=\Tr\wh\Phi$ to get a $\GL_\bd$-equivariant
 sheaf $\van(\phi)$. Part (i) of the following result is proved by a straightforward
computation.

\begin{prop} Let $F=\C\langle\xxx\rangle$ and $\Phi\in F\br$. Fix
$\bd\geq 1$
and put $\phi=\Tr\wh\Phi$. 

\vi The map of Proposition \ref{DGqis}(ii) induces a morphism of
complexes, cf. \eqref{qbv},
$$\Tr_\bd:\
\Big(\oo_\bd\big(\FD(F,\Phi)\big),\,\Delta_\dd+\d_\dd\Big)
\;\too\;
\Big(\KK(\Rep_\dd F,\g g\g l_\bd,\phi),\, \Delta_{\Rep_\dd}+i_{d\phi} +
\ka\Big);
$$

\vii 
In $\C[\Lie \gld]^{\Ad\GL_\bd}$,
one has
$$
\CHI
\big(\KK(\Rep_\dd F,\g g\g l_\bd,\phi),\, \Delta_{\Rep_\dd}+i_{d\phi} +
\ka\big)=\chi^{\GL_\bd}(\van(\phi)).
$$
\end{prop}

\begin{proof}[Idea of proof of part (ii)] 
First, we use the Euclidean volume  on $\Rep_\dd F$
to get an isomorphism $\La^\hdot\T(\Rep_\dd F)$ $\cong
\La^{r-\hdot}\T^*(\Rep_\dd F),$ where $r=\dim\Rep_\dd F.$
Contraction with the 1-form ${d\phi}$ goes, via this  isomorphism,
to the wedging $\beta\mto\beta\wedge{d\phi}.$ 
We simplify the notation and write $\wedge_{d\phi}$,
resp. $i_{d\phi}$, for the tensor product of
the latter operation with the identity map
on $\Sym\g g,$ resp. for $i_{d\phi}\o \Id_{\Sym}$. This way,
one easily obtains an identification of
complexes
$$
\Big(\KK(\Rep_\dd F,\g g\g l_\bd,\phi),\, \Delta_{\Rep_\dd}+i_{d\phi}+
\ka\Big)\cong
\Big(\big(\Om^\hdot(\Rep_\dd F)\o\C[\gld]\big)^{\GL_\bd},\, d+\wedge_{d\phi}+
\ka^\dag\Big),
$$
where the terms $d+\ka^\dag$ account for the standard equivariant
differential on the de Rham complex for $\gld$-equivariant cohomology.

Finally, one proves using e.g. \cite{Kap}, that the  complex
with differential $(d+\ka^\dag)+\wedge_{d\phi}$
is exactly the de Rham complex for computing $\GL_\bd$-equivariant
cohomology of the vanishing cycles sheaf $\van(\phi).$
\end{proof}

\begin{rem}  Recall that whenever
a closed subscheme $ Y$ of a smooth scheme
$ X$ has a symmetric perfect obstruction theory,
there is a Borel-Moore homology class
$[ Y]\vir \in H^\text{BM}_0( X),$  the {\em virtual fundamental class},
cf. \cite{Be}. In particular, the  critical locus of a
regular function $\phi\in\C[X]$
has a well-defined  class
$[\crit(\phi)]\vir \in H^\text{BM}_0(X)$.

More generally, let $G$ be a reductive group, let
$X$ be a smooth $G$-variety and  $\phi$  a  $G$-invariant regular
 function on $X$. Then, there is  a well defined
 virtual fundamental class $\crit(\phi)$ in the $G$-equivariant
Borel-Moore homology of $X$.
By the main result of \cite{Be}, \cite{BF},
one has 
\beq{BF}
\int^G_{[\crit(\phi)]\vir}1=\chi^G(\van(\phi)).
\eeq
\end{rem}

\begin{problem}\label{qu} Let a reductive group $G$ act linearly
on a vector space $E$
and let $\phi\in\C[E]^G$ be a $G$-invariant polynomial.
Formulate and prove some sort of Gauss-Bonnet type formula
$$
\chi^G(\van(\phi))= \int_{E}^G e^\phi.
$$
\end{problem}
\begin{problem}\label{ququ} Let $F=\C\langle\xxx\rangle$ and $\Phi\in F\br.$
Express the Euler characteristic of the complex
$\Big(\oo_\dd\big(\FD(F,\Phi)\big),\,\Delta+\d\Big)$
in terms of the
asymptotics of the  `{\em partition function}'
$$Z(\Phi,q)=\sum_{\dd=1}^\infty q^\dd\cdot\chi^{\GL_\bd}(\van(\Tr\wh\Phi)).
$$
\end{problem}

{\large\section{{\textbf{Calabi-Yau condition}}}\label{cy_sec}}
\subsection{Mirror symmetry motivation}\label{motiv}
One of the motivations for studying Calabi-Yau
geometry comes from mirror symmetry. The latter
deals with families of CY varieties in  the vicinity of 
certain  special values of parameters
where the varieties in the family
degenerate to acquire singularities, cf.  \cite{KS1}.

Specifically, let  $(X,X')$ be a pair of  mirror dual CY varieties,
not necessarily smooth in general. 
Mirror symmetry
predicts a natural correspondence
\beq{mirror1}
{\small\left\{\!\begin{array}{c}
\text{Smooth CY deformations}\\
\text{of $X$}
\end{array}\!
\right\}}
\en\stackrel\sim\leftrightarrow\en
{\small\left\{\!
\begin{array}{c}
\text{Smooth crepant}\\
\text{resolutions of $X'$}
\end{array}\!
\right\}}
\eeq

More precisely, let $\mm_{\mathbf{C}}(X)$ denote
an appropriate local
 moduli space of smooth
complex deformations of $X$ as a
 CY variety. Dually, one has  $\mm_{\mathbf{K}}(X')$,
a so-called  stringy moduli space of K\"ahler resolutions of $X'$.
This  moduli space has not been defined in mathematical terms yet,
but it is expected to be
 closely related to the space of stability conditions
on $X'$, in the sense of T. Bridgeland,  cf. \cite[\S2]{B}.

In the special case where $X'$ has at most finitely many nonisomorphic
 {\em terminal crepant} resolutions of singularities
$\pi: \wt X_j \to X',\,j=1,\ldots,r,$ it is  expected
that there is an open dense subset of   $\mm_{\mathbf{K}}(X')$ of the form
$(\cup_{j=1}^r {\mathbf{K}}_j)/\G$,
where ${\mathbf{K}}_j$ denotes  the complexified K\"ahler cone
of the manifold $\wt X_j ,\,j=1,\ldots,r,$
and where $\G$ is a discrete group of automorphisms
of the corresponding super-conformal field theory.

According to mirror symmetry, one expects to have a  natural 
bijection
\beq{mir2}
\xymatrix{
\mm_{\mathbf{C}}(X)\quad
\ar@{<->}[rr]^<>(0.5){c\leftrightarrow c'}&&\quad
\mm_{\mathbf{K}}(X').
}
\eeq

Write $X_c, X'_{c'}\,$; 
$c\in \mm_{\mathbf{C}}(X), c'\in
\mm_{\mathbf{K}}(X'),$ for a pair of smooth CY manifolds
corresponding to each other under the bijection.
According to a mathematical formulation
of  the mirror conjecture  due
to Maxim Kontsevich, there should be an equivalence
between  the bounded derived category
of coherent sheaves on $X_c,$ on one hand,
and the Fukaya category associated with $X'_{c'}$,
on the other.

One problem with \eqref{mir2} is that there are many singular CY varieties
either without {\em any} smooth deformation or without  {\em any} smooth crepant
resolution, or both. Thus, the moduli spaces in \eqref{mir2}
may be very small indeed, even empty. It is believed, however,
that any  singular CY variety does have a smooth deformation,
as well as a smooth crepant resolution, provided that
deformation, resp. resolution, is allowed to be
{\em noncommutative}, cf. \cite{VdB6}. Furthermore, 
 in order for the bijection \eqref{mir2} to hold one 
{\em must}, in effect, include such
 `noncommutative' deformations, resp. resolutions,
in the corresponding moduli spaces, since
the bijection may  send an ordinary variety into
a noncommutative one, and vice versa.
Kontsevich's  conjecture should also extend to
the noncommutative setting.

Now, it is widely accepted (following the philosophy
advocated by Drinfeld, Kontsevich, and others) that a noncommutative
variety should be thought of entirely in terms 
of the corresponding derived category of coherent sheaves;
in any case, this is the only thing that needs to be
defined in order for the Kontsevich's 
form of the mirror conjecture to make sense.
Thus, one is led to consider 
triangulated categories, to be called {\em CY categories},
which should imitate all the essential
features of the category $D^b\Coh(X)$
for a CY variety $X$.  In particular, one might want to
imitate the following  geometric properties:
\begin{itemize}
\item{$X$ is smooth;}
\item{$X$ is compact;}
\item{The canonical bundle of $X$ is trivial.}
\end{itemize}

An approach to the definition of a  CY category has been proposed
by Kontsevich. First of all, any such category is equivalent
to $D_\perf(A),$  a suitably
defined derived category of finitely generated projective 
left DG modules over
a DG (or, $A_\infty$-) algebra $A$. Then, according to Kontsevich,
the geometric properties above are adequately
encoded in the following properties of the DG algebra
$A$.

\begin{defn}[M. Kontsevich]\label{konts_def}
A $\oper{DG}$ algebra $A$ is called

\npb{{\em homologically smooth} if,
 viewed as a DG {\em bi}module over itself,  $A$ has  a bounded re-
solution by finitely generated projective $\oper{DG}$
 $A$-bimodules, cf. \cite{KS2}, Definition 8.1.2;}
\vskip 1pt

\npb{{\em compact} if $H^\hdot(A)$, the total cohomology of $A$,
 is finite dimensional, cf.  \cite{KS2},\hfill Definition 8.2.1;}\hfill
\vskip 2pt

\noindent$\bullet\en\;${\em CY of dimension $d$} if\;\footnote{this
 definition is temporary and it will be altered shortly.} $A$ is  homologically smooth,
 compact and, moreover, the shift

functor
$M\mto M[d],$ on $D_\perf(A),$  is a {\em Serre functor}, i.e.,
there are bifunctorial isomo-

rphisms
\beq{serrr}
\Hom(M,N)\cong \Hom(N, M[d])^*, \quad\forall M,N\in D_\perf(A).
\eeq
\end{defn}
This definition has to be made precise, of course, by introducing first the
notions 
of  a Serre functor, of a finitely presented DG algebra, and
of a  `projective' DG module, etc. 
Finitely presented DG algebras
are defined as compact objects in an
appropriate category of DG algebras, cf. \cite{TV}.
For projective DG modules see Definition \ref{sf} in \S\ref{selfdual_sec};
`finitely generated and projective'  DG modules are defined, under the name
 `semi-free objects' in \cite[Appendices A,B]{Dr}. We won't go into
 more
details
since we will only use  Definition \ref{konts_def}
in the cases where such notions are
standard and are well known.

\subsection{Definition of Calabi-Yau algebras}
\label{CY_defi} The setting of the present paper differs from
the one suggested by Definition \ref{konts_def}
 in several ways. First of all, we will be working
with ordinary algebras rather than with DG, or $A_\infty$-, algebras.
This restriction, though quite essential, is more of a technical
nature. It corresponds, geometrically, to restricting oneself
to  varieties $X$ such that the category
$D^b\Coh(X)$ has a {\em tilting generator}, see Definition 
\ref{tilt}.
Such a condition holds in many interesting examples, and it
is assumed in order to avoid too much technology from category theory
and homological algebra.
We believe that most of what we will be doing
below may be extended, in principle,  to the general $A_\infty$-setting.

More importantly, the formalism we are going to develop
corresponds, geometrically, to working with
{\em noncompact}, that is, with what are usually called `open'
CY varieties. Thus, our algebras $A$ are
typically neither graded nor 
finite dimensional and, therefore, are {\em not} compact in the
sense of Definition \ref{konts_def}. For such an algebra,
formula \eqref{serrr} involves the duals of infinite dimensional
$\Hom$-spaces. As a consequence, there is no  Serre functor on
$D_\perf(A)$,
the derived category of perfect complexes,
in general (sometimes, but not always, there is a {\em relative}
version of Serre functor, cf. \S\ref{ser_sec}).
So, the very
definition of Calabi-Yau algebras  has to be modified.

Our  definition of Calabi-Yau algebras is best understood in the
general context of    duality. 
There are actually two versions of duality.
The first one involves the functor
$M\mto M^!:=\RHom_{\bimod{A}}(M, A\o A)$
on an appropriate  derived
category of (DG-) $A$-bimodules. This is the derived analogue
of \eqref{DUAL}. Here,
 $\RHom$ is taken with respect
 to the outer bimodule structure on $A\o A$, as in \eqref{DUAL}. The inner structure 
then gives  $M^!$ a natural $A$-bimodule structure again.
For any  morphism $f: M\to N$ (in the derived category) of
$A$-bimodules one
has, by functoriality, the induced morphism $f^!: N^!\to M^!.$

The second version, the Grothendieck duality, involves the
 notion of a
{\em rigid dualizing complex}, an object of the derived
category of (DG-) $A$-bimodules.
In the special case of 
an (ordinary, not DG) {\em Gorenstein} 
algebra $A$, i.e., an algebra of finite injective dimension,
two versions of duality are known to  coincide,
 cf. \cite{Y},\cite{VdB1}.
In that case,  the {\em inverse} of the rigid dualizing complex for $A$ is 
quasi-isomorphic to 
$A^!:=\RHom_{\bimod{A}}(A, A\o A)$, up to shift. This is so, for instance,
for any homologically smooth algebra, provided a  rigid dualizing
 complex for that algebra exists.

\begin{examp}\label{cyx}
Let $X$ be an  algebraic variety of dimension $d$,
and write $\eps: X\into X\times X$ for the diagonal imbedding.
If $X$ is smooth, then the standard Koszul resolution
of the  sheaf $\eps_*\oo_X$ yields
\beq{X_vol}
{\scr E}{xt}_{\oo_{X\times X}}^r(\eps_*\oo_X, \,\oo_{X\times X})=
\begin{cases} \La^d\T_X &\op{if}\enspace r=d\\
0&\op{if}\enspace r\neq d.
\end{cases}
\eeq

Thus, 
we deduce that for
$A=\C[X]$,  the coordinate ring of  a  smooth affine variety $X$,
one has that ~$A^!=\La^d\T(X)[d]$ is indeed the inverse of the canonical
class of $X$.\hfill$\lozenge$
\end{examp}

Motivated by the above discussion,  we 
now replace the notion of
Calabi-Yau algebra given in Definition \ref{konts_def} by the following,
cf. also  \cite{VdB1} and Remark \ref{vdb_crit} below.

\begin{defn}\label{vdb_def} A homologically smooth
($\operatorname{DG}-$)
algebra  $A$ is said to be \cyd (= {\em Calabi-Yau algebra
of dimension} $d\geq 1$) provided 
one has an  $A$-bimodule quasi-isomorphism
\bex f:\ A\iso A^![d]\quad\operatorname{such\en that}\quad
f=f^![d].$\eex
\end{defn}

From now on, we will use
this definition of Calabi-Yau algebras, which is  {\em not} equivalent to Definition
\ref{konts_def}, in general. However, in \S\ref{CYCYpf} we will prove,
cf. also \cite{KS2}, Conjecture ~10.2.8.

\begin{prop}\label{cy=cy_conj} Let $A$ be a homologically smooth,
finitely presented, and {\em compact} $A_\infty$ algebra with unit.
Then, $A$ is  CY algebra of dimension $d$  in the sense of Definition \ref{vdb_def}
iff  the CY condition in
 Definition  \ref{konts_def} holds.
\end{prop}

For an 
ordinary (not DG) algebra $A$,
 the existence of
a (quasi-) isomorphism $A^![d]\cong A$ amounts
to  $A$-bimodule isomorphisms
\beq{cy_def}
\Ext_{\bimod{A}}^k(A,A\otimes A)\stackrel{f}\cong
\begin{cases} A &\op{if}\enspace k=d\\
0&\op{if}\enspace k\neq d.
\end{cases}
\eeq

Given an  associative algebra $A$  that
satisfies  condition \eqref{cy_def},
one may consider 
\beq{Iso}
\Iso_{\bimod{A}}\big(A,\,\Ext_{\bimod{A}}^d(A,A\otimes A)\big),
\eeq
the set  of all $A$-bimodule isomorphisms $f:A\iso \Ext_{\bimod{A}}^d(A,A\otimes
A)$.
The image  of $1_A\in A$ under such an  isomorphism
will be
 referred to as a {\em volume} on $A$.
This is  a {\em central} element
 $\pi=f(1_A)\in\Ext_{\bimod{A}}^d(A,A\otimes A),$ i.e.,
such that
 $a\cdot\pi=\pi\cdot a$ for any $a\in A.$
In terms of $\pi$, the corresponding map $f$ can  be written as
 $f: a\mto a\cdot\pi$.

\begin{notation}\label{cent} Given an associative algebra
 $A$, we
let ${\mathcal Z}(A)$ denote the center of $A$,
and ${\mathcal Z}(A)^\times\sset {\mathcal Z}(A)$ 
the multiplicative group of invertible
elements of the center.
\hphantom{.}\hfill$\lozenge$
\end{notation}

For any  pair $M,N,$ of $A$-bimodules,
 multiplication by  elements of ${\mathcal Z}(A)$
 gives each of the  vector spaces
$\Ext^r_{\bimod{A}}(M,N), \,r=0,1,\ldots,$
a ${\mathcal Z}(A)$-module structure.
In particular,  one has
$\Hom_{\bimod{A}}(A,A)=~{{\mathcal Z}(A).}$ 
We see that, for 
 an associative algebra $A$  that
satisfies \eqref{cy_def},
the set \eqref{Iso} is a ${\mathcal Z}(A)^\times$-torsor.
Thus, we conclude that the set of volumes on such an 
algebra $A$ is also   a ${\mathcal Z}(A)^\times$-torsor.

\begin{rem}\label{self-d} 
Following M. Kontsevich we observe that,
for an algebra $A$ satisfying  \eqref{cy_def},
the assignment $f\mto f^![d]$ gives
an  involution of the set \eqref{Iso}.
By functoriality, this involution commutes with
the  ${\mathcal Z}(A)^\times$-action.
We deduce that there exists a unique 
$z\in {\mathcal Z}(A)^\times$ such that $z^2=1$ and such that for any element
$f$ from the the set  \eqref{Iso}
we have $f^![d]=z\cdot f.$ 
It follows from our Theorem \ref{main}(ii) below that
one has in effect $z=1$, provided the algebra $A$ is
{\em friendly}
in the sense of Definition \ref{ww}.
Thus, in the definition of CY algebras,
the requirement that  $f^![d]=f$ is {\em a posteriori} superfluous, 
at least for
friendly algebras.

It was pointed out to me by M. Kontsevich that,
for  general DG or $A_\infty$ algebras, the condition  $f^![d]=f$
in Definition \ref{vdb_def} is essential and is, moreover,
quite subtle;
it should involve `higher homotopies'.
\end{rem}

\subsection{Crepant resolutions}\label{crepant}
In the special case $A=\C[X]$, considered
in Example \ref{cyx},
 giving a volume  amounts, according to \eqref{X_vol}, to
 giving a nowhere vanishing global regular section of
$\La^d\T_X$,  $d=\dim X$.
The existence of
such a section is nothing but the standard definition
of a Calabi-Yau variety.

The class of examples arising from smooth {\em affine} varieties
has an important generalization to some nonaffine varieties. Specifically, 
 an algebraic variety $X$ is said
to be  projective over an affine variety,
if there exists an affine variety $Y$ and a projective
morphism $X\to Y$
(typically, one takes $X$ to be a resolution of singularities of
a singular affine variety $Y$).

\begin{prop}\label{CYCY} Let $X$ be  a smooth connected
variety which is projective
 over an affine variety, and let $\CE\in D^b(\coh X)$ be
a tilting generator,
 cf. Definition \ref{tilt}. Then, we have 
\vskip 1pt

\vi The algebra $\End\CE:=\Hom_{D^b(\coh X)}(\CE,\CE)
$ is a finitely generated ${\mathcal Z}(\End\CE)$-module, and
multiplication by  regular functions induces
natural algebra
isomorphisms
 $$\Gamma(X,\oo_X)\cong \DZ(\coh X)\cong {\mathcal Z}(\End\CE).$$

\vii The algebra $\End\CE$ is a CY algebra of dimension $d$  iff
$X$ is a  Calabi-Yau manifold of dimension $d$. 
\end{prop}
\noindent
Here,  $\DZ(\coh X)$ denotes the center
of $D^b(\coh X)$, that is, the endomorphism algebra
of the identity functor $\Id: D^b(\coh  X)\to D^b(\coh  X).$

 Proposition   \ref{CYCY} will be deduced from Theorem \ref{CM} of \S\ref{CYCYpf}.
Closely related to this, is also the following 
special case of Theorem \ref{CM}.

\begin{cor}\label{crep_cor} Let $A$ be an algebra of   finite Hochschild
dimension and let $Z$ be a central commutative subalgebra
such that 

\npb{$\Spec Z$ is a smooth affine Calabi-Yau variety;}

\npb{$A$ is a finitely generated $Z$-module.}

Then, the following  conditions are equivalent:
\begin{enumerate}

\item $A$ is a  Calabi-Yau  algebra of dimension  $d:=\dim(\Spec Z)$;

\item $A$ is a  projective $Z$-module, and
 there  exists
a nondegenerate   trace  $\phi: A\to Z$.
\end{enumerate}
\end{cor}

\noindent 
Recall that a  $Z$-linear map $\phi: A\to Z$  is called a 
{\em nondegenerate   trace} if 
 the assignment $1\mto \phi$ extends to 
an  $A$-bimodule isomorphism $A\iso \Hom_Z(A,Z)$,
where the space $\Hom_Z(A,Z)$ is equipped with the
natural $A$-bimodule structure given by
$(a\cdot f\cdot b)(x):= f(bxa),$ for any $f\in \Hom_Z(A,Z)$
and $a,b,x\in A$. Any nondegenerate   trace is automatically 
{\em symmetric}, i.e., one has
$\phi(ab)=\phi(ba),\,\forall a,b\in A.$

\subsection{Batalin-Vilkovisky structure on Hochschild
 cohomology}\label{BVHOCH}  Given an associative (DG-)
algebra $A$, we let
$H_\idot(A,A):=\Tor^{\bimod{A}}_\idot(A,A),$
resp. $H^\hdot(A,A):=\Ext_{\bimod{A}}^\hdot(A,A)$, denote Hochschild homology,
resp. cohomology, of $A$. 
Hochschild homology comes equipped with the standard {\em Connes differential}
 $\BB: H_\idot(A,A)\to H_{\hdot+1}(A,A),$
a noncommutative analogue of the de Rham differential on differential forms.

Hochschild  cohomology has a natural structure of  {\em Gerstenhaber
algebra}
given by  the cup-product and   Gerstenhaber bracket
$H^{k-1}(A,A)$ $\times H^{l-1}(A,A)\to H^{k+l-1}(A,A)$.
Furthermore, for any $l\geq k$, there is a natural {\em contraction} action
$H^k(A,A)\times H_l(A,A)\map H_{l-k}(A,A)$ that makes
  Hochschild homology a Gerstenhaber {\em module} over
the  Gerstenhaber algebra $H^\hdot(A,A).$

Let $A$ be a Calabi-Yau (DG-) algebra 
of dimension $d,$ with a volume $\pi$.
According to  M.
Van den Bergh \cite[Thm.~1]{VdB1},  for any $A$-bimodule $M$,
there is a {\em functorial}
Poincar\'e duality type isomorphism (that depends on the choice of
$\pi$)
\beq{vdb}
{\mathbb D}:\
\Tor^{\bimod{A}}_\idot(A,M)\iso
\Ext_{\bimod{A}}^{d-\hdot}(A,M).
\eeq
\begin{rem}\label{vdb_crit} For $M=A\o A$, the above isomorphism
reduces to condition \eqref{cy_def}. Thus,  the existence of an isomorphism
of functors
$\Tor^{\bimod{A}}_\idot(A,-)\cong
\Ext_{\bimod{A}}^{d-\hdot}(A,-)$ may serve
as a replacement of the condition $A^![d]\cong A$.
It follows that the notion of a Calabi-Yau
 DG algebra is an invariant
of the pair $\big(D_\perf(\bimod{A}),\, -\stackrel{_L}\o_A-\big)$, viewed as
a triangulated category with monoidal structure given
by the (derived) tensor product of bimodules.
Note that $A$, the  diagonal bimodule, may be characterized as the
unit of the monoidal structure; on the contrary, the bimodule $A\o A$
(hence, also condition \eqref{cy_def} that involves this bimodule) 
apparently has {\em no}
intrinsic  characterization in terms involving the
triangulated category with monoidal structure  only.
\erem

We use the  (degree reversing) isomorphism ${\mathbb D}$ in \eqref{vdb}
in the special case where $M=A$, and
transport the  Connes differential on  Hochschild homology
to obtain  a degree $-1$ differential,
 $\Delta_\pi,$ on  Hochschild {\em co}homology. Thus, we get a commutative diagram
$$
\xymatrix{
H_\idot(A,A)\ar@{=}[d]_<>(0.5){{\mathbb D}}\ar[rr]^<>(0.5){\BB}&&
H_{\idot+1}(A,A)\ar@{=}[d]_<>(0.5){{\mathbb D}}\\
H^{d-\hdot}(A,A)\ar[rr]^<>(0.5){\Delta_\pi}&&H^{(d-\hdot)-1}(A,A).
}
$$

The following result
is a noncommutative generalization of the standard
 BV algebra structure on the algebra of polyvector
 fields on a
 Calabi-Yau manifold, cf. eg. ~\cite{Sch}.
\begin{thm}[BV structure on Hochschild cohomology]
\label{BV} Let $A$ be a   CY algebra of dimension $d$  with
volume   element $\pi\in \Ext_{\bimod{A}}^d(A,A\o A).$ Then,

\vi The map ${\mathbb D}$  in \eqref{vdb}
intertwines contraction  and  cup-product maps, i.e., we have
$$ {\mathbb D}(i_\eta c)=\eta\cup{\mathbb D}(c),
\quad\forall c\in H_\idot(A,A),\,\eta\in H^\hdot(A,A).
$$

\vii 
The differential 
$\Delta=\Delta_\pi$ makes Hochschild cohomology
a Batalin-Vilkovisky algebra; in other words, for the Gerstenhaber bracket
 on Hochschild cohomology one has
$$\{\xi,\eta\}= \Delta(\xi\cup  \eta)-\Delta(\xi)\cup  \eta -
(-1)^{|\xi|}\cd \xi\cup  \Delta(\eta),
\qquad \forall \xi,\eta\in H^\hdot(A,A).$$
\end{thm}

This Theorem was first announced without proof in
\cite[Claim 4.5.6]{CBEG}. The proof will be given in \S\ref{bv_pf} below.

In low degrees, we have  $H^0(A,A)={\mathcal Z}(A)$,
and  $H^1(A,A)=$ $\Der(A,A)/\Inn(A,A),$
a quotient  of the space of  derivations $A\to A$ by the subspace of
inner derivations. The Gerstenhaber bracket on $H^1(A,A)$ is
 induced by the commutator of 
derivations. The above theorem implies in particular
that the operator $\Delta$ associated with
a   volume $\pi$ on $A$ gives rise to
a natural divergence map
$\div: \Der(A,A)\to {\mathcal Z}(A),\,\xi\mto\div\xi.$ This map vanishes on inner derivations
and is a Lie cocycle, ie., one has
$$\div[\xi,\eta]=\xi(\div\eta)-\eta(\div\xi),\quad\forall
\xi,\eta\in\Der(A,A).
$$

\begin{rem}\label{tradler} T. Tradler has defined in \cite{Tr} a similar BV structure on
Hochschild cohomology 
of a compact
$A_\infty$-algebra with an inner product. Tradler's construction
does not overlap with ours since our setting does not involve
 inner products, cf. however Remark ~\ref{Raf}. 

Roughly speaking, the construction
 in \cite{Tr} uses the isomorphism $A\cong A^*$ arising
from  an inner product on $A$ to obtain an isomorphism
$H^\hdot(A,A)\cong H^\hdot(A,A^*)$, and to 
transport Connes' differential from $H^\hdot(A,A^*)$
via the latter isomorphism.
The isomorphism $H^\hdot(A,A)\cong H^\hdot(A,A^*)$ is,
however, very different from the one in \eqref{vdb}
since it only involves Hochschild  {\em co}homology
and does {\em not} reverse cohomology degrees: $(\hdot)\mto d-(\hdot).$
In particular, it does not give rise to any divergence
map $\div\!:\Der(A,A)\to {\mathcal Z}(A)$.
\erem

\subsection{General setting}\label{general} In order to explain a
universal
construction of CY algebras we need first to introduce the notion
of a smooth algebra.
There are actually {\em two totally different} notions
of smoothness in noncommutative  geometry.
The corresponding two  types of smooth algebras play different roles.

The algebras of the first type, called just `{\em smooth}' in
Definition \ref{ww} below, serve as a proper generalization
for the notion of a {\em free} associative algebra.
The condition of being smooth in that sense is
extremely restrictive.  Free tensor algebras
and path  algebras of quivers are smooth, and any smooth algebra
is, in a way, close to being
free, cf. \cite{CQ} for some examples. On the contrary,  the coordinate
ring of a 
{\em smooth}
affine algebraic variety $X$ 
is {\em not}   smooth in the sense of Definition  \ref{ww} unless $\dim X\leq 1$.

The second, much weaker, notion of smoothness
involves homological smoothness.  The corresponding 
class  of {\em friendly algebras} introduced in
Definition \ref{ww} below is, we believe,
an adequate  noncommutative  generalization of the class of
coordinate rings of smooth affine schemes of finite type,
in conventional algebraic geometry.

\begin{defn}\label{ww} An associative  algebra $F$ is said to be

\npb{{\em coherent}, resp. left or right coherent, if the
kernel  of
any morphism
 between finite rank free  $F$-bimodules, resp. left or right  $F$-modules,
is finitely generated;}
\vskip 2pt

\npb{{\em friendly}  if it is homologically smooth, finitely
presented,
and  coherent.}
\vskip 2pt

\npb{{\em smooth}, if  $\Om^1 F$
is a finitely generated, and projective $F$-bimodule.}
\end{defn}

\begin{rem} \vi
For a smooth algebra $F$, the schemes $\Rep_\dd F$ are smooth
for all $\dd\geq 1$, see eg. \cite{Gi2}, \S 20.
 A free product
of smooth algebras is a smooth  algebra.
Any  smooth algebra  is both left and right
coherent,
but not necessarily coherent, in general, cf. eg.,
\cite[\S19.3]{Gi2}.

\vii
Any bimodule noetherian algebra  is  coherent,
resp. any left or right noetherian algebra is
left or right coherent. 
 Path algebras of finite
quivers are smooth but not  left or right noetherian, in general.
 It is  unknown  to the author
 whether or
not a free finitely generated algebra on $n >1$ generators is
({\em bimodule}) coherent.
\erem

A key role in what follows will be played by a {\em reduced
contraction}  map introduced in \cite{CBEG}; this map
 is analogous to contraction of differential forms
with a vector field, in differential geometry. To define 
reduced contraction, recall first that
 any  double derivation $\th\in\dder F$
gives, by definition, an $F$-bimodule map
 $\Om^1 F\to F\o F.$ This map can be uniquely
extended to a map
$i_\th: \Om^\hdot F\to\Om^\hdot F\o\Om^\hdot F,$
a degree $-1$ super double derivation
of the graded algebra $\Om^\hdot F$, cf. \cite{CBEG}.

\begin{defn} 
The {\em reduced} contraction map
is defined by the pairing
\beq{reduced}
\bi:\ (\Om^p F)\br\times
\dder F\to \Om^{p-1} F,\quad \al\times\th\mto \bi^\al(\th)=
\bi_\th\al:=m^\opp(i_\th\al).
\eeq
Here, for any graded algebra $B^\hdot=\bigoplus_{p_{_{}}\in\Z} B^p$ with
multiplication map $m: B^\hdot\o B^\hdot\to B^\hdot,$ we put
$\dis
m^\opp(a\o b):=(-1)^{p\cdot q}\cdot b\cdot a\in B^{p^{^{}}+q},$
for any $a\in B^p,\,b\in B^q.$
\end{defn}

Depending on whether we fix the first or  second
argument, we may view the pairing 
in \eqref{reduced}
either as a map $\bi^\al:  \dder F\to \Om^{p-1} F,$
or as  a  map $\bi_\th:(\Om^p F)\br\to \Om^{p-1} F.$
Similar convention will be applied for the nonreduced contraction
$\al\times\th\mto i_\th\al=i^\al(\th)$ as well.

In the most important special case of $p=1$,
 the reduced contraction map has the following explicit form
\beq{i}
\bi^\al:\
\dder F\to F,\quad
\th \mto \th''(b)\cdot a\cdot
\th'(b),\quad\text{if}\quad\al=a\,db\in \Om^1F,\;a,b\in F.
\eeq

The definition of reduced contraction is motivated, in part, by the following 

\begin{examp} Let $F=\C\langle\xxx\rangle$, and choose a potential
 $\Phi\in F\br$. In the case of
 an {\em exact} cyclic 1-form $\al$, 
it is easy to
see that for $\bi_\th\al$, in $F$, we have 
$$\bi_{\th_j}(d\Phi)=\pa \Phi/\pa x_j,
\qquad\text{for}\en\al=d\Phi\en\text{and}\en \th=\th_j=\pa/\pa x_j,\en\; j=\nn,
$$
where  the notation is the same as in formulas \eqref{partial} and
 \eqref{ThF}.\hfill$\lozenge$
\end{examp}

\begin{examp}\label{qexamp1} Let
$F:=\C[x,x\inv]*\C\langle y,z\rangle=\C\langle x^{\pm1},y,z\rangle,$
and define double derivations $\pa/\pa x,\,\pa/\pa y,\,\pa/\pa z:
F\to F\o F$  similarly to the previous example.
The algebra $F$ is smooth, and  $\dder F$ is a free $F$-bimodule
with basis $\{\pa/\pa x,\,\pa/\pa y,\,\pa/\pa z\}$.

For any $\Phi\in F\br,$ we introduce the algebra
\beq{QA}
\A(F,\Phi):=\C\langle x^{\pm1},y,z\rangle
\big/\llb\bi_{_{\pa/\pa x}}d\Phi,\,\bi_{_{\pa/\pa y}}d\Phi,\,\bi_{_{\pa/\pa
z}}d\Phi\rrb.
\eeq

We consider a special case where $\Phi=xyz-q\cdot yxz +f,$ for
some {\em Laurent} polynomial $f\in\C[x,x\inv]$,
cf. Example \ref{qexamp}.
The  relations in the corresponding algebra \eqref{QA} read
$$x\cd y\cd x\inv=q\cd y,\quad
x\cd z\cd x\inv=q\cd z,\quad
yz-q\cd zy=g(x),\quad
g(x) :=\mbox{\large$\frac{f'(x)}{x}$}\in\C[x,x\inv].
$$

This algebra is a deformation of $\U_q({\mathfrak{sl}}_2)$,
the quantized enveloping algebra of the Lie algebra ${\mathfrak{sl}}_2$.
If $q\in \C$ is a root of unity, then the algebra $\A(F,\Phi)$ 
has a large center $\CZ(F,\Phi)$, cf.  Example \ref{qexamp}.

On the other hand, we may view $q$ as an indeterminate and put $q=e^\eps$.
This way, one makes $F$ a $\C[[\eps]]$-algebra, $F_\eps$.
We also  have  $\A_\eps(F,\Phi)$,
the corresponding $\C[[\eps]]$-algebra quotient.
One can show that the center  of the
algebra  $\A_\eps(F,\Phi)$ is generated by a single
element $\phi$, a deformation of the Casimir element
in $\U({\mathfrak{sl}}_2).$ 
The algebra $\A_\eps(F,\Phi)/\llb\phi\rrb$ is
the flat  noncommutative  deformation of the coordinate ring  of type $\mathbf{A}$
Kleinian singularity studied by T. Hodges \cite{Ho}.\hfill$\lozenge$
\end{examp}

\subsection{Universal construction of CY algebras}\label{potential_sec}
We are going to present a universal construction, based
on noncommutative symplectic geometry, that produces all CY algebras of
dimension $d$. Our  construction will be
a generalization of   the construction
of the DG algebra $\FD(F,\Phi)$.

For a graded algebra $D$, any (super)-derivation $\xi: D\to D$
gives rise to a map $L_\xi:\dder D\to\dder D$,
resp. $L_\xi: \Om^\hdot  D\to\Om^\hdot  D,$
called {\em Lie derivative}.
Thus, given a DG 
algebra $(D,\xi)$, the spaces $\dder D$ and $\Om^1 D$ acquire
natural  structures of DG $D$-bimodules, with differential $L_\xi$.

\begin{defn}[cf. \cite{CBEG}]\label{sympl} Let $(D,\xi)$ be a \DG
algebra and $\om\in (\Om^2 D )\br$
a cyclic 2-form  such that $L_\xi\om=0$ 
and  $d\om=0$.

The 2-form $\om$
is called
{\em symplectic}, resp., {\em homologically symplectic}, if  the
map $\bi^\om: \dder  D\to\Om^1  D$, resp.,
the map $H_\idot(\bi^\om): H_\idot(\dder  D, L_\xi)\to
H_\idot(\Om^1  D,L_\xi)$, is a bijection.\hfill$\lozenge$
\end{defn}

The ingredients of our 
construction involve the following
\vskip 5pt

\noindent
{\sc{Symplectic DG data:}}\en
An integer $n\geq 0,$ and a
 triple $(D, \om, \xi),$ where 
\begin{align}
\mathsf{(i)}\;\en&\parbox[t]{80mm}{\text{$D=\bigoplus_{k\geq 0}
D_k$,
is a {\em smooth} graded  algebra  such that $\Om^1 D,$ viewed as a
graded}}\nonumber\\
&\parbox[t]{80mm}{\text{$D$-bimodule, is generated by homogeneous elements of degree $0,1,\ldots,
n;$}}\nonumber\\
\mathsf{(ii)}\;\en&\parbox[t]{80mm}{\text{$\om\in
(\Om^2 D)\br$ is a homogeneous symplectic
 2-form
of grade degree $n$;}}\label{D}\\
\mathsf{(iii)}\;\en&\text{$\xi: D_\idot\to D_{\idot-1}$ is a
 super-derivation  such that $L_\xi\om=0$,
$\,\xi^2=0,$ and $\C\cap \xi(D)=0.$}\nonumber
\end{align}
\vskip 2pt

Thus, conditions \eqref{D}(i)-(iii) insure that
 $(D, \xi)$ is a DG  algebra; in
\eqref{D}(ii), the  grade degree of the 2-form
is taken with respect to the grading $(\Om^2 D)\br=\bigoplus_{r\geq 0}
(\Om^2_r D)\br$
 induced from the grading
on $D$. 

We say that a data  $(D, \om, \xi)$ is
{\em homologically symplectic}
if, in the condition  \eqref{D}(ii), the requirement that the 2-form
$\om$ be symplectic is replaced by a weaker requirement
that  $\om$  be only  homologically symplectic.

We have the distinguished double derivation
$\De: D\to D\o D,$ cf. \eqref{De}.
The corresponding
reduced  contraction $\bi_\De: (\Om^p D)\br\to \Om^{p-1} D$
is a grading preserving map that commutes with
any super-derivation on $D$,  with
$\xi$ in particular. In the special case $p=1$,
see \eqref{i},  one finds 
$\bi_\De(a\,db)=ab-ba,$ 
for any 1-form $\al=a\,db\in (\Om^1 F)\br.$

According to  \cite{CBEG}, there exists an element
 $\omh\in  D$ such that, in $\Om^1 F$,
one has  $d\omh=\bi_\De\om$.
It was shown in  \cite{CBEG}  that such an  element  is unique up to a
 summand from $\C
\sset D$.
 Further, it is clear that the 1-form $\bi_\De\om$ is  homogeneous 
of grade degree $n$. Thus, one may choose $\omh\in D_n$;
we will always make such a choice from now on.
This determines $\omh$ uniquely, unless $n=0$.
\smallskip

\noindent
{\bf{Main construction.}}\en
For any  symplectic data $(D, \om, \xi),$ as in
\eqref{D}, we  define
 a DG algebra structure on the free product
 $\FD:=D *^{\,}  \C[t]$ as follows.

Let $\deg t:=n+1,$
and assign all elements of the subalgebra $D\sset D *^{\,}  \C[t]$
their natural degree.
Further, we use the element $\omh\in D_n$ to define  a degree $-1$ super-derivation
\beq{ddata}\d: \
\FD_\idot\too  \FD_{\idot-1},
\quad t\mto\omh,\quad u\mto \xi(u), \en\forall u\in D\sset
 \FD.
\eeq

By  \eqref{D}(iii), we have  $L_\xi\om=0$.
It follows that $\xi(\omh)\in \C$. Therefore,
 the condition that $\C\cap \xi(D)=0$, see \eqref{D}(iii), yields
$0=\xi(\omh)=\d^2(t).$ We deduce that $\d^2=0$. Thus,
we get a DG algebra $(\FD,\d)$. We will use
the following notation for this DG algebra, resp. its 0-th
homology,
which is an associative algebra,
$$
 \FD=\FD(D, \om, \xi):=D *^{\,}  \C[t], \quad\text{resp.},\quad\A=\A(D, \om, \xi):=H_0(\FD,\d).
$$

The following result, along with Conjecture \ref{uniconj} below,
provides a universal construction of CY algebras.

\begin{thm}\label{main} \vi Let  $(D, \om,
\xi)$ be a symplectic data such that
\eqref{D}$\mathsf{(i)}$-$\mathsf{(iii)}$ hold and, moreover, such that
$\dis H_k(\FD,\d)$ $=0$  for all $k\neq 0$.
Then, the algebra 
$\dis\A(D, \om, \xi)$ is a CY algebra of
dimension $d=n+2$.
\vskip 3pt

\vii Let $A$ be  a friendly algebra that satisfies condition
 \eqref{cy_def} for $d\geq 2$.
Then, $A$ is a CY algebra of dimension $d$; furthermore,
there exists a {\em homologically} symplectic data $(D, \om, \xi)$ with  ${n=d-2}$
such  that
$$
A\cong \A(D, \om, \xi)=H_0(\FD,\d)
\quad\oper{and,\en moreover,}\quad
H_k(\FD,\d)=0,\quad\forall k\neq 0.
$$
\end{thm}

In the special case $d=2$, we have $d-2=0$, so $D=D_0$; hence,
Theorem \ref{main}(ii) gives

\begin{cor}[CY algebras of dimension $2$]\label{cytwo}
 Any friendly CY  algebra of dimension $2$
has the form $A=D/D\omh D$, where $D$ is a smooth  algebra,
$\om\in (\Om^2 D)\br$ is a symplectic $2$-form, and $\omh\in D$ is
an element such that $d\omh=\bi_\De\om$.\qed
\end{cor}

\begin{conj}\label{uniconj} In Theorem \ref{main}(ii), one can
replace `homologically symplectic' by `symplectic'.
\end{conj}

The  proof of  Theorem \ref{main} is rather long; it
occupies the rest of \S\ref{cy_sec}.
The reader may skip the  proof and go directly to \S\ref{alg}.

\subsection{Preparation for the proof  of Theorem \ref{main}}\label{prep}
For any  algebra $B$ and a cyclic $p$-form $\al\in(\Om^p B)\br$, the map
$\bi^\al: \dder B\to \Om^{p-1} B$ is easily seen to be
a morphism of $B$-bimodules.
Conversely, we have

\begin{lem}\label{FF} Let $B$ be a {\em smooth} algebra. Then,
the assignment $\al\mto\bi^\al$ 
yields the following bijections
\begin{align}
&(\Om^1 B)\br\iso
\Hom_{\bimof{B}}(\dder B,B);\label{bi1}\\
&(\Om^2 B)\br\iso
\{f\in \Hom_{\bimof{B}}(\dder B,\Om^1 B)\en\big|\en f=f^\vee\}.
\label{bi2}
\end{align}
\end{lem}
\begin{rem}
A map $f$ such that  $f=f^\vee$ is called 
{\em selfadjoint}. The meaning of the  equation  $f=f^\vee$ is as follows.

For any  $B$-bimodule morphism $f: M\to N$,
one has the dual morphism $f^\vee: N^\vee\to M^\vee$.
In the special case where $M$ is a   finitely generated,  projective
$B$-bimodule
and $N=M^\vee$, we have $N^\vee=(M^\vee)^\vee=M,$ canonically.
Thus, the dual morphism is a map $f^\vee: M\to M^\vee$, again.
Now, in the case of a  smooth algebra $B$, one can take
$M:=\dder B$ $=(\Om^1 B)^\vee$, a   finitely generated  projective
$B$-bimodule such that $M^\vee=\Om^1 B$.
\end{rem}
\begin{notation}\label{e} Write
 $B^e=B\o B^\opp$, where $B^\opp$ denotes the opposite algebra.
\hfill$\lozenge$
\end{notation}

\begin{proof}[Proof of Lemma \ref{FF}]
For any pair of $F$-bimodules, that is,  left $F^e$-modules, $P,Q$, one has
a diagram
\beq{Pvee}
\frac{P^\vee\o_F Q}{[F, P^\vee\o_F Q]}\cong P^\vee\o_{F^e}Q\cong \Hom_{\Lmod{F^e}}(P,F\o
F)\o_{F^e}Q
\stackrel{u}\too\Hom_{\bimod{F}}(P,Q).
\eeq

If  $P$ is   finitely generated and projective 
then so is $P^\vee$, and the map $u$ in \eqref{Pvee}
is an isomorphism.

Assume now that $F$ is   smooth.
Then, 
taking $P:=\dder F$ and $Q:=F$, we get
$P^\vee=((\Om^1 F)^\vee)^\vee=\Om^1 F$.
Furthermore, the map $u$ is an isomorphism, hence the composite map in \eqref{Pvee}
yields an isomorphism
$\Om^1 F/[F,\Om^1 F]\iso\Hom_{\bimod{F}}(\dder F,F).$
It is straightforward to verify that 
this map is indeed given by the assignment $\al\mto\bi^\al$,
and \eqref{bi1} follows.

To prove \eqref{bi2}, we take
$P=\dder F$ and $Q=\Om^1 F$. We have
$(\dder F)^\vee\o_F\Om^1 F=\Om^2 F.$ Hence,  using \eqref{Pvee} we obtain
$$ \Om^2 F/[F,\Om^2 F] \cong
\Hom_{\bimof{B}}(\dder B,\Om^1 B).
$$

Let $\phi^\al: \dder B\to\Om^1 B$ denote the map
associated to a 2-form $\al\in \Om^2 F$ via the
above isomorphisms. It is then straightforward to check
by comparing an explicit formula for $\phi^\al$ with the one
for the map $\bi^\al$
that $\phi^\al=\bi^\al$ holds iff the map $\phi^\al$ is selfadjoint.
\end{proof}

By a DG  algebra we always mean a graded  algebra
 $\FD=\bigoplus_{r\geq 0}\FD_r,$
with differential $\d: \FD_\idot\to\FD_{\idot-1},$ such
that $\d(\FD)\cap \C=0.$
We   need the following result, to be proved in 
~\S\ref{prop_pf}.
\begin{prop}\label{kunneth}
Let $(\FD,\d)$ be a smooth DG  algebra.
Put $A=H_0(\FD,\d)=\FD_0/I$,
where $I:=\d(\FD_1)$ is a  two-sided ideal
in $\FD_0$.  Then, we have

\noindent
\vi If the map $\FD\qis A$ is a quasi-isomorphism,
then each of the following 3 maps is also a quasi-isomorphism
$$
\Om^1 \FD\stackrel{1}\onto\Om^1 (\FD|A)\stackrel{2}\onto\Om^1 A,\quad\oper{and}\quad
\dder \FD\stackrel{3}\onto\dder (\FD|A).
$$

\noindent
\vii  Assume that the map $\Om^1 (\FD|A)\onto\Om^1 A$ is a
 quasi-isomorphism
and, moreover, either of the following two conditions  holds:

\npb{Each homogeneous component of the algebra  $\FD$ is complete
in the $I$-adic topology;}

\npb{The algebra  $\FD$ is {\em bi}graded,
$\FD=\bigoplus_{r,k\geq 0}\FD_r(k)$, the differential
$\d: \FD_\idot(k)\to\FD_{\idot-1}(k)$ preserves the second grading and, moreover, we have
$\FD_\idot(0)=\C.$}
\vskip 2pt
\noindent
Then, the projection  $\FD\onto A$ is a quasi-isomorphism.
\end{prop}

To summarize,   with the assumptions of  Proposition \ref{kunneth}(ii),
 we have
\beq{p}
\xymatrix{
\FD\ar@{->>}[rr]^<>(0.5){p}_<>(0.5){\text{qis}}&& A&
\ar@2{<->}[r]&&
\Om^1 (\FD|A)\ar@{->>}[rr]^<>(0.5){p_{\FD|A}}_<>(0.5){\oper{qis}}&& \Om^1 A.
}
\eeq
\subsection{Proof  of part (i) of Theorem \ref{main}}\label{sketch}
The proof  proceeds in three steps.

\step{1.} 
We consider a diagram, in which $\jmath_D$ is the imbedding from ~\eqref{morph},
\beq{sympl_diag}
\xymatrix{
D\o D\ar@{=}[d]^<>(0.5){\Id}\ar[rrr]^<>(0.5){\ad
:\,\;1\o1\,\mto\, \De}&&&\dder D
\ar[d]^<>(0.5){\bi^\om}\ar[rrrr]^<>(0.5){i^{d\omh}:
\;\th\mto \th'(\omh)\o \th''(\omh)}&&&&
D\o  D\ar@{=}[d]^<>(0.5){\Id}\\
D\o D\ar[rrr]^<>(0.5){\nu:\,\;1\o1\,\mto\, d\omh}&&&
\Om^1D\ar[rrrr]^<>(0.5){\jmath_{_D}=i_\De}&&&&D\o  D.
}
\eeq

\begin{lem}\label{lem1} \vi Diagram \eqref{sympl_diag} commutes; furthermore the vertical map
$\bi^\om$ 
commutes with the Lie derivative ~$L_\xi$.

\viii The rows of the diagram are obtained from each
other by applying the duality $(-)^\vee;$
 in particular,
we have $\nu^\vee=-i^{d\omh}$ and $\ad^\vee=\jmath_D$,
cf. \eqref{DUAL}.
\end{lem}

\proof Commutativity of the left square can be checked
on the element $1\o1$. This amounts to
the equation $\bi_{\De}\om=d\omh,$ which is the definition of
$\omh$.  Commutativity of the right square
is equivalent to the equation
$i_\De\bi_\th\om=-i^{d\omh}(\th)$. Using Proposition \ref{symmetric}(i),
see \S\ref{nchess},
we compute (we use the notation $(x\o y)\y=y\o x$):
$$i_\De\bi_\th\om=-(i_\th\bi_\De\om)\y=-(i_\th d\omh)\y=
-(\th(\omh))\y=i^{d\omh}(\th),
$$
and commutativity of the diagram follows.

To prove that the map  $\bi^\om$ commutes with the Lie derivative
$L_\xi$,
we use  a standard identity
$$ \bi_{L_\xi\th}= \bi_\th \ccirc L_\xi-L_\xi\ccirc\bi_\th,
\quad\forall\xi\in\Der (D,D),\;\th\in\dder D.
$$
Since $L_\xi\om=0$, 
the above identity yields
$\dis \bi^\om(L_\xi\th)=L_\xi(\bi^\om(\th)),$ as required.
Part (i) follows.

The dualities $\ad^\vee=\jmath_D$ and $(i^{d\omh})^\vee=\nu,$ in part
(ii), 
both follow from
\eqref{morph}. The equation $(\bi^\om)^\vee=\bi^\om$ is due
to \eqref{bi2}. This proves part (ii).
\qed
\vskip 2pt

\step{2.} The definition $\FD=D *  \C[t]$ yields
 a $\FD$-bimodule isomorphism
\beq{omt}
\FD\o \FD\;\boplus\;\Om^1 (D|\FD) \iso
\Om^1 \FD,\quad(p'\o p'')\oplus 
\be  \mto
p''\cd dt\cd  p'+ \be
\eeq
Tensoring \eqref{omt} by $A$ on each side, we deduce
$$
\Om^1 (\FD|A)=
A\o_\FD\left(\FD\o \FD\;\boplus\;\Om^1 (D|\FD)\right)\o_\FD A
=A\o A\;\boplus\;\Om^1 (D|A).
$$
Thus, we obtain an $A$-bimodule isomorphism
\beq{odda}
\cone\big[p_{\FD|A}: \Om^1 (\FD|A)\to \Om^1 A\big]\cong
A\o
A\;\boplus\;\Om^1 (D|A)\;\boplus\;\Om^1 A.
\eeq

The  natural DG module structure on the
cone of a morphism
makes the LHS of \eqref{odda} a DG $A$-bimodule.
Hence, the direct sum  in the RHS acquires a  DG bimodule structure
as well. To write the resulting differential
explicitly, it is convenient to
interpret the RHS in  \eqref{odda} as a
3-term complex of
 $A$-bimodules with differential $\nu$,
\beq{FDG}
\xymatrix{
A\o A\;\ar[rrrr]^<>(0.5){\nu_{_{1}}:\,\;1\o1\,\mto\, d\omh}&&&&
\;\Om^1 (D|A)\ar[rr]^<>(0.5){\nu_{_0}=p_{D|A}}&&\Om^1 A.
}
\eeq

The middle term in \eqref{FDG} is itself a DG $A$-bimodule with 
differential $L_\xi$, the Lie derivative with respect to the
odd super-derivation $\xi: D_\idot\to D_{\hdot-1}$.
We  view the other two terms
in \eqref{FDG} as  DG bimodules with zero differentials.
This makes \eqref{FDG}   a double complex,
with two (anti)-commuting differentials, $\nu$ and $L_\xi$.
One can verify that the differential
in the RHS of \eqref{odda}
obtained by  transporting the natural
differential from the LHS equals
 $\nu+L_\xi,$ the total differential.

Now, the assumptions of Theorem \ref{main}(i)
say that $H_r(\FD,\d)=0$ for all $r>0$.
Thus, Proposition \ref{kunneth}(i) implies
the map $p_{\FD|A}: \Om^1 (\FD|A)\to \Om^1 A$ is a quasi-isomorphism,
cf. \eqref{p}.
Using the isomorphism in
\eqref{odda} we conclude  that 
\eqref{FDG}, viewed as a  DG $A$-bimodule with respect to the total
differential $\nu+L_\xi$, is  {\em acyclic}.
\vskip 2pt

\step{3.} 
Observe that the composite map in either
row of diagram \eqref{sympl_diag} sends
the element $1\o 1\in D\o D$ to $\omh\o1-1\o\omh\in  D\o D$.
The element $\omh$ projects to zero in $A$.
Therefore, 
applying the functor $A\o_D(-)\o_DA$ to  diagram \eqref{sympl_diag}
yields a diagram whose rows are  {\em complexes} of  $A$-bimodules.
Clearly, the  3 term complex in the bottom row of the resulting diagram
is obtained from   \eqref{FDG} by
composing the  map $\nu_0=p_{D|A}$ in  \eqref{FDG}
with the tautological imbedding $\jmath_A:\Om^1 A\into A\o  A$.

In more detail, the grading $D=\bigoplus_{r\geq 0} D_r$ induces
gradings on $\Om^1D$ and  on $\dder D$.
By condition \eqref{D}(i), we have that
 $\Om^1 D$ is a projective $D$-bimodule
generated by finitely many homogeneous elements
of degrees  $0,\ldots, n.$ It follows that
$\dder D$ is a $D$-bimodule
generated by finitely many homogeneous elements
of degrees  $-n,\ldots,0$. Thus, we have the grading
$\Om^1 D
=\bigoplus_{r\geq 0}\Om^1_rD$, and also 
the dual grading  $\dder D=\bigoplus_{r\geq -n}\dder_r D$.
 The symplectic 2-form $\om$ has degree $+n$, 
so the corresponding reduced contraction
map acts as  $\bi^\om: \dder_\idot{D}\to\Om^1_{\idot+n}D$.

Tensoring with $A$, we get $A$-bimodule
 direct sum decompositions
\beq{on}
\Om^1(D|A)=\boplus_{0\leq r\leq n}\Om^1_r(D|A),
\quad\text{resp.},\quad
\dder(D|A)=\boplus_{-n\leq r\leq 0}\dder_r(D|A),
\eeq
where each
direct summand 
is a finitely generated, projective $A$-bimodule.
Thus, applying the functor
$A\o_D(-)\o_DA$ to diagram ~\eqref{sympl_diag}, one obtains the
following commutative diagram
{{\small
\beq{ssdiag}
\xymatrix{
A\o A
\ar[r]^<>(0.5){\ad}
\ar@{=}[d]^<>(0.5){\Id}&\dder _0(D|A)\ar[r]^<>(0.5){\d_1^\vee}
\ar[d]^<>(0.5){\bi^\om}&
\ldots\ar[r]^<>(0.5){\d^\vee_n}&\dder_{-n}(D|A)\ar[rrr]^<>(0.5){\th\mto \th'(\omh)\o\th''(\omh)}
 \ar[d]^<>(0.5){\bi^\om}
&&&A\o A\ar@{=}[d]^<>(0.5){\Id}\\
A\o A\,\ar[r]^<>(0.5){^{1\o1\mto d\omh}}&
\,\Om^1 _n(D|A)\ar[r]^<>(0.5){\d_n}&
\ldots\ar[r]^<>(0.5){\d_1}&\Om^1 _0(D|A)
\ar[rrr]^<>(0.5){\jmath_{_A}\,\circ\, p_{_{D_0|A}}=i_\De}&&&A\o A.
}
\eeq}}
Each row of this diagram is a
complex of finitely generated, projective $A$-bimodules
and the vertical maps provide an isomorphism between the two complexes.

By Step 2, the complex  \eqref{FDG}, viewed
as a DG $A$-bimodule
with respect to the
total differential $\nu+\xi$ is acyclic. Therefore, from the short exact sequence
$\Om^1 A\into A\o  A\onto A$, we see that
the complex in the bottom row in \eqref{ssdiag}
has the only nonzero homology in degree zero, and this
homology group equals $A$.
Hence,  the bottom row in \eqref{ssdiag} provides a projective
$A$-bimodule resolution of
$A$. 

By Lemma \ref{lem1}(i), all the  vertical maps in diagram
\eqref{ssdiag} are $A$-bimodule isomorphisms, moreover,
these maps  commute with the Lie derivative $L_\xi$.
We deduce that the complex in the top row of  diagram
\eqref{ssdiag} is isomorphic to the complex in the bottom row.
Hence, the top row is a projective
$A$-bimodule resolution of
$A$ as well.

By definition,  applying the functor
$\Hom_{\bimof{A}}(-, A\o A)$ to the  projective $A$-bimodule resolution in
 the bottom row of  diagram \eqref{ssdiag}
yields a complex that represents the object $A^!=\RHom_{\bimod{A}}(A,
A\o A)$.
But 
applying $\Hom_{\bimof{A}}(-, A\o A)$ to the bottom row of 
the diagram gives, up to shift by
$d=n+2$, the top row of the same diagram, by
Lemma \ref{lem1}(ii).
We conclude that the vertical maps in  diagram \eqref{ssdiag}
provide a quasi-isomorphism $A^![d]\iso A$.
Moreover, this isomorphism is {\em selfdual}, since the map
$\bi^\om$ is selfadjoint,
due to Lemma \ref{FF}(ii).
Thus, we have constructed a  selfdual
quasi-isomorphism $A^![d]\iso A$, and $A$ is a  CY algebra of dimension $d$.
\qed

\subsection{Sketch of proof of Theorem \ref{main}(ii)}\label{proof_ii}
Let  $A$ be a friendly CY algebra of dimension $d=n+2$.
The case $n=0$ is somewhat special, so below we only consider
the case where $n > 0.$ We will
freely use the results of \S\ref{selfdual_sec}.

To prove the first statement in part  (ii) of the theorem, think of the quasi-isomorphism 
$f: A[d]\iso A^!,$ from Definition \ref{vdb_def},
as a morphism of resolutions given by the 
rows of  diagram \eqref{PPA}.
This  diagram is  {\em selfdual} by Theorem \ref{ddpf}. Therefore,
$f^!$, the dual morphism,
is  represented by the {\em same} diagram.
We deduce that $f=f^![d],$ hence
 $A$ is a CY algebra of dimension $d$.

The rest of the  proof  of part (ii) of Theorem \ref{main} proceeds in three steps.

\step{1.} Since $A$ is finitely presented,
one can find  a free and finitely generated  algebra
$B$ and a finitely generated two-sided ideal
$I'\sset  B$ such that $B/I'\cong A$.
Lemma \ref{imj}(iv) provides us with an $A$-bimodule resolution of the form 
\beq{complex1}
0\too A\o A\,\stackrel{\bj'}\too\,
P'_n\stackrel{d'_n}\too\ldots\stackrel{d'_2}\too P'_1\stackrel{\bp'}\too
I'/(I')^2\too0,
\eeq
where  all $P'_i$
are finitely generated projective $A$-bimodules.

Let $F$ be the $I'$-adic completion of $B$, and
let $I\sset F$ be the closure of $I'$. Thus,  we have $F/I=B/I'=A$
and, moreover, $F$ is
a smooth algebra.
Further, we have  $I/I^2=I'/(I')^2$, is a finitely
generated projective $A$-bimodule.

It is easy to prove the following

\begin{lem}\label{easy} There exist finitely
generated,  projective $F$-bimodules $P_r,\,r=1,\ldots,n,$
and a complex of $F$-bimodules
\beq{complex2}
\xymatrix{
F\o F\;\ar@{^{(}->}[r]^<>(0.5){\bj}&
P_n\ar[r]^<>(0.5){d_n}&
P_{n-1}\ar[r]^<>(0.5){d_{n-1}}&\ldots\ar[r]^<>(0.5){d_2}&P_1\ar@{->>}[r]^<>(0.5){\bp}&
I,
}
\eeq
such that, for any $r=1,\ldots,n,$ we have $P'_r=A\o_F P_r\o_FA,$
and, moreover, the functor $A\o_F(-)\o_F A$ takes the 
complex in \eqref{complex2} to \eqref{complex1}.\qed
\end{lem}

We may use the bimodules  from Lemma \ref{easy} to define a  graded
algebra
$$
D':=T_F\big(P_1\oplus P_2\oplus\ldots\oplus P_n\big)=
(T_FP_1)*_{_F}(T_FP_2)*_{_F}\ldots*_{_F}(T_FP_n).
$$
We equip this  algebra with a grading
$D'=\bigoplus_{r\geq 0} D'_r,$   such that  the bimodule $P_r$
is placed in degree $r$.  Thus, we have $D'_0=F$.
The  differential
in the complex  \eqref{complex2} can be 
uniquely extended to a  super-derivation 
$\xi: D'_\idot\to D'_{\idot-1}$ such that $\xi(F)=0$ and
$\xi^2=0$. 

Let $D=\bigoplus_{r\geq 0} D_r$ be the graded
algebra obtained by the $I$-adic completion of each
homogeneous component $D'_r,\, r=1,2,\ldots,$
of the graded algebra $D'$. The super-derivation
$\xi$ extends, by continuity, to a super-derivation on $D$.
Thus, we have made $D$ a DG algebra, with differential $\xi: D_\idot\to
D_{\idot-1}$.
Note that $D_0=F$.

We also introduce a larger graded algebra
$$
\FD:=D*_{_F}T_F(F\o F)\;\supset\;
D'*_{_F}T_F(F\o F)=
T_F\big(P_1\oplus P_2\oplus\ldots\oplus P_n \oplus (F\o F)\big),
$$
where the summand $F\o F\sset\FD$ is placed in degree
$n+1$ and elements of $D$ are assigned their natural degrees.
We may extend the differential $\xi$ on $D$ to
a differential $\d: \FD_\idot\to\FD_{\idot-1}$ using the
map $\bj: F\o F\to P_n$ from \eqref{complex2}.
It is clear that $D$ becomes a DG subalgebra in $\FD$.
Furthermore,  each of the algebras
$D',D,$ and $\FD$ is smooth, by \cite{CQ}, Proposition 5.3(3).

Consider a free product,  $F*  \C[t]$,  viewed as a graded
algebra such that $F$ is a degree zero subalgebra and $\deg t:=n+1.$
We observe  that the assignment $t\mto 1\o 1$ extends to
a graded $F$-algebra isomorphism
$F*  \C[t]\iso T_F(F\o F).$ 
 Using this isomorphism, 
 we get
\beq{frpr}
\FD=D*_{_F}T_F(F\o F)=D *_F \big(F*  \C[t]\big)=D*  \C[t],\qquad\deg t=n+1.
\eeq

From now on, we will identify the DG algebra $\FD$ with
$D*  \C[t].$
\begin{lem}\label{cyc0} 
 The composite  $\FD\onto\FD_0=F\onto F/I= A$ is a quasi-isomorphism.
\end{lem} 

\begin{proof} Using  \cite{CBEG}, Lemma 5.2.3, 
 one finds, cf.  \eqref{omt},
\beq{fits}
\Om^1\FD=\Om^1(F|\FD)\;\boplus\; \left(\bigoplus_{r=1}^n
\FD\o_FP_r\o_F\FD\right)\;\boplus\;   \FD\cd dt\cd \FD.
\eeq
It is clear that the differential in $\Om^1 (\FD|F)$
induced from the differential in the DG algebra $\FD$ goes,
via this identification, to the differential in the complex
\eqref{complex2}.

Applying the functor  $A\o_F(-)\o_F A$
to \eqref{fits} and using Lemma
\ref{easy} we obtain the following commutative diagram
of $A$-bimodules
\beq{gi4}
{\small\xymatrix{
A\o A\,
\ar@{^{(}->}[rr]^<>(0.5){1\o1\mto d\omh}\ar@{=}[d]^<>(0.5){\Id}&&\,
P'_n\ar[r]^<>(0.5){d'_n}\ar@{=}[d]&
\ldots\ar[r]^<>(0.5){d'_2}&P'_1\ar[r]^<>(0.5){d'_1}\ar@{=}[d]&\Om^1(F|A)\ar@{->>}[rr]^<>(0.5){p_{F|A}}\ar@{=}[d]&&\Om^1 A\ar@{=}[d]
\\
\Om^1_{n+1}(\FD|A)\,\ar@{^{(}->}[rr]^<>(0.5){\d_{n+1}}&&
\,\Om^1_n(\FD|A)\ar[r]^<>(0.5){\d_n}&
\ldots\ar[r]^<>(0.5){\d_1}&\Om^1_1(\FD|A)\ar[r]^<>(0.5){\d_1}&
\Om^1_0(\FD|A)\ar@{->>}[rr]^<>(0.5){p_{\FD_0|A}}&&\Om^1 A.
}}
\eeq

The complex in the top row of the above diagram is an exact sequence
by construction. Hence, the complex in the bottom row is  an exact
sequence
as well, and
we may apply
Proposition \ref{kunneth}(ii). The lemma follows.
\end{proof}

\step{2.} It is 
clear from the construction that  $\d(t)\in D_n$.
We put $\omh:=\d(t)\in D$.

\begin{lem}\label{cyc1}
 There exists a unique cyclic {\em closed}
2-form $\om\in \big(\Om^2 D\big)\br,$ which is homogeneous of degree
$n$ and  such that, in
$\Om^1D$, we have $d\omh=\bi_\De\om$.
Furthermore, we have $L_\xi\om=0$.
\end{lem} 

\begin{proof} 
We show first  that 
\beq{dom_der}
d\omh\in [D,\,\Om^1 D].
\eeq

To prove this, we consider a chain of isomorphisms
\beq{HHH}
H_{n+1}(\Om^1 \FD/[\FD,\Om^1 \FD];\,\d)
\stackrel{1}\cong H_{n+2}(\FD,\FD;\,
 b+\d)\stackrel{2}\cong H_{n+2}(A,A)\stackrel{3}\cong H^0(A,A)=\CZ(A).
\eeq
Here, $H_\idot(\FD,\FD;\, b+\d)$ denotes Hochschild homology of
$\FD$, the latter being viewed as a DG algebra, that is, the homology of
the Hochschild chain complex equipped with the differential
$b+\d$, the sum of the Hochschild  differential
$b$ and the  differential  $\d$ induced by the differential in the
DG algebra $\FD$ itself. 
The first isomorphism above is a standard result about
Hochschild homology of smooth algebras.
The second isomorphism follows from the quasi-isomorphism
$\FD\qis A$, and the third  isomorphism is due to \eqref{vdb}.

It is easy to see that the class in
$H_{n+1}(\Om^1 \FD/[\FD,\Om^1 \FD],\,\d)$
corresponding to the element $1\in\CZ(A)$ on the right
of \eqref{HHH} must be
represented
by the 1-form $dt$. It follows, in particular, that this 1-form 
must be  a {\em cycle}, i.e., that the equation $\d(dt)=0$
holds in $\Om^1 \FD/[\FD,\Om^1 \FD]$.
This implies $\d(dt)=d(\d(t))=d\omh\in [\FD,\Om^1 \FD]$.

Finally, we observe that $\d(t)\in P_n\sset D.$
Hence, we have $d\omh\in \Om^1 D$. 
Further, from \eqref{fits} we deduce
$$
\Om^1 \FD\big/[\FD,\,\Om^1 \FD]=(\FD\o_F\Om^1D)\big/[D,\,\FD\o_D\Om^1D]
\;\;\boplus\;\; \FD\cd dt.
$$
It follows easily
that $\Om^1 D\cap [\FD,\Om^1 \FD]=[D,\Om^1 D],$ and \eqref{dom_der}
is proved.

Next, we apply
\cite{Gi3}, Proposition 5.5.1, which yields

\begin{claim}\label{inject} For any  algebra $B$ such that $ H_2(B,B)=0$,
 the   map $\bi_\De$ below is {\em bijective}
\vskip 4pt

\noindent
$\dis\qquad\qquad\qquad
\xymatrix{
(\Om^2 B)\br^{\oper{closed}}\;\ar[rr]^<>(0.5){\bi_\De}_<>(0.5){\sim}&&\;\{\al\in [B,\Om^{1}
 B]\en\big|\en d\al=0\}.}$\hfill$\lozenge$
\end{claim}

The algebra $B:=D$ is smooth, and the 1-form
$\al=d\omh$ is clearly closed.  Hence, we have $H_r(D,D)=0,\,\forall r\geq 2,$
and Claim \ref{inject} combined with \eqref{dom_der} yields  the existence of a closed 2-form
$\om$ such that $\bi_\De\om=d\omh.$

Finally, since $\bi_\De$ commutes with $L_\xi$ and
anti-commutes with $d$, we compute
$$
\bi_\De(L_\xi\om)=L_\xi(\bi_\De\om)=L_\xi(d\omh)=-d(\xi(\omh))=
-d(\d^2(t))=0.
$$
Here, the map $\bi_\De$ that appears in the leftmost term is bijective,
by Claim \ref{inject}.
Hence, $\bi_\De(L_\xi\om)=0$ implies 
 $L_\xi\om=0$, and Lemma \ref{cyc1} follows.
\end{proof}

\step{3.} In the previous step,
we have constructed a data
$(D, \om,\xi)$ that satisfies all the conditions from
\eqref{D}(i)-(iii) except possibly the nondegeneracy condition, 
requiring the map $\bi^{\om}: \dder D\to\Om^1 D$ be bijective.
It is in effect unlikely that 
the  nondegeneracy condition may be achieved without  further modification
of the data $(D, \om,\xi)$. 

Nonetheless, we are going to prove 

\begin{lem}\label{hom} Let $(D, \om,\xi)$ be the data constructed above.
Then, the contraction map
$\bi^\om: \dder D\to\Om^1 D$ is a quasi-isomorphism.
\end{lem}
\begin{proof}[Sketch of Proof] Let $p:\FD\onto A$ be the natural algebra projection.
We consider the following commutative diagram
\beq{induced}
\xymatrix{
\dder (D|\FD)\ar[d]^<>(0.5){p\o\Id_{\dder D}\o p}
\ar[rrrr]^<>(0.5){\bi_\FD:=\Id_\FD\o\bi^\om\o\Id_\FD}&&&&\Om^1 (D|\FD)
\ar[d]^<>(0.5){p\o\Id_{\Om^1D}\o p}\\
\dder (D|A)
\ar[rrrr]^<>(0.5){\bi_A:=\Id_A\o\bi^\om\o\Id_A}&&&&\Om^1 (D|A).
}
\eeq

We claim that each of the vertical maps in the diagram is a
quasi-isomorphism. To prove this, view
$\dder D$  as a left DG $D^e$-module, cf. Notation \ref{e}.
Since $D$ is smooth, this DG $D^e$-modules is d-projective
in the sense of Definition \ref{sf}.
Hence, applying Lemma \ref{claim2} in the case
where $M=\dder D$ and $f$ is the quasi-isomorphism $p\o p: \FD^e\to A^e$,
we deduce that the vertical map on the
left of  diagram
\eqref{induced} is a quasi-isomorphism. Similar argument applies
to  the vertical map on the right of  diagram
\eqref{induced}, and our claim follows.

Next, we claim that the morphism $\bi_A$ in the bottom
line  of  diagram
\eqref{induced} is a quasi-isomorphism.  To prove this, we observe that
the morphism in question is part of diagram \eqref{ssdiag}.
Specifically, remove two vertical identity maps
from the sides  of diagram \eqref{ssdiag}. The collection
of the remaining  vertical  maps is then nothing but
a more detailed description of our map
$\bi_A$. Therefore, we see
that  proving that  our map is a quasi-isomorphism is equivalent to
proving that 
the  collection
of {\em all} vertical
maps in  \eqref{ssdiag}, including the two identity
maps that we have removed earlier, provides
 a quasi-isomorphism between the complexes
in the rows of that diagram.

Now, the complex in the bottom  row of  diagram \eqref{ssdiag}
may be identified with  the complex in the bottom  row of
 diagram \eqref{gi4}. The latter  complex is acyclic by Step 2.
Hence,  the complex in the bottom  row of  diagram \eqref{ssdiag}
provides a projective resolution of $A$, the diagonal
$A$-bimodule. The  complex in the top row of  diagram \eqref{ssdiag}
is obtained by applying the duality functor $(-)^\vee$ to
that resulution. The algebra $A$ is a CY algebra of
dimension $n+2$.
We conclude that the  complex in the top row of  
 \eqref{ssdiag} is also quasi-isomorphic to $A$. It follows
that the morphism between the rows  of  diagram \eqref{ssdiag}
is a quasi-isomorphism. This completes the proof of our claim
that the morphism $\bi_A$ in the bottom
line  of  diagram
\eqref{induced} is a quasi-isomorphism.

Thus, we have proved that 3 out of 4 maps in 
 diagram
\eqref{induced} are quasi-isomorphisms. It follows that the remaining 4th
map, the map $\bi_\FD$ in the top row, is a quasi-isomorphism as well.
Observe further that both $\dder (D|\FD)$ and $\Om^1 (D|\FD)$
are d-projective left DG $\FD^e$-modules, since
$D$ is smooth. Hence,  Lemma \ref{claim2}(ii)
 implies
that there exists
a quasi-inverse DG $\FD^e$-module morphism $h:
 \Om^1 (D|\FD)\to\dder (D|\FD)$
such that each of the maps $h\ccirc\bi_\FD$
and $\bi_\FD\ccirc h$ is homotopic to the identity.

We claim that the map $h$ takes the
subspace $\Om^1D\sset\Om^1 (D|\FD)$ into 
the subspace $\dder D\sset\dder (D|\FD)$.
To see this,  recall first that 
$\Om^1 D$ is a $D$-bimodule
generated by  homogeneous elements
of degrees  $0,\ldots, n,$ and
$\dder D$ is a $D$-bimodule
generated by  homogeneous elements
of degrees  $-n,\ldots,0$, cf. \eqref{on}. Furthermore, the  reduced contraction
map acts as follows $\bi^\om: \dder_\idot{D}\to\Om^1_{\idot+n}D$.
So, the  quasi-inverse map acts as
$h: \Om^1_\idot{D}\to\dder_{\idot-n}$.
Now, we have $\deg t=n+1.$ We conclude that, for any homogeneous 
generator $\al\in\Om^1_\idot{D}$, the element $h(\al)\in\dder_\idot{D}$
can involve neither $t$ nor $dt$. Thus, the map $h$
restricts to a morphism 
$h_D:\Om^1D\to\dder D$.

We know that there exists a homotopy
$\kappa: \Om^1_\idot(D|\FD)\to\Om^1_{\idot+1}(D|\FD)$ such
that we have $h\ccirc\bi_\FD-\Id_{\Om^1(D|\FD)}=[\d,\kappa]$,
where $\d$ stands for the differential
on $\Om^1_\idot(D|\FD)$. The map $\kappa$ has degree $+1$.
Therefore, we may repeat the argument of the preceeding paragraph
to conclude that, for any homogeneous 
generator $\al\in\Om^1_r{D}$,
 the element $\kappa(\al)\in\Om^1_{r+1}(D|\FD)$
contains neither $t$ nor $dt$,
{\em provided we have} $\deg\al=r<n$.
This argument breaks down, however, 
for generators of degree $n$.

To overcome this difficulty we refine the choice
of a quasi-inverse map $h$ as follows.
We adapt the proof of
Theorem \ref{ddpf} given in \S\ref{splice_sec}
and show that one can find
a {\em selfdual} map $h=h^\vee[d]: \Om^1 (D|\FD)\to\dder (D|\FD)$
which is a quasi-inverse of the (selfdual)
map $\bi_\FD:\dder (D|\FD)\to\Om^1 (D|\FD)$.
It is clear that, for $h=h^\vee[d]$,
the action of $h$ on elements of degree $n$ is completely determied
by the corresponding action on elements of lower degrees.
Using this, one proves that there exists a homotopy
 $\kappa$  that preserves the
subspace $\Om^1D\sset\Om^1 (D|\FD).$

This completes the proof of the lemma, hence, the proof
of Theorem \ref{main}(ii).
\end{proof}

{\large\section{{\textbf{Quiver algebras and McKay correspondence in dimension 3}}}\label{alg}}
\subsection{$R$-algebras}
It will be necessary for applications below
to replace the  ground field $\C$ by an arbitrary finite dimensional semisimple
$\C$-algebra $R$, as our ground ring.  An algebra $F$ equipped with
an algebra imbedding $R\into F$ will be referred to as
an $R$-algebra.

All standard constructions of  noncommutative  calculus have their 
{\em relative} analogues in the context of $R$-algebras.
Given
an $R$-algebra $F$, we let
$\Om^1_RF:=\ker(F\o_R F\to F)$ be the $F$-bimodule
of {\em relative} differentials,
and  $\Om^\hdot_RF:=T^\hdot_F(\Om^1_RF)$,
the DG algebra of relative differential forms.
In  $\Om^\hdot_RF,$ we have the de Rham differential $d$ that descends
to the super-commutator quotient $(\Om^\hdot_RF)\br=
\Om^\hdot  _RF/[\Om^\hdot  _RF,\Om^\hdot  _RF]_\text{super}.$

From now on, we let $R:=\k I$ be the commutative algebra of $\C$-valued functions on 
a finite set $I$, with pointwise
multiplication. For any $i\in I$, let $1_i\in\k I$ be the
characteristic function of the one-point set $\{i\}\sset I$.

Given  an  $R$-bimodule $M$, we let
$\dis M^R:=\{m\in M\en\big|\en r\cdot m=m\cdot
r,\,\forall r\in R\}$ denote the {\em center} of $M$.
 The element ${\mathbf e}:=\sum_{i\in I} 1_i\o 1_i$
belongs to  $(R\o R)^R$.

Let $F$ be an $R$-algebra.
The {\em inner} bimodule structure makes the space $(F\o F)^R$
an $F$-bimodule.
This bimodule is
generated by the element
${\mathbf{e}}\in (R\o R)^R\sset (F\o
F)^R$, cf. \cite{CQ}. There is
a natural  $F$-bimodule isomorphism  $(F\o F)^R\cong(F\o_RF)^\vee$, cf. \eqref{DUAL}.

Write $\dder_RF=(\Om^1_RF)^\vee$ for
the bimodule of {\em relative} double derivations, i.e.,
of derivations $\th: F\to F\o F,\,u \mto \th'u\o \th''u,$
such that $\th(R)=0.$ Also put $\Th^\hdot_RF=T_R^\hdot(\dder_RF).$

Observe  that, for $u\in F\o F$,
the inner derivation $\ad u: F\to F\o F$ annihilates
the subalgebra $R\sset F$
iff $u\in (F\o F)^R$. In particular, following \cite{CBEG},
\cite{VdB4},
we put $\De:=\ad{\mathbf{e}}$. 
Thus,
$\De\in\dder_RF$ is the
 relative $R$-counterpart  of \eqref{De}.
  
We write $\Lmof{F}$ for the category
of {\em finitely generated} left $F$-modules.
Let $A*_RB$ denote a free product (over $R$) of two $R$-algebras,
$A$ and $B$.

\subsection{Quiver algebras}\label{Q_sec}
Let $Q$ be a finite quiver with  vertex set $I$.
We let $\QQ$ be the 
double of $Q$ obtained by adding, for each edge
$x\in Q$, a reverse edge $x^*$. Further, let
 $\QQQ$  be the quiver obtained from $\QQ$ by attaching
an additional edge-loop, $t_i$, for every vertex $i\in I$.

One can  extend various constructions of noncommutative
calculus to the quiver setting. 
To this end,
write $\k Q$ for the path algebra of $Q$
and  view $R=\k I\sset \k Q$ as the subalgebra formed by trivial paths.
The commutator
quotient, $F\br$, of the algebra $F=\k Q$ may be identified with the
vector space spanned by {\em cyclic} paths in $Q$, cf. e.g.
[Gi1].

Generalizing formulas \eqref{dj} and  \eqref{partial}, 
 one defines,
 for any edge $x\in Q$, 
linear maps $F\br\to F\to F\o F,$ both
denoted $\pa/\pa x$ (the first map
has been introduced in  [Gi1] and [BLB1],
and the second map in [CBEG]). Thus, for any $\Phi\in  F\br,$
one has well-defined elements $\pa \Phi/\pa x\in F,\, x\in Q.$
This way, one can define an associative algebra
\beq{AQ}
\A(\k Q,\Phi):=\k Q/\llb \pa \Phi/\pa x\rrb_{x\in Q},
\eeq
the  quotient of the path algebra  by the two-sided
ideal generated by the collection of elements
$\pa \Phi/\pa x,\,x\in Q.$

\begin{rem}
In the special case of a quiver
with one vertex and $n$ edge-loops, we have
$R=\C,\,\C Q=\C\langle\xxx\rangle$, and formula \eqref{AQ} reduces
to \eqref{AF}.
\erem

For any quiver $Q$, there are natural algebra isomorphisms
\begin{align}\label{qqq}
&\Th^\hdot_R(\k Q)\,\iso\,\k\QQ, &x\mto x,\,\pa/\pa x\mto x^*,\quad\forall
x\in Q,\quad\text{resp.},\\
&(\Th^\hdot_R(\k Q))*_R R[t]\,\iso\,\k\QQQ, 
&\text{as above, and}\quad t\mto \sum\nolimits_{i\in I} t_i.\quad\nonumber
\end{align}
We remark that the isomorphism in the top line
sends the double derivation
$\De\in \Th^1_R(\k Q)$ to the element $\sum_{x\in Q} \,[x,x^*]\in\k\QQ$.
The isomorphism in the second line sends the element
$1_i\cdot t\cdot 1_i$ to  $t_i$, the  additional edge-loop  at the vertex $i\in I$.

\begin{examp}[Conifold algebra] Let $Q$ be a quiver with two vertices
and two edges, $x_i,\, i=1,2$, both directed in the same way.
Thus, $\QQ$, the double of $Q$, has two additional reverse edges,
$x^*_i,\, i=1,2$. The conifold algebra is defined as
$A=\A(\k\QQ, \Phi)$, where $\Phi=x_1x^*_1x_2x^*_2-x^*_1x_1x^*_2x_2.$

The center of the algebra $A=\A(\k\QQ, \Phi)$ is generated by 4
 elements which
satisfy one quadratic relation. Thus, the spectrum of
the centrer is a quadratic cone in $\k^4.$
\end{examp}

\begin{examp}[Cyclic quiver] Let $Q$ be 
a quiver with $n+1$ vertices and edges $x_i: i \to i+1\;\op{mod}\,(n+1)$,
which form an oriented
cycle of length $n+1$. Let $\Phi$ be the cycle, $\Phi:=x_1\,x_2\ldots
x_n\,x_{n+1}\;\op{mod}\,[\C Q,\C Q]$.
Then, $A=\A(\k{Q}, \Phi)$ is a quotient of $\k{Q}$ by the two-sided ideal
generated
by all paths of length $\ge n$. 
\end{examp}

\begin{rem} It is not difficult to show, using for instance Theorem \ref{CYconditions}
from \S\ref{cot} below, that the algebras
$A$ arising from
the above examples are Calabi-Yau algebras of dimension ~3.
\erem

For any quiver $Q$ and any potential $\Phi\in (\k Q)\br$, there is a
 quiver version of the DG algebra \eqref{Dalg}. 
It is defined to be the  algebra $\FD(\k Q, \Phi):=\k\QQQ,$
cf. \eqref{qqq}, equipped with a grading
such that $\deg t=2,$ and $\deg x=0,$ resp. $\deg x^*=1,$ for
any edge $x\in Q$. 

The differential $\d: 
\FD_\idot(F,\Phi)\to\FD_{\idot-1}(F,\Phi)$ is defined on generators by the assignment
$$\d:\
t\mto \De=\mbox{$\sum_{x\in Q}$}
\,[x,x^*],
\quad x^*\mto \pa\Phi/\pa x,\quad 
x\mto 0\quad\forall x\in Q.
$$

\subsection{$N=1$ super-symmetric
quiver gauge theory} An important class of examples is obtained from
a `nongraded'  analogue of the algebra
$\FD(\k Q, \Phi)$ as follows, see \cite{CKV}.

Fix $(Q,I)$, as above,  and view
$\k \QQQ= \k \QQ *_R R[t],$ cf. \eqref{qqq}, 
 as an algebra {\em without} grading.
Given an element $f\in R[t]$, in \cite{CKV},
the authors introduce the
following potential, cf. ~\eqref{RC2},
\beq{ga}
\Phi_f=t\cdot\sum_{x\in Q}
\,[x,x^*]\; +\;f(t)\,\in\, (\k\QQQ)\br.
\eeq

For the corresponding algebra $\A$, one finds
\beq{contact_quiv}
\A(\k\QQQ,\Phi_f)=
\frac{\k\QQQ}{
\left(\!\!\left(\begin{array}{c}
\mbox{$\sum_{x\in Q}$}\,[x,x^*]=\frac{df}{dt},\\
t\en\text{is a central
element}\end{array}
\right)\!\!\right)}.
\eeq

To get a better understanding of the RHS in this formula, observe that
one has obvious algebra isomorphisms
$R[t]=(\oplus_{i\in I}\C 1_i)\bigotimes\C[t]\iso\bigoplus_{i\in I}\k[t_i],\,
1_i\o t\mto t_i$.
Therefore, choosing an element $f\in R[t]$ amounts to
choosing  an $I$-tuple of polynomials $\{f_i\in\k[t_i]\}_{i\in I}$
with $\C$-coefficients. With this understood, \eqref{contact_quiv} reads
$$
\A(\k\QQQ,\Phi_f) =
\frac{\k\QQQ}{
\left(\!\!\left(\begin{array}{c}
\sum_{h(x)=i}xx^*-\sum_{t(x)=i}x^*x=\frac{\pa f_i}{\pa t_i},\;i\in
I,\\
t_ix=xt_j\en\text{for}\; x: j\to i.
\end{array}
\right)\!\!\right)}.
$$

\begin{examp}[Special case: $f=0$] In this case  we get
$\Phi_f=t\cdot\sum_{x\in Q}
\,[x,x^*]$.
Notice that this
the potential 
is a generalization of the potential
\eqref{abc} for $c=0$.
The corresponding  quiver generalization of the algebra
isomorphism in \eqref{bas_iso} reads
\beq{Pieq}
\A(\k\QQQ,\,\Phi_0)\,\cong
\,\C\QQ/\llb \De\rrb\;\o_R\; R[t]=\Pi(Q)\o_R R[t],
\eeq
where $\Pi(Q):=\C\QQ/\llb\sum_{x\in Q}[x,x^*]\rrb$
stands for the
{\em preprojective algebra} of $Q$. \hfill$\lozenge$
\end{examp}

\begin{examp}
More generally, given an $I$-tuple of constants
$\mu_i\in \k,\,i\in I,$
put $f=\frac{1}{2}\sum_{i\in I}\mu_i t_i^2$. We get the following 
'centrally extended' preprojective algebra,
that has been  considered earlier in \cite{ER}:
\vskip 4pt
$\qquad\qquad\dis\A\Big(\k\QQQ,\, t\cd\mbox{$\sum_{x\in Q}$}
\,[x,x^*]+\mbox{$\frac{1}{2}\sum \mu_it_i^2$}\Big)\cong\left.\frac{\k\QQQ}{
\left(\!\!\left(\begin{array}{c}
\mbox{$\sum_{x\in Q}$}\,[x,x^*]=\mbox{$\sum_{i\in I}$}\,\mu_i t_i,\\
\mbox{$\sum$}t_i\en\text{is a central
element}\end{array}
\right)\!\!\right)}\right._{\hskip 10mm\lozenge}$
\end{examp}

In the special case where $Q$ is an extended Dynkin quiver,
the algebra $\A(\k\QQQ,\Phi_f)$ has a geometric origin.
To explain this, recall that the vertex set $I$ of
an extended Dynkin quiver can be identified canonically
with the set of simple roots of the corresponding affine
root system, cf. eg.
\cite{CBH}. This way, the vector space $R=\C I$ becomes
identified with the affine Cartan algebra. There is a
vector $\delta_\text{im}\in R$, the minimal imaginary root,
whose coordinates are given by  Dynkin labels.
Let ${\mathfrak{h}}\sset R$ be the codimension 1
subspace orthogonal to  $\delta_\text{im}$.
This subspace may be identified with the Cartan algebra
of the  finite root system
associated to the Dynkin diagram obtained by removing the extending
vertex.

Now, let $\X_0\to X_0$ be the minimal resolution
of  the Kleinian singularity associated with the
Dynkin diagram. Thus, $\X_0$ is a smooth 2-dimensional variety that
has a standard universal deformation
$p:\X\to{\mathfrak{h}}$ parametrized by the vector space
${\mathfrak{h}}$, the base of the deformation.
Here, $\X$ is a smooth variety and $p$ is a flat morphism.

To any  polynomial
map $f: \C\to {\mathfrak{h}}$ one associates a 
smooth 3-dimensional noncompact  CY manifold
 $\X_f$.  It is defined via
the following base change diagram
$$
\xymatrix{
\X_f\ar[rr]\ar[d]^<>(0.5){p_f}&&\X\ar[d]^<>(0.5){p}\\
\C\ar[rr]^<>(0.5){f}&&{\mathfrak{h}}.
}
$$

By construction, $\X_f$ is fibered over $\C$ 
 with smooth 2-dimensional fibers.
One can apply the  equivalence of categories,
due to Kapranov-Vasserot \cite{KV},
along the fibers of the projection $p_f$, cf.
~\cite{Sz}, \cite{Z}.
This way, one establishes an equivalence
between the $D^b\big(\Lmof{\A(\C\QQQ,\Phi_f)}\big)$ and an appropriate
derived category of coherent sheaves on $\X_f$.
Here, we have identified the polynomial map $f: \C\to {\mathfrak{h}}$
with an  element  $f\in R[t]$ via the imbedding
${\mathfrak{h}}\into R$, and write $\Phi_f$ for the
corresponding potential \eqref{ga}.

\subsection{McKay correspondence in dimension 3.}\label{mckay}
Let $V=\C^3$ be a 3-dimensional vector space.
We  fix a nonzero
volume element $\pi\in \La^3V^*$ and  a finite
subgroup $\G\sset SL(V,\pi)$.

\begin{rem} The interested reader may consult e.g. \cite{FFK} 
for a complete classification of finite subgroups
in $SL_3(\C).$
\erem

Associated to $\G$ is its {\em MacKay quiver}
$Q=Q_\G$ whose vertex set $I$ is the set of
(isomorphism classes of) irreducible representations of $\G$.
Write $L_i$ for the  irreducible representation with label $i\in I$.
By definition, for any two vertices $i,j\in I,$
we have
\beq{edges}\#\{\text{edges}: i\to j\}=\dim\Hom_\G(L_i,L_j\o V)=
\dim\Hom_\G(V^*\o L_i,L_j),
\eeq
where $V^*$ is the dual representation.
So, the set  of edges $i\to j$ may (and will) be identified 
with a $\k$-linear basis of the vector space $\Hom_\G(L_i,L_j\o V).$

To any oriented triangle $\btr$ in $Q_\G$ one may
canonically attach a complex
number $\la(\btr)\in\C$ as follows.
Let $i,j,k\in I$ be the vertices of $\btr$ and
$x_{ij},\, x_{jk},\,x_{ki}$ the corresponding edges,
 the sides of the triangle.
Thus, we have three $\G$-intertwiners
$$x_{ij}: L_i\to L_j\o V,\quad x_{jk}: L_j\to L_k\o V,\quad
x_{ki}: L_k\to L_i\o V.
$$

We consider the composite intertwiner
\beq{intert}
\xymatrix{
L_i\ar[r]^<>(0.5){x_{ij}}&
L_j\o V\ar[rr]^<>(0.5){x_{jk}\o\Id_V}&&
L_k\o V^{\o 2}\ar[rr]^<>(0.5){x_{ki}\o\Id_V^{\,\o 2}}
&& L_i\o V^{\o 3}\ar@{->>}[rr]^<>(0.5){\Id_{L_i}\o \pi}&&
L_i.
}
\eeq
Here, we view  $\La^3V^*$ as being a subspace in $(V^{\o 3})^*$,
so that the volume $\pi$ becomes a map $V^{\o 3}\onto\C$.

The composite intertwiner in \eqref{intert}
is a scalar operator $x_{ijk}: L_i\to L_i$, by  Schur's lemma. Furthermore,
one verifies easily that
 the trace of the composite
is independent of the choice of the source-vertex, $i$, so one may
  put  $\la(\btr):=\Tr(x_{ijk})\in\C$.

Let $F=\k Q_\G$ be the path algebra of $Q_\G$. We observe  that,
for any oriented triangle in $Q$ as above, the  cyclic path
$x_{ij} x_{jk}x_{ki}$ gives a well-defined element
$\btr\in F\br$. We introduce the following  cubic potential
\beq{FG}
 \Phi_\G=\sum_{\{\text{oriented triangles}\;\btr\in Q\}}\,
\la(\btr)\cd\btr\in F\br.
\eeq

Thus, $\A\big(\k Q_\G,\, \Phi_\G\big)$ is a quadratic algebra which is,
moreover,
a CY algebra of dimension 3, by
Theorem \ref{mca} below. It follows  that   $\A\big(\k Q_\G,\,
\Phi_\G\big)$ is a
{\em Koszul algebra}, by \cite{BM}, \cite{BT}.

\begin{rem} \vi The corresondence between edges of the quiver
$Q_\G$ and  intertwiners $x_{ij}$
depends on a choice of basis in the vector space
$\Hom_\G(L_i, L_j\o V)$. Therefore, the
coefficients $\la(\btr)$ and the potential $\Phi_\G$ depend
on various choices. A more careful analysis shows
that the resulting algebra $\A(\C Q_\G,\Phi_\G)$ is,
nevertheless, independent
of those choices, up to isomorphism.

\vii Oriented triangles in  McKay quivers have first appeared
in the work of A.~ Ocneanu; they are sometimes referred to
as {\em Ocneanu cells}, cf. e.g. \cite{CT}.
\erem

Recal the idempotent
 $1_i\in\k I\sset \k Q$ corresponding to a vertex $i\in I$.
In particular, one has the  idempotent $1_o$, where
$o\in I$ is the distinguished vertex
associated with the trivial representation, $L_o$, of $\G$.

 The following result which
 is an $SL_3(\C)$-generalization of a similar
result for finite subgroups $\G\sset SL_2(\C)$, proved in
\cite{CBH}.

\begin{thm}\label{mca} \vi The algebra $\A\big(\k Q_\G,\, \Phi_\G\big)$
is a CY algebra of dimension 3 and there is a Morita equivalence 
$$\C[V]\#\G\;\stackrel{\text{Morita}}\eqq\;\A\big(\k Q_\G,\, \Phi_\G\big).$$

\vii There is a canonical algebra isomorphism
$$\C[V]^\G\;\cong\; 1_o\cd\A\big(\k Q_\G,\, \Phi_\G\big)\cd 1_o.
$$

\viii For any $i\in I,$ there is a canonical isomorphism
of right modules 
$$
\Hom_{\G}(L_i,\C[V])\;\cong\; 1_i\cd\A\big(\k Q_\G,\,
\Phi_\G\big)\cd 1_o.
$$
\end{thm}
\noindent
Here, the object on the LHS of the last isomorphism is viewed as a right
$\C[V]^\G$, the object on the RHS  is viewed as a right
$1_o\cd\A\big(\k Q_\G,\, \Phi_\G\big)\cd 1_o$-module,
and the module structures on the two sides of the
isomorphism   of part (iii) are compared
via the algebra isomorphism of part (ii).

\begin{rem} N. Nekrasov has kindly informed me that part (i)
of the theorem was, in some form,  known to physicists.
\erem

Using the result of Bridgeland-King-Reid \cite{BKR}, from Theorem
\ref{mca}, one obtains part (ii) of the following corollary
(part (i)  of the corollary is immediate):

\begin{cor} \vi There is an equivalence of {\em abelian}
categories
$$\coh V/\G \cong \Lmof{\big(1_o\cd\A(\k Q_\G,\, \Phi_\G)\cd 1_o\big)};$$

\vii For any smooth crepant resolution $X\onto V/G$,
there is an equivalence of {\em triangulated} categories
$$D^b(\coh X)\cong D^b\big(\Lmof{\A(\k Q_\G,\, \Phi_\G)}\big).$$
\end{cor}

\subsection{Proof of Theorem \ref{mca}}
The proof is based on the technique of
tensor categories that has been applied earlier
in \cite{MOV} in a similar situation.

\step{1.} Let  $\Lmof{SL_3}$ be the tensor
category of finite dimensional rational representations of the
group $SL_3=SL_3(V,\pi)$.
The trivial representation, $ L_o$, is
the unit of the tensor structure. 
Further, let $W:=V^*$, so we have
 $\pi\in \La^3(V^*)=\La^3W\sset  W^{\o 3}$.
Thus, the volume element  determines a morphism
$\pi:  L_o\to W^{\o 3}$.

We consider
the following  canonical 
$SL_3$-module morphisms
$$
V\iso V\o L_o\stackrel{\id\o\pi }\tooo V\o W^{\o 3}=
(V\o V^*)\o W^{\o 2}\stackrel{tr}\too L_o\o W^{\o 2}= W^{\o 2}.
$$
The composite morphism yields an
isomorphism 
$\gamma:\ W^*=V\iso \La^2 W\sset W^{\o 2}.$

We may also view the volume element 
$\pi\in W^{\o 3}$ as a cubic potential for the tensor
algebra $TW\cong\C\langle x,y,z\rangle,$ see
Example \ref{bas_iso}. The map $\gamma$ then takes the form
of a `directional derivative'
$\gamma: W^*\ni u\mto \frac{\pa\pi }{\pa u}\in TW$.

With this notation, a coordinate free version of \eqref{bas_iso}
reads
\beq{sym_inv}
\Sym W=TW/\llb\La^2 W\rrb=TW/\llb\gamma(W^*)\rrb.
\eeq
Here, we treat the tensor algebra $TW=\C\oplus W \oplus W^{\o
2}\oplus\ldots,$ as well as each side of the isomorphism in \eqref{sym_inv},
 as a $\Z$-graded {\em  algebra ind-object} in the tensor
category $\Lmof{SL_3}$,
and where $\llb\La^2 W\rrb=\llb\gamma(W^*)\rrb$ denotes
the tensor ideal  generated by
the $SL_3$-subrepresentation of  skew 2-tensors.

\step{2.} Let $I$ be a finite set, that will later on become the set
of (isomorphism classes of) simple $\G$-modules.
The category 
$\scc_I:=\Lmof{\C I}$ is   a semisimple
abelian category which is clearly equivalent to the category of finite
dimensional $I$-graded vector spaces. For each $i\in I$, write
 $\C_i$ for the corresponding  1-dimensional 
 simple object of $\scc_I$
concentrated in grade $i$. Similarly, one has
the category $\bimof{\C I}$ with simple objects
$\C_{i,j},\, (i,j)\in I\times I.$
The latter category  has a natural monoidal
structure given by
 convolution;  the monoidal
structure is given by the formula
$E\o E'=\bigoplus_{i,k\in I} (E\o E')_{i,k},$ where
$$
(E\o E')_{i,k}:=\bigoplus_{j\in I} (E_{i,j}\o E'_{j,k}),
\quad\forall
E=\boplus_{i,j\in I}E_{i,j},
\;E'=\boplus_{i,j\in I}E'_{i,j}\in \bimof{\C I}.
$$

Following \cite{MOV}, we also consider the category
$Fun(\scc_I,\scc_I)$, of all additive functors $\scc_I\to\scc_I$.
The category $Fun(\scc_I,\scc_I)$ has a natural monoidal
structure given by composition of functors,
and the category $\scc_I$ is, tautologically,
a module category over the monoidal
category $Fun(\scc_I,\scc_I)$.
Furthermore, we have an obvious equivalence
of monoidal categories
$Fun(\scc_I,\scc_I)\cong \bimof{\C I}.$

\step{3.}
Let $\Lmof{\G}$ be the abelian category
of finite dimensional representations of our finite
subgroup $\G\sset SL(V,\pi)$, let $I$ be the set of 
(isomorphism classes of) its simple objects.
Let $\C\G\cong \oplus_{i\in I} \End_\C L_i$ be the group algebra of $\G$.
For each $i\in I$, we choose and fix $p_i\in\C\G$, a 
 minimal idempotent in the simple 
direct summand $\End_\C L_i$. 
Thus, we have $L_i=\C\G\cdot p_i$. Set $p:=\sum_i p_i.$ 
Then $p$ is an idempotent, and we have
an algebra isomorphism $p(\C\G)p\cong \C I.$ Furthermore, the
following functor is a Morita 
equivalence 
$$
\Lmof{\G}\too \Lmof{\C I}=\scc_I,\quad
M\mto \oplus_{i\in I}\Hom_{\G}(L_i,M)= p\cd M.
$$
We deduce  the following equivalences of the corresponding monoidal categories 
\beq{ttt}
Fun(\Lmof{\G},\Lmof{\G})\cong
Fun(\scc_I,\scc_I)\cong \bimof{\C I}.
\eeq

We also have the category $\bimof{\G}$, of finite dimensional
$\C\G$-bimodules, which
is a monoidal category with respect to the tensor product
of bimodules.
One has an equivalences of monoidal categories 
\beq{psi}
{\mathsf G}:\ \bimof{\G}\iso \bimof{\C I},\quad U\mto 
{\mathsf G}(U)=p\cd U\cd p=\oplus_{i,j} p_i \cd U\cd  p_j.
\eeq

The category  $\Lmof{\G}$ is naturally a module
category over  the monoidal category $\bimof{\G}$. Specifically,
any object $U\in  \bimof{\G}$ gives
an additive functor $M\mto U\o_{\C\G}M,$
on $\Lmof{\G}$.
This yields an equivalence
\beq{TT}: \bimof{\G}\iso Fun(\Lmof{\G}, \Lmof{\G}),
\quad U\mto U\o_{\C\G}(-).
\eeq

The category $\Lmof{\G}$ is also a module
category over  the tensor
category $\Lmof{SL_3},$ in other words,
tensoring with an $SL_3$-representation  gives
a functor $\Lmof{\G}\to\Lmof{\G}$. This way, we get
a monoidal functor ${\mathsf P}:\Lmof{SL_3}\to Fun(\Lmof{\G},\Lmof{\G}).$

Further, to any  $U\in \Lmof{SL_3}$
one can associate a bimodule $U\#\G\in \bimod{\C\G}$;
it is defined to be the vector space $U\o \C\G$ equipped
with the action $\gamma(u\o g)\gamma':= \gamma(u)\o \gamma g\gamma',$
for any 
$\gamma,\gamma',g\in \G.$

This way, we obtain the following commutative diagram 
of monoidal functors, in which ${\mathsf F}$
stands for the composite of the functor ${\mathsf P}$ with the
equivalences
in \eqref{ttt},
\beq{diag_fun}
\xymatrix{
\Lmof{SL_3}\ar[d]_<>(0.5){U\mto U\#\G}
\ar[drr]^<>(0.5){{\mathsf P}}\ar@/^1pc/[drrrr]|-{\;{\mathsf F}\;}&&&&\\
\bimof{\G}\ar@{=}[rr]^<>(0.5){\eqref{TT}}
\ar@/_1pc/[rrrr]|-{\;{\mathsf G}\;}&&
 Fun(\Lmof{\G},\Lmof{\G}) \ar@{=}[rr]^<>(0.5){\eqref{ttt}}
&&\bimof{\C I}.
}
\eeq

\step{4.} 
Let $Q=Q_\G$ be the McKay quiver for $\G$. 
One can view any graded $\C I$-algebra, e.g., the path algebra $\C Q$
or its quotient $\A(\C Q,\Phi_\G)$, as a graded algebra-ind-object
in the monoidal
category $\bimof{\C I}$.
 We claim that there is an
isomorphism of algebra-ind-objects
\beq{symw}
{\mathsf F}(\Sym W)\cong\C Q/\llb\pa\Phi_\G/\pa
x\rrb_{\{x\in Q\}}= \A(\C Q,\Phi_\G),
\eeq
where $\Sym W=L_o\oplus W \oplus \Sym^2W\oplus\ldots,$ is viewed
 as a $\Z$-graded algebra-ind-object in  $\Lmof{SL_3}.$

To prove this, we first notice that ${\mathsf F}(L_o)=\C I,$ the diagonal
$\C I$-bimodule. Next, we 
apply the functor ${\mathsf F}$ to the 
object $W\in \Lmof{SL_3}$.
Using the definition of the McKay quiver, cf. \eqref{edges},
we find ${\mathsf F}(W)=\bigoplus_{\{a: i\to j\;|\; a\in Q\}} \C_{i,j},$
where $\C_{i,j}$ denotes a 1-dimensional space placed
in grade $(i,j)\in I\times I$.

Now, view the tensor algebra $TW=\C\oplus W \oplus W^{\o
2}\oplus\ldots,$
 as a $\Z$-graded algebra ind-object in  $\Lmof{SL_3}.$
Following \cite{MOV},
we deduce that
${\mathsf F}(TW)=\C Q,$ the path algebra of $Q$ viewed as 
 a graded algebra ind-object in  $\bimof{\C I}.$

Further, it follows from definitions 
that the  composite morphism $\pi: L_o\iso\La^3W\into W^{\o 3}\into TW$
goes, under the functor ${\mathsf F}$, to the map
${\mathsf F}(\pi): \C I\to \C Q$ that sends the unit
$1 \in \C I$ to the sum in the RHS of formula
\eqref{FG}, that is, to the potential $\Phi_\G$. Therefore, the  functor
${\mathsf F}$
sends the morphism
$\gamma: W^*\to \La^2W\into TW,\, u\mto \pa \pi/\pa u$
to the morphism ${\mathsf F}(\gamma): x\mto \pa \Phi_\G/\pa x$.

Recall that ${\mathsf F}$, being a tensor functor, is automatically exact.
We conclude that
$${\mathsf F}(\Sym W)={\mathsf F}\left(\frac{TW}{\llb\gamma(W^*)\rrb}\right)=
\frac{{\mathsf F}(TW)}{\llb{\mathsf F}(\gamma(W^*))\rrb}
=\frac{\C Q}{\llb\im{\mathsf F}(\gamma)\rrb}=
\frac{\C Q}{\llb\pa\Phi_\G/\pa
x\rrb}_{\{x\in Q\}},
$$
and \eqref{symw} follows.
From this,
using \eqref{psi} and commutativity of diagram \eqref{diag_fun},
we get
$$ p((\Sym W)\#\G)p={\mathsf G}((\Sym W)\#\G)={\mathsf F}(\Sym W)=\A(\C Q,\Phi_\G).$$
The Morita equivalence from part (i) of Theorem \ref{mca}
follows.

Next, for 
 any $i,j\in I$, by definition, we deduce
$\dis
{\mathsf G}(L_i\boxtimes L_j)=p(L_i\boxtimes L_j)p=\C_{i,j}.
$
Thus, 
\begin{align*}
\Hom_{\G}(L_i,\Sym W)
&=\Hom_{\bimof{\G}}(L_i\boxtimes L_o,\;(\Sym W)\#\G)\\
&\iso\Hom_{\bimof{\C I}}({\mathsf G}(L_i\boxtimes L_o),\;{\mathsf G}(\Sym W))\\
&=\Hom_{\bimof{\C I}}(\C_{i,o},\, \A(\C Q,\Phi_\G))=1_i\cd\A(\C Q,\Phi_\G) \cd 1_o.
\end{align*}
 In particular, for $i=o$, we get
$(\Sym W)^\G=1_o\cdot \A(\C Q,\Phi_\G) \cdot 1_o.$
The isomorphisms of
 Theorem \ref{mca}(ii)-(iii) follow.
\qed

\begin{rem} The proof above applies without modification
in the more general setting where the tensor category $\Lmof{SL_3}$
is replaced by its $q$-analogue $\Lmof{\U_q({\mathfrak{s}\mathfrak{l}}_3)}$,
the category of finite dimensional representations
of the quantized enveloping algebra $\U_q({\mathfrak{sl}}_3)$,
where $q\in\C^\times$. This is a semisimple  braided monoidal category,
provided $q$ is not a root of unity.

Module categories over  $\Lmof{\U_q({\mathfrak{s}\mathfrak{l}}_3)}$ play the
role of $q$-analogues of categories $\Lmof{\G}$
arising from finite
subgroups $\G\sset\SL_3$. The module categories over  $\Lmof{\U_q({\mathfrak{s}\mathfrak{l}}_3)}$ have
been extensively studied, cf. eg. \cite{O} and references therein.
 A quantum analogue of the
McKay correspondence associates to any such  module category
$\scc_q$ a certain  data involving, in particular, the
corresponding McKay  quiver $Q(\scc_q)$, and a collection
of complex numbers $\la_q(\btr),$ one 
for each oriented triangle $\btr\in Q(\scc_q)$.
Thus, one can define a $q$-deformed potential
$\Phi_q:=\sum\la_q(\btr)\cd\btr.$

In the special case where the category
$\scc_q$ is a deformation of the
category  $\Lmof{\G}$ for a finite
subgroup $\G\sset\SL_3$
the corresponding  algebra
$\A(Q,\Phi_q)$ provides a CY deformation of
the algebra $\A(Q,\Phi_\G)$.

It would be very interesting to analyze these $q$-deformed
CY algebras in more detail.
\end{rem}

{\large\section{{\textbf{Calabi-Yau algebras of dimension  3}}}\label{cy3}}
\subsection{The algebra $\A(F,\al)$.}\label{algA}
Fix a finite dimensional,
semisimple algebra $R$ and an  arbitrary
$R$-algebra $F$.
Given    a cyclic 1-form
$\al\in(\Om^1_RF)\br$, we consider $\bi^\al:\dder_R F\to F,$
 the reduced contraction map. This is a morphism
of $F$-bimodules, hence, the image of that morphism  is a two-sided
ideal  $\imm\sset F$.

\begin{defn}\label{defal}
For any {\em smooth} $R$-algebra $F$ and a cyclic 1-form $\al$ such that 
\beq{cond}
d\al=0\qquad\text{and}\qquad i_{\De}\al =0,
\eeq
holds, we define an associative algebra
$\A(F,\al):=F/\imm.$
\end{defn}

Recall that, for a smooth algebra  $F$, all the schemes
$\Rep_\dd F,\,\dd=1,2,\ldots,$ are smooth. 
From the definition of the algebra
$\A(F,\al)$, we see that
$$ 
\Rep_\dd\A(F,\al)=\textsf{Zero locus  of the 1-form}\en
\Tr\wh\al\quad\sset\; \Rep_\dd F.
$$

Observe that conditions \eqref{cond} hold {\em automatically} 
for any {\em exact} 1-form $\al=d\Phi,$ where $\Phi\in F\br.$ 
Abusing the notation, in such a case we will write $\A(F,\Phi)$ instead
of $\A(F,d\Phi)$. This is justified by the following simple result.

\begin{lem} Let $Q$ be a quiver with vertex set $I$ and put
$R=\k I$. 

For any  potential  $\Phi\in(\C Q)\br$, the algebra
$\A(\C Q, d\Phi)$, arising from Definition \ref{defal} for $\al=d\Phi$,
is isomorphic to the algebra
$\A(\C Q,\Phi)=\C Q/\llb \pa \Phi/\pa x\rrb_{x\in Q},$
from \eqref{AQ}.\qed
\end{lem}
\vskip 4pt

Fix an algebra $F$, a 1-form $\al$ and put $A=\A(F,\al)$. We apply the functor
$A\o_F (-)\o_F A$ to the $F$-bimodule map $\bi^\al$.
Clearly, one has
$A\o_F (\imm)\o_F A=\imm/(\imm)^2$.
Thus, we obtain
the following   sequence of maps
\beq{map}
d\ccirc\bi:\
\xymatrix{
 \dder_R(F|A)\en\ar@{->>}[rr]^<>(0.5){\bi^\al}&&
\en \imm/(\imm)^2\en \ar@{^{(}->}[rr]^<>(0.5){d_{F|A}}&&\en
\Om^1_R(F|A).}
\eeq
Injectivity of the map
 $d_{F|A}$  follows from the   {\em cotangent  sequence} \eqref{CQ}.

Further, we splice the cotangent sequence \eqref{CQ}
with the tautological short exact sequence
$\Om^1_RA\into A\o_R A\onto A.$
This way, using  \eqref{map}, we obtain a sequence of maps
\beq{ext3}
0\too\xymatrix{
( A\o  A)^R\ar[r]^<>(0.5){\ad}&
 \dder_R(F| A)\ar[r]^<>(0.5){d\circ\bi}&
\Om^1_R(F| A) \ar[r]^<>(0.5){p_{F| A}}&
 A\o_R  A\ar[r]^<>(0.5){\op{mult}}& A
}\too 0.
\eeq

\begin{lem}\label{id_cond} For any cyclic 1-form $\al$ that
satisfies the equation $i_{\De}\al =0$ from \eqref{cond}, the sequence
\eqref{ext3} is a complex
of $A$-bimodules.
\qed
\end{lem}

We call
the complex formed by the first
three, resp. four, terms in \eqref{ext3}
 the {\em  cotangent complex} (to be denoted
${\mathbb L}\Om^1_R(F,\al)$),
resp.,   {\em  extended
cotangent complex}, of $\A(F,\al)$.

\begin{examp}[Quiver case]
Let $Q$ be a quiver with vertex set $I$ and put
$R=\k I$. 

Let $E$
be a $\k$-vector space with basis formed by the
edges of $Q$, and let $E^*=\Hom_\C(E,\C)$ be the dual
vector space, with the dual basis $\{x^*,\, x\in Q\}$. The spaces $E$
and $E^*$ have obvious
$R$-bimodule structures. Moreover, there is a canonical algebra
isomorphism $\C Q=T_R E.$

Now, fix a potential  $\Phi\in (\C Q)\br$, let $\al=d\Phi,$
and $A=\A(\C Q, \Phi)$. For any pair
of edges $x,y\in Q$, we have an element
$\frac{\pa^2\Phi}{\pa x\pa y}\in A\o A$.
We define an $A$-bimodule map
\beq{imap3}
A\otimes_R E^*\otimes_R  A\to
A\otimes_R E\otimes_R  A,\quad
a\o x^* \o b\mto
\sum\limits_{y\in Q}\, a
\left({\Large\frac{\pa^2\Phi}{\pa x\pa
y}}\right)'\o y\o
\left({\Large\frac{\pa^2\Phi}{\pa x\pa y}}\right)''\!b.
\eeq
which is a quiver version of the map \eqref{imap}.
We also introduce 
an $A$-bimodule map
$$
j:\ (A\o A)^R\to A\otimes_R E^*\otimes_R  A,\quad
a\o b\mto\sum_{x\in Q} (bx\o x^*\o a - b\o x\o x^*a).
$$
The dual map  is given by the formula
$$j^\vee: A\otimes_R E\otimes_R  A\too
A\o_R A,\quad
a\o x\o b\mto ax\o b - a\o xb.
$$

 With this understood, we
have

\begin{prop} In the special case where $F=\C Q$  and $\al=d\Phi$,
the map
$d\ccirc\bi$ in \eqref{map} reduces to  \eqref{imap3}; furthermore, the complex in
\eqref{ext3} takes the form
$$
0\to(A\o A)^R\stackrel{j}\too
\xymatrix{
  A\otimes_R E^*\otimes_R  A\ar[r]^<>(0.5){\eqref{imap3}}&
A \otimes_R  E\otimes_R  A}
\stackrel{j^\vee}\too 
A\otimes_R A\to A\to 0.
$$
\end{prop}
\end{examp}

\subsection{The DG algebra $\ffd(F,\al)$}\label{dg}
Fix a smooth algebra $F$ and a cyclic 1-form $\al\in(\Om^1_RF)\br$ such
that \eqref{cond} hold. Motivated by  isomorphisms
\eqref{RCT}, \eqref{natiso}, we introduce
a  graded algebra
\beq{FDA}
\FD(F,\al):=(\Th^\hdot_RF)*_RR[t],
\eeq
a free product over $R$ of the graded
 algebra  of noncommutative
polyvector fields and a polynomial algebra $R[t],$
where $\deg t=2$. 
We define further a super-derivation
$\d: \FD_\idot(F,\al)\to\FD_{\idot-1}(F,\al)$ by the following
assignment,
cf. \eqref{d_def},
\beq{ddi}\d:\
t\mto\De,\quad f\mto0\en \forall f\in F=\FD_0(F,\al),\quad
\th\mto\bi_\th\al \en \forall \th\in\dder F=\FD_1(F,\al).
\eeq

From the second equation in \eqref{cond} we deduce
$\d^2(t)=\d(\De)=\bi_\De\al=0$. It follows that $\d^2=0$, hence,
we have made $\FD_\idot(F,\al)$ a DG algebra.

There is an alternative `Hamiltonian' interpretation of the 
above defined differential
$\d$ in the spirit of \S\ref{cy_rep}.
To explain this,  recall first that
the graded algebra  $\Th^\hdot_R F$
comes equipped with a canonical
noncommutative {\em odd} symplectic 
cyclic 2-form $\om\in \Om^2_R((\Th^\hdot_RF))\br,$ cf. \cite{CBEG}, \cite{VdB4},
\cite{GS2},
which is analogous to the canonical (odd) symplectic structure
on the odd  cotangent bundle of a manifold, cf. \S\ref{info}. 

We consider the following maps
\beq{bij}
(\Om^1_RF)\br\into(\Om^1_R(\Th^\hdot_RF))\br
\stackrel{\om}\iso \Der_R(\Th^\hdot_RF,\,\Th^\hdot_RF),\qquad
\beta\mto \bth_\beta.
\eeq
Here, the first map
is induced
by the natural algebra imbedding $F\iso \Th^0_RF\into \Th^\hdot_RF$,
and the second map is an inverse of the bijection $\xi\mto i_\xi\om$
given by contraction with the symplectic 2-form $\om\in
(\Om^2_R(\Th^\hdot_RF))\br$.

We observe that the grading on $\Th^\hdot_RF$
makes $(\Om^2_R(\Th^\hdot_RF))\br$ a graded vector space,
and we have $\deg\om=1$. Therefore, contraction with $\om$,
hence, the composite in \eqref{bij}, is
a parity reversing map. 
Since any  $\al\in  (\Om^1_RF)\br$  is an {\em even} 1-form of degree zero,
the image of $\al$  via  the composite in \eqref{bij}
is an {\em odd} Hamiltonian super-derivation
 $\bth_{\al}\in \Der(\Th^\hdot_RF)$, of  degree $-1$.

The following result provides a
Hamiltonian point of view on
the DG algebra $\big(\FD_\idot(F,\al),\d\big)$.

\begin{prop}\label{hamprop} \vi The  triple $(\Th^\hdot_RF, \om, \bth_{\al})$
is a  symplectic data 
in the sense that conditions \eqref{D}$\mathsf{(i)-(iii)}$ hold,
except possibly the condition that $\bth_{\al}(\Th^\hdot_RF)\cap R=0$.

\vii The DG  algebra $\big(\FD_\idot(F,\al),\d\big)$ is
isomorphic to $\FD_\idot(\Th^\hdot_RF, \om,
\bth_{\al})$, the  DG  algebra 
associated to the symplectic data in $\mathsf{(i)}$ by the
`universal construction' of \S\ref{potential_sec}.
\end{prop}

\begin{rem} \vi Given  a data $(D,\om,\xi)$, the  condition
that  $\xi(D)\cap R=0$, that appears in Proposition \ref{hamprop}(i),
 is an $R$-analogue of the
requirement from \eqref{D}$\mathsf{(iii)}$
saying that $\xi(D)\cap\C=0$.  We will see in the course of the proof of the
proposition
that  the condition that $\xi(D)\cap R=0$ implies
the equation $\bi_\De\al=0$, cf. \eqref{cond}, but is not equivalent to
it, in general. However, as we have seen earlier
from the explicit formula for the differential $\d$,
the latter condition
is still
sufficient to insure that $\d^2=0.$

\vii
The special features of the  symplectic data  $(\Th^\hdot_RF, \om, \bth_{\al})$, as compared
to more general data $(D,\om,\xi)$  considered in   \S\ref{potential_sec}, are:
\vskip 1pt
\npb{The DG algebra  $D=\Th^\hdot_RF$ comes from  noncommutative  polyvector fields;}

\npb{The symplectic derivation $\xi=\bth_{\al}$ is Hamiltonian, and the
corresponding potential belongs to the subspace
$F\br=(\Th^0_RF)\br\sset (\Th^\hdot_RF)\br.$\hfill$\lozenge$}
\end{rem}

\begin{rem} One can state and prove a generalization of Proposition
\ref{cotangent_thm} saying that there is an isomorphism 
between the complex $0\to\Om^1(\FD|\A)\to\A\o\A\to0$
and the extended cotangent complex ~\eqref{ext3}.
\hfill$\lozenge$
\end{rem}
\begin{proof}[Sketch of proof of Proposition \ref{hamprop}] The proof of part (i)
is straightforward and is left for the reader.

To prove (ii), for any 1-form $\al\in (\Om^1_RF)\br$,  on the algebra
$\Th^\hdot_RF$,
we introduce
a super-derivation $\arr\al$
defined by
the following assignment in degrees 0 and 1, respectively:
\beq{cco}
\arr\al:\ f\mto 0\en\forall f\in F=\Th^0_RF,
\quad \th\mto\bi_\th\al\en\forall \th\in\dder_RF=\Th^1_RF.
\eeq
It is easy to verify that this assignment can be uniquely
extended to a super-derivation 
$\arr\al: \Th^\hdot_RF\to\Th^{\hdot-1}_RF$.

On the other hand, one has the Hamiltonian derivation $\bth_\al$,
see \eqref{bij}. We claim that
\beq{arr}
\bth_\al=\arr\al,\qquad\forall
\al\in (\Om^1_RF)\br.
\eeq

The proof of this equality copies the computation from
the proof of \cite{CBEG}, Lemma 5.4.2. The double derivation denoted
by $\arr\al$ in {\em loc. cit.} is, in effect, a double analogue
of the derivation  $\arr\al: \Th^\hdot_RF\to\Th^{\hdot-1}_RF,$
given by \eqref{cco}. Thus, the computation in our
present setting is slightly simpler than the one in \cite{CBEG}, Lemma
5.4.2; we leave it to the reader.

Recall next that, in the algebra $\Th^\hdot_RF$,
one has $\omh=\De$, by \cite{CBEG}, Theorem 5.1.1.
The statement of the proposition now follows immediately from
\eqref{arr}, by comparing formula \eqref{ddi} with \eqref{ddata}
and \eqref{cco}.
\end{proof}

Next, we consider the two-sided
ideal $I:=\imm\sset F$, and  let $\wh F$ be the  corresponding
$I$-adic completion of the algebra $F$. 
For each $r=0,1,2,\ldots,$
let $\wh\FD_r(F,\al)$ be the $I$-adic completion of 
the $F$-bimodule $\FD_r(F,\al)$, the $r$-th homogeneous
component of the graded algebra $\FD(F,\al)$.
The space $\wh\FD_\idot(F,\al):=\bigoplus_{r\geq 0}\wh\FD_r(F,\al)$
has an obvious graded algebra structure, and we have
$\wh\FD_0(F,\al)=\wh F.$
Further, the derivation $\bth_{\al}$ annihilates the ideal $I$,
hence, extends to a continuous derivation
$\wh\FD_\idot(F,\al)\to\wh\FD_{\idot-1}(F,\al)$.
 This way,
we obtain a completed symplectic
 data $\big(\wh\FD(F,\al), \om, \bth_{\al}\big)$,
and the corresponding completed  DG algebra
$\ffd(F,\al):=\wh\FD(F,\al)*_RR[t]$.

To summarize, we  have a diagram 
which is  a generalization of diagram \eqref{HD},
\beq{Ares}
\xymatrix{
\FD(F,\al)\ar@{^{(}->}[r]&\ffd(F,\al)
\ar@{->>}[r]&H_0\big(\ffd(F,\al),\,\d\big)=H_0\big(\FD(F,\al),\,\d\big)=\A(F,\al).
}
\eeq

\subsection{Main result}\label{cot} Here is
 our  main result on Calabi-Yau
 algebras of dimension 3.

\begin{thm}\label{CYconditions} Let $F$ be a smooth $R$-algebra. Fix
$\al\in (\Om^1_RF)\br$ such that conditions  \eqref{cond} hold.
Then, the following properties of the algebra
$\A:=\A(F,\al)$ are  equivalent:
\vskip 2pt
\vi $\,\A$  is a Calabi-Yau algebra
of dimension $3$;

\vii The map $\dis(\A\o \A)^R \to \dder_R\A, \, a'\o a''\mto \ad(a'\o a''),$
is a bijection.

\viii  The extended cotangent complex \eqref{ext3}  is a resolution
of $\A$.

\iv The DG algebra $\ffd(F,\al)$ is acyclic in positive degrees.
\end{thm}

If the equivalent conditions above hold,
then the class $\pi\in\Ext^3_{\bimod{\A}}(\A,\A\o \A)$
of the extension  \eqref{ext3} provided by Theorem \ref{CYconditions}(iii)
 is a {\em volume} for
the CY algebra $\A$, defined below  formula \eqref{Iso}.

\begin{rem}\label{rem3} In  Theorem
\ref{CYconditions}(iv), the completed  DG algebra $\ffd(F,\al)$
can not be replaced, in general, by $\FD(F,\al)$, its non-complete analogue.
However, the proof of the theorem shows that the following
implication is true:
{\em  If the DG algebra $\FD(F,\al)$ 
 is acyclic in positive degrees then $\A(F,\al)$ is a CY  algebra
of dimension $3$.}\qed
\end{rem}

\begin{cor}\label{corthm} Let $F$ be a smooth $R$-algebra  and  $\Phi\in F\br,$ a potential.
Then, the  algebra $\A=\A(F,\Phi)$   is a CY algebra
of dimension $3$ if and only if either of the following
 equivalent conditions holds

\begin{enumerate}
\item  $H^0(\A,\A\o \A)=0$ and $H^1(\A,\A\o \A)=0$;

\item The extended cotangent complex \eqref{ext3}  is a resolution
of $\A$. \qed
\end{enumerate}
\end{cor}

The conjecture below is a more precise form of the
claim made in \S\ref{pot} saying that {\em any} CY algebra 
 of dimension $3$ has the form $\A(F,\al)$.

\begin{conj}\label{proj}
Let  $A$ be a friendly  Calabi-Yau $R$-algebra  of dimension $3$.

Then,  there exists
a  smooth  $R$-algebra $F$ and  a cyclic 
1-form $\al\in (\Om^1_RF)\br$ such that
 we have
$\dis A\cong \A(F,\al)$ and, moreover, the
1-form $\al$ satisfies  conditions \eqref{cond}.
\end{conj}

\begin{proof}[Proof of Theorem \ref{CYconditions}]   First of all, we observe that
condition (ii) of  Theorem \ref{CYconditions} is equivalent
to condition (1) of Corollary \ref{corthm}. This shows that
the theorem implies 
the corollary; it also proves the implication
(i)$\;\Rightarrow\;$(ii) in the theorem.

 To prove other  implications we assume for simplicity that $R=\C$.
We will freely use
the notation and results of \S\ref{cy3_pf}.  

Fix a  smooth algebra
$F$ and let $A=\A(F,\al)$.
By \eqref{DD1}-\eqref{DD2},
we know that 
$\DPA$ and $\Om^1(F|A)$ are finitely generated
 projective $A$-bimodules dual to each other.
Thus,  the first $4$ terms in \eqref{ext3} form
a selfdual sequence, to be denoted ${\mathsf P},$ of projective $A$-bimodules.
We conclude that if the sequence of the Theorem is exact,
then ${\mathsf P}$ is a  selfdual projective resolution of
$A$, the diagonal $A$-bimodule.
It follows that, in this case, for any $A$-bimodule $M$, we have
$$\Ext_{\bimod{A}}^\hdot(A,M)=H^\hdot\big(\Hom_{\bimod{A}}({\mathsf P},M)\big),
\quad
\Tor^{\bimod{A}}_\idot(A,M)=H^\hdot({\mathsf P}\o_{A^e}M).
$$

Therefore, from the selfduality of ${\mathsf P}$ one deduces  canonical
isomorphisms \eqref{vdb}.
This implies,  as explained in Remark \ref{vdb_crit},
 that $A$ is a CY algebra of dimension 3, hence completes
the proof of the
implication (iii)$\;\Rightarrow\;$(i) of the theorem.

To prove that  (ii)$\;\Rightarrow\;$(iii),
assume that
 the map $\ad: A\o A\to \Der(A,A\o A)$ is an isomorphism.
Thus, the exact sequence in \eqref{chain2}  reads
$$
0\to A\o A\stackrel{j}\too\Der(F,A\o A)
\stackrel{d\bi^\al}\too
\Om^1(F|A)\stackrel{p}\too\Om^1A\to 0,
$$

We splice this exact sequence
with the tautological short exact sequence
$\Om^1A\into A\o A$ $\onto A$. This
yields the resolution in part (iii) of the Theorem. 

Finally, the equivalence (iii)$\;\Leftrightarrow\;$(iv)
follows from Proposition  \ref{kunneth} (to be proved
later in \S\ref{prop_pf}) applied to the DG algebra $\FD=\FD(F,\al)$.
\end{proof}

\begin{rem} Applying Proposition \ref{kunneth} in a similar situation
provides 
a short alternative proof of a result due to Anick \cite{An}.

Specifically, let  $A=\bigoplus_{r\geq 0}A_r$ be a graded algebra.
Anick has
associated with $A$  two objects: a certain complex
$C_\idot A$, nowadays called {\em Anick's complex},
and also a DG algebra $K_\idot A$,  referred to in \cite{EG} as
{\em noncommutative Koszul complex}.
The result proved by Anick says that the  DG algebra $K_\idot A$
is acyclic in positive degrees iff Anick's complex is a resolution
of $A$, the latter being viewed as an $A$-bimodule. 
 Anick's result clearly follows from Proposition  \ref{kunneth}
 applied to the DG algebra $\FD:=K_\idot A$ using the observation
that $C_\idot A=\Om^1(K_\idot A),$ cf. \cite{EG}.
\erem

\subsection{Graded CY algebras of dimension 3}\label{gr_sec}
The special case of {\em graded} CY algebras 
of dimension 3 has been studied in detail 
in \cite{BT}, \cite{Bo}. We only briefly recall the main points.

Throughout this section \S\ref{gr_sec}, we fix a finite set $I$ and put $R=\C I$.
We also fix $V=\bigoplus_{i>0} V_i$, a finite dimensional
positively graded $R$-bimodule. 
Thus, $T_RV$, the tensor algebra of the $R$-bimodule $V$,   is
 a smooth nonnegatively
graded $R$-algebra such that $T^0_RV=R$.

R. Bockland proved that Conjecture \ref{proj} holds in the graded case
(even without the assumption that the algebra be friendly). Specifically,
 we have
\begin{thm}[Bockland, \cite{Bo}] Let $I\sset \bigoplus_{i\geq 2} T^i_RV$ be a graded two-sided
ideal such that the algebra $A=T_RV/I$ is a CY algebra of dimension 3.
Then, there exists a homogeneous potential
 $\Phi\in (T_RV)\br$, of degree $\deg\Phi>2$ such that
 $I=\llb\frac{\pa\Phi}{\pa v}\rrb_{\{v\in V\}}$.\qed
\end{thm}
\begin{rem} Bockland only considered
a special case where $V=V_1$ is a homogeneous graded
bimodule placed in degree 1, but the proof in \cite{Bo} works
in the general case as well.
\erem
 
In the graded case, our proof of Theorem \ref{CYconditions} is based on the
graded case of Proposition \ref{kunneth}(ii). This way, from Theorem
\ref{CYconditions},
cf. also Corollary \ref{corthm},
we obtain
\begin{cor}\label{CYgraded} Let  $(R,V)$ be as above,
and let $\Phi\in (T_RV)\br$,
be a homogeneous potential of degree $\deg\Phi>2$.

Then,  the following properties of the algebra
$\A:=\A(T_RV,\Phi)$ are  equivalent:
\vskip 1pt
\begin{enumerate}
\item   $\,\A$  is a Calabi-Yau algebra
of dimension $3$;

\item  $H^0(\A,\A\o \A)=0$ and $H^1(\A,\A\o \A)=0$;

\item    The extended cotangent complex \eqref{ext3}  is a resolution
of $\A$;

\item   The (noncompleted) DG algebra $\FD_\idot(T_RV,\Phi)$ is acyclic
in positive degrees.
\qed
\end{enumerate}
\end{cor}

The following immediate corollary of Theorem \ref{CYconditions}
and Corollary \ref{CYgraded}
has been first established by Berger and Taillefer in the recent paper \cite{BT}. 

\begin{cor}\label{corBT} Let   $(R,V)$ be as above.
Let $\Phi\in (T_RV)\br$,
be a homogeneous potential of degree $\deg\Phi=m>2$,
such that $\A(T_RV,\Phi)$ is a graded CY algebra of dimension ~3.

Then, for any (not necessarily homogeneous) potential
 $\Psi=\Psi_{m-1}+\Psi_{m-2}+\ldots+\Psi_3\in(T_RV)\br,$
of degree $<m$, the algebra
$\A(T_RV,\Phi+\Psi)$ is a  CY algebra of dimension 3.
\end{cor}

\begin{proof} Let $F_k:=\bigoplus_{r\leq k} (T_RV)_r,
k=0,1,\ldots,$ be the
increasing filtration on the tensor algebra $F=T_RV$
induced by
the grading $T_RV=\bigoplus_{r\geq 0} (T_RV)_r$. The filtration
on $F$ induces an increasing filtration on the free algebra
$\FD(F,\Phi+\Psi)$. The latter filtration is compatible with the
differential $\d$ on $\FD(F,\Phi+\Psi)$. Furthermore, it is clear
that, for the correspoding associated graded algebra, one has
a DG algebra isomorphism
$$\gr\FD(F,\Phi+\Psi)\cong\FD(F,\Phi),$$
where the DG algebra on the right is associated with the
{\em homogeneous} potential $\Phi$.

The increasing filtration on the  DG algebra
$\FD(F,\Phi+\Psi)$ gives rise to  a standard converging spectral sequence 
$$
E_2=H_\idot\big(\gr\FD(F,\Phi+\Psi),\,\d\big)=
H_\idot\big(\gr\FD(F,\Phi),\,\d\big)\en\Rightarrow\en
E_\infty=\gr
H_\idot\big(\FD(F,\Phi+\Psi),\,\d\big).
$$

Now, Corollary \ref{CYgraded} implies that
the DG algebra $H_\idot\big(\FD(F,\Phi),\,\d\big)$ is acyclic
in positive degrees, since $\A(F,\Phi)$ is a CY algebra by 
our assumptions. It follows, that the $E_2$-term,
hence also the $E_\infty$-term, in the
above spectral sequence vanishes in positive degrees.
We deduce that the DG algebra $H_\idot\big(\FD(F,\Phi+\Psi),\,\d\big)$
is acyclic
in positive degrees as well. Thus,  Theorem \ref{CYconditions}
implies that $\A(F,\Phi+\Psi)$ is a CY algebra of dimension 3,
cf. Remark ~\ref{rem3}.
\end{proof}

\begin{notation} Given
a graded $R$-bimodule $E=\bigoplus_{r\geq 0}E_r,$ where  $E_r$
is a finite dimensional $R$-bimodule for any $r=0,1,\ldots,$
define the matrix Hilbert series $h(E;t)$, an $I\times I$-matrix
with entries
$$
h_{ij}(E;t):=\sum_{r\geq 0} t^r\cd\dim(1_i\cd E_r\cd 1_j)\in \Z_{\geq
0}[[t]],\qquad i,j\in I,
$$
cf. \cite{EG}, \S2.2 for more details.
Let $\bone$ denote the identity $I\times I$-matrix. 
\end{notation}

\begin{defn}\label{cartan} Given a  pair $(V,\Phi)$, where
$V$ is a finite dimensional
positively graded $R$-bimodule and $\Phi\in (T_RV)\br$ is
 a homogeneous potential, we define

$$\pp(V,\Phi;t):=
\bone-h(V;t)+t^m\cd h(V;t\inv)-t^m\cd \bone,\quad m:=\deg\Phi>2.
$$
This is a polynomial with matrix coefficients, called 
{\em Cartan polynomial} for the pair ~$(V,\Phi)$.
\end{defn}

From the exactness of the extended cotangent complex,
using the Euler-Poincar\'e principle,
one derives the following result that has been also
proved by Bockland, \cite{Bo}.

\begin{prop} Let  $(V,\Phi)$ be a pair as above such that
the  algebra  $A=\A(F,\Phi)$  is a Calabi-Yau algebra
of dimension $3$. 

Then, the  (matrix) Hilbert  series of $A$ equals 
$h(A;t)=\pp(V,\Phi;t)\inv$.\qed
\end{prop}

Observe that the grading on $V$ induces an additional $\Z_+$-grading on
each of the following algebras
$\FD(F,\Phi), \, \Sym\FD(F,\Phi), $ and $\oo_\dd\big(\FD(F,\Phi)\big).$
We write $\chi\Big(\oo_\dd\big(\FD(F,\Phi)\big);t\Big)\in
\Z[[t]]$
for the graded Euler characteristic of the DG algebra
$\oo_\dd\big(\FD(F,\Phi)\big)$. We have the following $\zeta$-function type formula
\beq{zet}
\chi\Big(\oo_\dd\big(\FD(F,\Phi)\big);t\Big)=
\prod_{s=1}^\infty \frac{1}{\det\big(\pp(V,\Phi;t^s)\big)}.
\eeq
This formula follows from \cite{EG}, formula (2.2.9)(ii), applied to the free graded
algebra $\FD(F,\Phi)$.

Formula \eqref{zet}
is useful for finding  {\em graded} Euler
 characteristics of cyclic homology.
 Specifically, a grading on an algebra $A$ gives,
for each $p\geq 0$,
 an
extra grading on 
 $\overline{HC}_p(A)$, the $p$-th {\em reduced}  cyclic homology group
of $A$. The graded Euler
 characteristic is then defined as 
$\sum_{p\geq 0} (-1)^p\cdot h\big(\overline{HC}_p(A);t\big).$
Using the same notation as in \cite{EG}, \S3.7, we have

\begin{prop}\label{eg}
Let  $\A:=\A(T_RV,\Phi)$
be a graded CY algebra of dimension 3. 
Then, 

\noindent
\vi The graded Euler characteristic of the  algebra
$\oo_\dd\big(\overline{HC}_\idot(\A)\big)$ is given by formula
~\eqref{zet}.

\noindent
\vii We have  $\dis H_p\big(\FD(T_RV,\Phi)\br,\,\d\big)=0$ for all $p>2$; furthermore, one has
$$H_3(\A,\A)\cong \overline{HC}_2(\A)\cong
H_2\big(\FD_\idot(T_RV,\Phi)\br,\,\d\big).
$$
\end{prop}
\begin{proof} The projection
$\FD_\idot=\FD(T_RV,\Phi)\onto\A(T_RV,\Phi)$
is a quasi-isomorphism for any CY algebra
 $\A:=\A(T_RV,\Phi)$  of dimension
3, by   Corollary \ref{CYgraded}(4). 
Therefore, according to \cite{Lo}, the reduced  cyclic homology of such an algebra
may be computed using the complex
$\big(\FD_\idot/(R+[\FD_\idot,\FD_\idot]_\text{super}),\,\d\big)$.
Part (i) now follows from
formula \eqref{zet}. 

Part (ii) is immediate from \cite{EG}, Lemma 3.6.1.
\end{proof}

It is known that
any graded 
 Calabi-Yau algebra of dimension 3 is Koszul and Gorenstein of dimension $3$,
but the converse is not necessarily true, in general,
see \cite{BT}, \cite{DV}, \cite{BM}.

\begin{rem}\label{Raf} For any Koszul  algebra $A=\bigoplus_{r\geq 0}
A_r$,
let $A^!=\Ext^\hdot_{\Lmod{A}}(\C,\C)$ be the
Koszul dual quadratic algebra.
The following result, pointed out to me by  R. Rouquier,
will be proved after the remark.
\begin{lem}\label{kos_lem} For any Koszul algebra $A$ such that
$\dim A^!<\infty$ one has a natural graded algebra isomorphism
$$
H^\hdot(A^!,A^!)\cong H^\hdot(A,A).
$$
\end{lem}

Assume now that $R=\C$ and let $A$ be a graded 
 Calabi-Yau algebra of dimension 3. Since $A$ is Koszul and  Gorenstein of dimension $3$,
we have $A^!_3=\C$, and $A^!_j=0$ for all $j>3$; in particular,
$A^!$ is a finite dimensional algebra.  Moreover,
the multiplication map  gives a nondegenerate
inner product  $A^!_j\times A^!_{3-j}\to A^!_3\iso\C$, on $A^!$.

R. Rouquier conjectures that the  BV algebra structure
on $H^\hdot(A^!, A^!)$ constructed by Tradler using the inner product
on $A^!$, cf. 
Remark \ref{tradler}, goes under the isomorphism
of Lemma \ref{kos_lem}, to the BV algebra structure
on $H^\hdot(A, A)$ given by Theorem \ref{BV}(ii).
\erem

\begin{proof}[Proof of Lemma \ref{kos_lem}]
The isomorphism of the lemma can be deduced  from the description of
Hochschild cohomology of a Koszul algebra in terms of Koszul complexes,
cf. \cite{VdB2}. 

Alternatively, to prove the lemma,
one may observe that Hochschild cohomology  of an algebra $A$
has a definition in terms of $D=D(\bimod{A})$, the derived category of
$A$-bimodules. Specifically,
write $\Id_D$ for the identity
functor on $D$ and, for each $r\in\Z$, let ${[r]}_D$ denote the
corresponing shift functor on  $D$.
Then, one has 
$$H^\hdot(A,A)=\Ext^\hdot_{\bimod{A}}(A,A)\cong\oplus_{r\geq 0}
\Hom_D(\Id_D, {[r]}_D).
$$

Now, according to Koszul duality \cite{BGS},
given a Koszul algebra $A$ such that $A^!$ is a finite dimensional
algebra,
 one has a triangulated equivalence $D^b(\bimod{A})\cong D^b(\bimod{A^!})$.
One can use this equivalence and the above definition
of Hochschild cohomology  to deduce the isomorphism of Lemma \eqref{kos_lem}.
\end{proof}

The reader is referred to \cite{Bo}, \cite{BT}, and \cite{DV}
for various examples of  graded 
 Calabi-Yau algebra of dimension 3.  These include, in particular,
 quadratic and cubic Sklyanin
algebras considered in Example  \ref{sklya},
as well as 
Yang-Mills algebras  considered in 
\cite{CDV2}.

{\large\section{{\textbf{Fundamental groups,  Chern-Simons, and a noncommutative \newline
analogue of the Weil
 representation}}}\label{3d}}
\subsection{Group algebras as CY algebras (after M. Kontsevich)}
 Recall that a compact connected $C^\infty$-manifold $X$ is
called {\em aspherical} if
its universal cover is contractible,
equivalently, if the   homotopy groups 
$\pi_r(X)$, vanish  for all $r\geq 2.$
Examples of 3-dimensional aspherical manifolds include
all hyperpolic spaces of the
form $X={\mathbb H}/\G$, where ${\mathbb H}$ is
the 3-dimensional Lobachevsky space of constant
negative curvature, and $\G$ is a discrete cocompact
group of isometries of  ${\mathbb H}$.

The following result along with its proof was kindly
communicated to me by M. Kontsevich.

\begin{prop} Let $X$ be an oriented 3-dimensional aspherical manifold.
Then, $\C[\pi_1(X)]$, the group algebra of the
 fundamental group of $X$, is a CY algebra of dimension ~3.
\end{prop}

In general, let $X$ be a space homotopy equivalent to a finite CW complex.
We consider the DG algebra $A(X)$ of singular chains in the based loop space
$\Omega(X,x),$ where
$x$ is a base point, with the multiplication coming from the composition of
loops.

\begin{thm}[M. Kontsevich]\label{max_pf} The  $\oper{DG}$ algebra $A(X)$
is 
homologically smooth.
Moreover, it is finitely presented in the sense of \cite{TV}.
\end{thm}
\proof[Sketch of Proof]
Following Kontsevich, we will construct an
 explicit DG algebra resolution of $A(X)$ as follows.

Consider a simplicial set $S$ with finitely many non-degenerate simlpices.
Associate with it a DG category ${\scr C}_\Z$ (over $\Z$, the ring of
integers),
whose objects are
$0$-simplices. The corresponding $\Hom$'s are defined
as follows.  Any non-degenerate simplex $s$ with vertices
$(s_0,...,s_n)$,
so
$\dim s=n>0,$
gives a morphism from $s_0\to s_n$ of degree $(1-n).$
The algebra of $\Hom$-spaces is, by definition,
 freely generated by the  morphisms above.
Thus,  our category is the path category associated to a
quiver.
The differential is defined by
$$d(s_0,...,s_n)= \sum_{0<i<n} (-1)^i\big((s_0,..,s_{i-1},s_{i+1},...,s_n)+
(s_0,...,s_i) \ccirc (s_i,...,s_n) \big).$$
 Let ${\scr C}'_{\Z}$ be obtained from ${\scr C}_{\Z}$ by inverting the
set of generating morphisms of degree 0
(corresponding to 1-simplices).

Theorem \ref{max_pf} is now a consequence of the following result.
\begin{prop} The complex  $\Hom_{{\scr C}'_{\Z}}(s,s')$  calculates the homology
of the
path space $\oper{Path}(s,s').$\qed
\end{prop}


\begin{cor}[M. Kontsevich]\label{p1}
Let $X$ be a  connected finite CW complex which is, moreover,
a $K(\pi,1)$ space. Then,

\vi $\C[\pi_1(X)],$ is a homologically smooth  algebra.

\vii  Assume, in addition, that $X$ is a Poincar\`e duality space, e.g., a
smooth compact oriented manifold, of dimension $n$.
Then,  $\C[\pi_1(X)],$ is a CY algebra of dimension $n$.\qed
\end{cor}

 According to Corollary \ref{p1}(i), examples of groups $\pi$ with
 homologically smooth and
finitely presented group algebra $\C[\pi]$ include
$S$-arithmetic groups, braid  groups, Teichmuller groups, and
 fundamental
groups of
knots.

\subsection{Potentials for the group algebra $\C[\pi_1(X)]$}
Fix an integer $\dd\geq 1$. We have
$$\Rep_\dd\C[\pi_1(X)]/\GL_\dd=\Rep_\dd\pi_1(X)/\GL_\dd\cong
\loc(X),$$
is the moduli stack of rank $\dd$ local systems, equivalently,
vector bundles 
on $X$ (viewed as a space with the ordinary  Hausdorff topology) with flat
connection. 

Thus, given an aspherical 3-manifold $X$,
any presentation of the form $\C[\pi_1(X)]=\A(F,\Phi)$
would identify the stack $\loc(X)$ with the stack of critical points
of the corresponding function $\Tr\wh\Phi.$
There are two possible conjectural approaches for finding 
such a presentation
$\C[\pi_1(X)=\A(F,\Phi)$.

One approach is to use the space $\conn$,
of
all $C^\infty$-connections $\alpha$ on a rank $\dd$  vector
bundle on $X$. We may view $\loc(X)$  as a substack in $\conn/\text{\em Gauge
group}$.
Specifically, it is well known that, in the
case of a trivial bundle, one has
$\loc(X)=\crit(CS)$, the critical set of the Chern-Simons functional
$$CS:\
\Om^1_{C^\infty}(X)\o\gl_\dd\to \C,\quad\al\mto
 CS(\alpha)=\int_X \Tr\left(\frac{1}{2}\al\wedge d\al +\frac{1}{3}\al^3\right),$$
 on the infinite dimensional vector space    of all
$\gl_\dd$-valued smooth 1-forms on $X$.
Thus, it might be plausible to look for a presentation
of the form $\C[\pi_1(X)]=\A(F, CS_\text{alg})$, where
$F=\C^\infty(X,\, \C\langle x_1,\ldots,x_\dd\rangle)$
is the algebra of smooth maps $X\to \C\langle x_1,\ldots,x_\dd\rangle,$
and where $CS_\text{alg}$ is an appropriate version
of the  Chern-Simons functional defined using Chen's theory
of iterated integrals.

A second approach is based on Heegaard splittings and a  noncommutative 
analogue of the Weil representation. In more detail,
let $X=X_+\cup X_-$ be a  Heegaard splitting of our
3-manifold as a union of two discs, $X_+$ and $X_-,$ such
that $X_+\cap X_-=C$ is a genus $g$ Riemann surface imbedded into $X$.
Let $\loc(C)=\Rep_\dd\pi_1(C)/\GL_\dd$ be the  moduli space of 
rank $\dd$ local systems on $C$.

Further, let $L_\pm:=\ker[\pi_1(C)\to\pi_1(X)]$.
The groups $L_\pm$ are commutative, they are rank $g$ lattices
in the group $\pi_1(C)$. The algebra map
$\C[L_\pm]\to\C[\pi_1(C)]$ goes, via the representation
functor, to the natural imbedding
$$\loc(X_\pm)=\Rep_\dd(\C[L_\pm])/\gld\into \loc(C)=\Rep_\dd(\C[\pi_1(C)])/\gld,
$$
where  $\loc(X_\pm)\sset \loc(C)$ stands for
the subset of those local systems
on $C$ which extend to the disc $X_\pm$. 

The algebra $\C[\pi_1(C)]$ has a  noncommutative  Poisson structure essentially
constructed by Goldmann \cite{Go}, cf. also \cite{Et}.
Accordingly, for any $\dd$, the smooth locus of the moduli space
$\loc(C)$ is known to be a symplectic manifold.
It is known
that each of the sets $\loc(X_\pm)$ is a Lagrangian subvariety and,
moreover, one clearly has $\loc(X)=\loc(X_+)\cap \loc(X_-)$.

Now, let $P:=\C\langle x_1,\ldots,x_g,y_1,\ldots, y_g\rangle.$
Travis Schedler has defined in \cite{Sc}, cf. also \cite{GS1}, a
Hopf algebra  $\U_q(P\br),$ a certain quantum deformation of
the necklace Poisson  algebra $\Sym(P\br)$.
In \cite{GS2}, we  extend Schedler's
original  construction and introduce a  noncommutative  version of
the Heisenberg representation of the algebra
$\U_q(P\br)$ in the vector
space $\C\langle x_1,\ldots,x_g\rangle\br$.
Noncommutative quadratic polynomials in the
variables $x_1,\ldots,x_g,y_1,\ldots, y_g$
form a Lie subalgebra in $\U_q(P\br)$ that turns out
to be isomorphic to the symplectic Lie algebra ${\mathfrak{sp}}_{2g},$
as in the classical case.
Restricting the Heisenberg representation to this
Lie  subalgebra yields a  noncommutative  version of the Weil representation,
see \cite{GS2} for more details.

There should be a `trigonometric counterpart' of  the above mentioned
Heisenberg and Weil representations,
where
 the algebra  $\C\langle x_1,\ldots,x_g,y_1,\ldots, y_g\rangle$
is replaced by $\C[\pi_1(C)]$.
Specifically, let $a_1,\ldots,a_g$ be a set of generators
of the abelian group $L_+$. We write
 $\C\langle L_+\rangle=\C\langle a^\pm_1,\ldots,a^\pm_g\rangle $  for the
algebra of `noncommutative Laurent polynomials'.
We expect that that there is
a Heisenberg type representation of the algebra
$\C[\pi_1(C)]$ in the vector
space $\C\langle L_+\rangle\br$.

Recall further that one has the mapping class group
${\mathbb M}_g$ that acts naturally on $\pi_1(C)$. 
We conjecture, that there is a
 vector space $\C\langle\!\langle L_+\rangle\!\rangle\br$,
 an appropriate completion of $\C\langle L_+\rangle\br$,
and a projective representation
$\wp: {\mathbb M}_g\too \PGL(\C\langle\!\langle L_+\rangle\!\rangle\br)$ that plays the role of 
the Weil representation of the Lie algebra ${\mathfrak{sp}}(\C^{2g})$;
observe that one has a natural group homomorphism
${\mathbb M}_g\to Sp(\Z^{2g})$ induced by the ${\mathbb M}_g$-action
on the symplectic lattice $H_1(X,\Z)$, the integral homology of $C$ 
equipped with  the intersection pairing.

Now, recall that the  Heegaard splitting
gives {\em two} lattices $L_+, L_-\sset \pi_1(C)$.
Furthermore, there exists
an element $u\in {\mathbb M}_g$ such that $L_-=u(L_+).$
Therefore, we have a linear map
$\wp(u): \C\langle\!\langle L_+\rangle\!\rangle\br\to
\C\langle\!\langle L_+\rangle\!\rangle\br,$ which is well defined
up to a nonzero constant factor.

\begin{conj} There is an algebra isomorphism
$\C[\pi_1(X)]\cong\A(F,\Phi)$, where
$$F=\C\langle\!\langle L_+\rangle\!\rangle,\quad\oper{and}\quad
\Phi=\wp(u)(1)\in \C\langle\!\langle L_+\rangle\!\rangle\br.
$$
\end{conj}

{\large\section{{\textbf{Calabi-Yau algebras and Calabi-Yau manifolds.}}}\label{cymnd_sec}}
\subsection{}\label{cymnd_subsec} This section is 
mostly devoted to the proof of Proposition \ref{CYCY}.
The argument  is based on Theorem \ref{CM} below,
which will be proved first.

We begin with our definition of a tilting object.
Let $X$ be a smooth quasi-projective variety,
and 
write $D^b(\coh X)$ for the derived category of
coherent sheaves on $X$. 

\begin{defn}\label{tilt} We call  $\CE\in D^b(\coh X)$
a {\em tilting generator}  if the following three conditions hold:
\vskip 2pt

\npb{$\Ext^i(\CE,\CE)=0$ for all $i>0;$}

\npb{$\Ext^\hdot(\CE,\CF)=0$ implies $\CF=0$, for any $\CF\in D^b(\coh X)$;}

\npb{The algebra $A_\CE:=\Hom(\CE,\CE)$ has finite Hochschild
dimension.}
\end{defn}

It is well-known that, given a tilting object
$\CE$, one has a triangulated equivalence
$D^b(\coh X)\iso D^b(\Lmof{A_\CE})$.
Let $\DZ(\coh X)$ be the center
of $D^b(\coh X)$, that is, the endomorphism algebra
of the identity functor $\Id: D^b(\coh  X)\to D^b(\coh  X).$

\subsection{Calabi-Yau condition and Serre functors.}\label{ser_sec}
For any left noetherian algebra $A$ we have 
an abelian category $\Lmof{A}$, the full subcategory
of the category $\Lmod{A}$ of all
left $A$-modules, whose objects are finitely generated   $A$-modules.
There are similar categories $\Rmod{A}$ and $\Rmof{A}$,   of 
right $A$-modules over a right noetherian algebra.

We fix a (left and right) noetherian algebra $A$.
Let $D^{^{_{_\heartsuit}}}(\Lmof{A}),$ where
${^{_{_\heartsuit}}}=\fpr,\fin,b,+,$
$-,$ denote the  derived category
of bounded complexes of  projective, 
resp., bounded complexes of injective,
 bounded complexes, bounded from below complexes, bounded from above
 complexes,
of  finitely generated $A$-modules.

Let  $Z\subset A$ be a {\em central} subalgebra.
Throughout this section, we will assume that 
\begin{align}
&\text{$\mathsf{(i)}\quad Z$ is a  finitely generated Cohen-Macaulay  $\k$-algebra 
without zero divisors;}\nonumber\\
&\text{$\mathsf{(ii)}\quad A$ is  finitely generated as a $Z$-module.}\label{zfin}
\end{align}

Let $\CK_Z\in D^\fin(\Lmod{Z})$ be the dualizing complex for $Z$.
Thus,  $\CK_Z$ has a single
nonvanishing cohomology group and 
there is  a  contravariant duality
\beq{Zdual}
\RHom_Z(-, \CK_Z):\ D^\fpr(\Lmof{Z})\longleftrightarrow
D^\fin(\Lmof{Z}).
\eeq

For any $Z$-module $K$, the space $\Hom_Z(A,K)$ has a natural
$A$-bimodule structure.
This $A$-bimodule is finitely generated, both as
a left $A$-module and as a right $A$-module,
provided $K$ is finitely  generated over $Z$.
Similarly, for any
$M,N\in D^{^{_{_\heartsuit}}}(\Lmof{A})$, the complex
$\RHom_A(M,N)$ may be viewed as a complex
of $Z$-modules.
Specifically,  we have
$$M\in D^\mp(\Lmof{A})\en\&\en N\in D^\pm(\Lmof{A})
\en\Longrightarrow\en
\RHom_A(M,N)\in D^\pm(\Lmof{Z}).
$$

We say that an  object $U\in D^b(\bimod{A})$  has 
{\em finite injective dimension}
provided it is quasi-isomorphic
to a complex $I$, of  $A$-bimodules, such that
$I\in D^\fin(\Lmof{A})$ as a complex of {\em left} $A$-modules,
and it is also  quasi-isomorphic
to a complex $J$, of  $A$-bimodules, such that
$J\in D^\fin(\Rmof{A})$ as a complex of {\em right} $A$-modules.

We put $\KA:=\RHom_Z(A,\CK_Z).$ 
This is a complex of $A$-bimodules which has a finite injective dimension,
by \eqref{Zdual}.
We  see that the  assignments
below give two 
 well-defined triangulated
functors
\begin{align}
M\mto \BD_A(M)&:=\RHom_A(M,\KA),
\quad &D^\pm(\Lmof{A})\stackrel{\BD_A}\longleftrightarrow
 D^\mp(\Rmof{A});\label{BD}\\
M\mto \BS(M)&:=\KA\lo_AM,
\quad &D^\pm(\Lmof{A})\stackrel{\BS}\too D^\pm(\Lmof{A}).\label{BS}
\end{align}

Furthermore, both   functors preserve the bounded derived category.
Also, we have
\beq{prin}
M\in  D^\fpr(\Lmof{A})\quad\Longrightarrow \BD_A(M), \BS(M)\in D^\fin(\Lmof{A}).
\eeq

\begin{defn}\label{serre_funct} A covariant triangulated
functor $S: D^\spad(\Lmof{A})\to D^\spad(\Lmof{A})$
is called  Serre functor on $ D^\spad(\Lmof{A})$ {\em relative to the
subalgebra} $Z$ provided one has functorial isomorphisms
\beq{serre_iso}\RHom_Z(\RHom_A(M,N),\,\CK_Z)\cong\RHom_A(N, S(M)),
\eeq
 for any pair of objects
$M,N\in D^\spad(\Lmof{A}).$
\end{defn}

A standard argument shows that if such a Serre  functor exists
it is unique up to a unique isomorphism, moreover,
this functor must be an equivalence.

\begin{lem}\label{serre_weak} Assume that  conditions
\eqref{zfin}$\mathsf{(i)}$-$\mathsf{(ii)}$  
hold for
a pair $(Z,A)$. Let $S:=\BS$ be the functor defined in  \eqref{BS}.
Then, there are
 functorial isomorphisms \eqref{serre_iso}, for any
$M\in D^-(\Lmof{A})$ and $N\in D^\fpr(\Lmof{A})$.
\end{lem}

\begin{rem} It should be emphasized that the functor $\BS$ does {\em not}
necessarily preserve the category
$D^\fpr(\Lmof{A})$, cf. \eqref{prin}; that is why, it is  {\em not}
a Serre functor on $D^\fpr(\Lmof{A})$, in general.
\erem

\begin{proof}
Let  $A$ and $B$ be two  noetherian algebras and let
 $U$ be a bounded complex of $(B\dash A)$-bimodules
such that its cohomology groups are finitely generated
both as left $B$-modules  and right $A$-modules.
Given $M\in D^-(\Lmof{A})$,
one has standard canonical isomorphisms
\begin{align}
\RHom_A(M, \RHom_B(U,N))&\iso\RHom_B(U\lo_A M, N),
\quad \forall N\in
D^+(\Lmof{B});\label{re1}\\
\RHom_B(N, U)\lo_A M&\iso\RHom_B(N, U\lo_A M),\label{re2}\\
&\qquad\quad \forall 
  N\in D^\fpr(\Lmof{B})\en\text{or}\en M\in D^\fin(\Lmof{B}).\nonumber
\end{align}

The first, resp. second, of the above isomorphisms is clear for  
$M=A$, resp. for $N=B$. Therefore, it holds in the
case where $M$, resp. $N$, is a finitely generated projective module.
This yields the  isomorphism for  $M$, resp. $N$,
being a  bounded from above, resp.   bounded,
complex of  finitely generated projective modules. This yields \eqref{re2}.
The general case of \eqref{re1}  follows since any object
in $D^- (\Lmof{A})$
is quasi-isomorphic to such a  bounded from above complex
of projectives.

Applying \eqref{re1} for $U=A$ and $N=\CK_Z$, for the functor
$\BD_A$ from \eqref{BD}, we deduce
\beq{y}
\RHom_A(M, \RHom_Z(A,\CK_Z))\cong
\RHom_Z(A\lo_A M, \CK_Z)=\RHom_Z(M, \CK_Z).
\eeq

Now,  for any two objects
$M\in D^-(\Lmof{A})$ and $N\in D^\fpr(\Lmof{A}),$ we have that
 $\KA\lo_A M\in D^-(\Lmof{A}),\,$ $\RHom_A(N,\KA)\in D^b(\Lmof{A}),$
 and
 $\BD_A(N)\lo_A M\in D^-(\Lmof{A}).$
Thus, we
compute
\begin{align*}
\RHom_Z(\RHom_A(N,\KA\lo_A M),\,\CK_Z)&\stackrel{\eqref{re2}}{\en=\en}
\RHom_Z(\RHom_A(N,\KA)\lo_A M,\,\CK_Z)\\
&{\en\stackrel{\text{def}}=\en}\RHom_Z(\BD_A(N)\lo_A M,\,\CK_Z)\\
&\stackrel{\eqref{re1}}{\en=\en}\RHom_A(M,\,\RHom_Z(\BD_A(N),\CK_Z))\\
&\stackrel{\eqref{y}}{\en=\en}\RHom_A(M,\,\RHom_A(\BD_A(N),\KA))\\
&{\en\stackrel{\text{def}}=\en}\RHom_A(M,\,\BD_A\BD_A(N))=\RHom_A(M,N).
\end{align*}

Applying the duality $\RHom_Z(-,\CK_Z)$ to each side completes the
proof of the  lemma.
\end{proof}

In the theorem below, we will also use the  functor
\beq{T}
 M\mto{\mathbb T}(M)=\RHom_{\bimod{A}}(A, A \otimes A)\lo_AM,
\quad D^\fpr(\Lmof{A})\to D^+(\Lmof{A}).
\eeq%

We say that an algebra $A$ has finite injective dimension
if it  has finite injective dimension
as an $A$-bimodule.

Our next result is closely related to the recent works
\cite{IR}, \cite{Br}, \cite{BrGS}, cf. also \cite{KR}.
The theorem below may be viewed as a  `relative' analogue
of \cite[Theorem 3.2]{IR}.

\begin{thm}\label{CM} Let $A$ and $Z$ satisfy the assumptions 
\eqref{zfin}$\mathsf{(i)}$-$\mathsf{(ii)}$  
above, and assume in $\mathsf{(i)\dash(iii)}$ below that $A$ has finite 
 injective dimension. Then, we have:
\vskip 2pt

\vi Each of the functors $\BS,\,\BD_A,$ and $\TTT$, cf. \eqref{T},
preserves the category $D^\fpr(\Lmof{A})$.

\vii The functor $\BS$ is  a Serre functor on  $D^\fpr(\Lmof{A})$
 relative to the subalgebra $Z$,
and the functor $\TTT$ is its  quasi-inverse.

\viii  The
following three conditions are equivalent:
\vskip 2pt

\npb{Calabi-Yau condition \eqref{cy_def} holds for the algebra $A$ (with $d:=\dim(\Spec Z)$);}

\npb{There is an isomorphism of functors $\BS(-)\cong (-)[d]$;}

\npb{There
is an $A$-bimodule isomorphism $A\cong \Hom_Z(A,\CK_Z[d]),$
moreover, $A$ is a maximal Cohen-Macaulay $Z$-module.}
\vskip 2pt

\iv If $A$ has  finite Hochschild dimension
then it has finite injective dimension.
\end{thm} 


\subsection{Proof of Theorem \ref{CM}.}\label{CM_pf}
The finiteness of   injective dimension of $A$ means
that $A$ is a Gorenstein algebra in the sense of  Yekutiely \cite[Example 4.3]{Y}.
Further, the  definition of the complex $\KA$ shows,
as explained in \cite[Proposition 5.7]{Y} and its proof, that  $\KA$ is
a {\em rigid dualizing complex} for the algebra $A$.

Now, \cite[Proposition 8.4]{VdB5} and
 \cite[Proposition 5.11]{Y} imply the following:
\vskip 2pt

\npb{Each of the two objects, $\KA, \,\RHom_{\bimod{A}}(A, A \otimes A)\in D^b(\bimod{A}),$
is a (two-sided) tilting complex; moreover, these bimodule
complexes are invertible and are mutually inverse to each other.}
\smallskip

Any two-sided tilting complex is by definition an
object of $D^\fpr(\Lmod{A})$, with respect to the left
$A$-action, and also an  object of $D^\fpr(\Rmod{A})$, with respect to
the
right $A$-action. In particular, we have $\KA\in D^\fpr(\Lmod{A})$, 
and therefore $\KA\lo_A M\in D^\fpr(\Lmod{A})$ for any
$M\in D^\fpr(\Lmod{A})$. Similarly, we deduce that
 $\RHom_{\bimod{A}}(A, A \otimes A)\lo_A M$ and $\BD_A(M)$
are in $D^\fpr(\Lmod{A}).$
This yields  part (i) of the theorem. 

Further, since
$\KA$ and $\RHom_{\bimod{A}}(A, A \otimes A)$ are mutually inverse bimodule
complexes, we conclude
that $\BS$ and $\TTT$ are mutually quasi-inverse equivalences,
and  part (ii) of the theorem follows.

To prove part (iii), we use  that $\TTT=\BS\inv$ to conclude that
$$
{\mathbb T}(-)\cong (-)[d]\quad
\Longleftrightarrow\quad
\BS\inv(-)\cong(-)[d]\quad\Longleftrightarrow\quad
\BS(-)\cong(-)[d].
$$
Further, using the definitions of the functors ${\mathbb T}$
and $\BS$ we find
\begin{align*}
\RHom_{\bimod{A}}(A, A \otimes A)\lo_A(-)={\mathbb T}(-)\cong (-)[d]\en&\Leftrightarrow\en
\RHom_{\bimod{A}}(A, A \otimes A)\cong A[d]\\
\KA\lo_A(-)=\BS(-)\cong(-)[d]\quad&\Leftrightarrow\quad\KA\cong A[d].
\end{align*}
Thus, combining all the equivalences above, we obtain
\beq{equiv3}
\BS(-)\cong(-)[d]\quad\Longleftrightarrow\quad
\RHom_{\bimod{A}}(A, A \otimes A)\cong A[d]\quad\Longleftrightarrow\quad
\KA[d]\cong A.
\eeq

Observe that the isomorphism in the middle of \eqref{equiv3} is just another
way of writing   equation \eqref{cy_def}, which is the
Calabi-Yau condition. Further, since $\KA=\RHom_Z(A,\CK_Z)$,
 the isomorphism  on the right  of \eqref{equiv3} may be rewriten in the
 form
$$
\Ext_Z^k(A,\CK_Z[d])\cong
\begin{cases} A &\op{if}\enspace k=0\\
0&\op{if}\enspace k>0.
\end{cases}
$$

This formula  is equivalent to saying that
$A$ is a Cohen-Macaulay $Z$-module and, moreover,
$\Hom_Z(A,\CK_Z[-d])\cong A$, which is the
last of the three conditions in the statement
of part (iii) of the theorem.
Thus, \eqref{equiv3} shows that all three conditions of part (iii)
are equivalent.

Assume now  that
$A$ has finite Hochschild dimension. The assumption implies
that 
any finitely generated $A$-module, or $A$-bimodule, has a finite 
resolution by finitely generated projective  $A$-modules.
It follows that 
$D^b(\Lmof{A})=D^\fpr(\Lmof{A})$.
Since , $\KA\in D^b(\Lmof{A})$, we conclude
that $\KA$ is an object of $D^\fpr(\Lmof{A})$

By definition of a dualizing complex,
we have $A=\RHom_A(\KA,\KA)$. Hence, applying \eqref{prin} to $M=\KA,$ 
we conclude that
$A=\RHom_A(\KA,\KA)=\BD_A(\KA)\in D^\fin(\Lmof{A})$.
Similarly, one shows that $A\in D^\fin(\Rmof{A})$.
Thus, the algebra $A$ has  finite injective dimension,  and the theorem is proved.
\qed

\subsection{Proof of Proposition \ref{CYCY}.}\label{CYCYpf}
For any object $\CF\in D^b(\coh  X)$, multiplication by a regular
function $f\in \Gamma(X, \oo_X)$ gives a morphism
$f\cdot : \CF\to\CF$.  This way, any  regular
function $f\in \Gamma(X, \oo_X)$ gives rise to 
an element $z_f\in \DZ(\coh X)$.

We will use the following unpublished result

\begin{lem}[R. Rouquier]\label{rouq_lem} \vi For any smooth,
quasi-projective, and connected variety $X$, the   assignment $f\mto z_f$ induces
an algebra isomorphism $\Gamma(X, \oo_X)\iso \DZ(\coh X)$.

\vii For any homologically smooth algebra $A$ without
zero divisors, the action of the center induces a
natural algebra isomorphism ${\mathcal Z}(A)\iso \DZ(\Lmof{A})$.\qed
\end{lem}

Now, let $\CE\in D^b(\coh X)$ be a tilting object such that
the algebra $A_\CE:=\End\CE$ has finite Hochschild dimension.
Then, according to \cite{BK}, one has  a triangulated equivalence
$D^b(\coh X)\cong D^b(\Lmof{A_\CE})$.
The  equivalence induces an isomorphism between the
centers of the two categories. Thus, we get
algebra isomorphisms
$$\G(X,\oo_X)\cong \DZ(\coh X)\cong\DZ(\Lmof{A_\CE})\cong {\mathcal Z}(A_\CE),$$
where the first isomorphism is due to Lemma \ref{rouq_lem}(i) and the
last isomorphism is due to Lemma \ref{rouq_lem}(ii).
Part (i) of the 
proposition is proved.

Recall next that $X$ is proper over an affine variety.
Thus, there is a proper  morphism
$\pi: X \to Y$, where $Y$ is an affine variety.
We may  assume, without loss of generality,
that $\pi$ is surjective and, moreover, that
$Y$ is smooth (by Noether normalization).
With these assumptions, we put $Z:=\Gamma(Y,\oo_Y),$
a finitely generated commutative algebra with smooth spectrum.

For any object $\CF\in D^b(\coh X),$
the sheaf $\pi_*{{\mathscr H}_{\!}om}(\CF,\CF)$
is $\oo_Y$-coherent, hence
the endomorphism algebra
$\Hom(\CF,\CF)$ is a finitely generated $Z$-module.
Applying this either to $\CF=\oo_X$ or to $\CF=\CE$, we obtain a diagram
of algebra maps
$$
\xymatrix{
Z=\Gamma(Y,\oo_Y)\ar@{^{(}->}[rr]^<>(0.5){\pi^*}&&
\Gamma(X,\oo_X)\ar[r]&{\mathcal Z}(A_\CE)\ar@{^{(}->}[rr]
&&A_\CE,
}
$$
such that $A_\CE$ is a finitely generated $Z$-module.

Next, put $d=\dim X$. On the  category $D^b(\coh X)$, we have the functor
$\BS_X: \CF\mto\CF\o_{\oo_X} \La^d\T^*_X$.
By the  Grothendieck-Serre duality theorem for 
the proper morphism $\pi: X \to Y$, the functor
$\BS_X$ is a Serre functor on
$D^b(\coh X)$ relative  to the subalgebra $Z$.
Transporting this functor via the equivalence
$D^b(\coh X)\cong D^b(\Lmof{A_\CE})$
we get a Serre functor $\BS: D^b(\Lmof{A_\CE})\to D^b(\Lmof{A_\CE}).$

The  definition of a tilting insures that the algebra $A_\CE$ has finite Hochschild
dimension.
 Hence, applying   Theorem \ref{CM}(iv) to $A=A_\CE$ and
$Z=\Gamma(Y,\oo_Y)$ we deduce that this algebra  has finite injective
dimension as well. Furthermore, in this case we have
$D^b(\Lmof{A_\CE})=D^\fpr(\Lmof{A_\CE}).$
Therefore, we are in a position to use Theorem \ref{CM}.
Thus, the equivalent conditions of part (iii) of the theorem
say that $A_\CE$ is a Calabi-Yau algebra iff
the  Serre functor $\BS: D^b(\Lmof{A_\CE})\to D^b(\Lmof{A_\CE})$
is isomorphic to the shift functor, which holds
iff the functor
$\BS_X: \CF\mto\CF\o_{\oo_X} \La^d\T^*_X$
is isomorphic to the shift functor. The latter 
condition is the  Calabi-Yau condition for the manifold $X$,
and we are done.\qed

\begin{proof}[Sketch of proof of Proposition \ref{cy=cy_conj}]
Let $A$ be a homologically smooth, finitely presented,
and compact $A_\infty$-algebra with unit. In this case, using that
$\dim H^\hdot(A)<\infty$, one proves by
induction on the length of the complex that we have
$$ \dim\Ext^\hdot_{D^\fpr(\bimof{A})}(M,N)<\infty,
\qquad\forall M,N\in D^\fpr(\bimof{A}).
$$

Further, the assumptions on $A$ insure that
the proofs of Lemma \ref{serre_weak}
and of Theorem \ref{CM} both apply in our present situation if we
let $Z\cong\k$  be the subalgebra generated by the unit element. This way, one proves that the
the category $D^\fpr(\bimof{A})$ has a Serre
functor $\BS$ and, moreover, there is an isomorphism of
functors $\BS\inv\cong{\mathbb T}$, where the functor 
${\mathbb T}$ is given by  \eqref{T}.
Thus, we conclude as in \eqref{equiv3} that the existence of an isomorphism
$\RHom_{\bimod{A}}(A, A \otimes A)\cong A[d]$
is equivalent to the existence  of an isomorphism
$\RHom_{D^\fpr(\bimof{A})}(A,\k)[d]\cong A,$
in $D^\fpr(\bimof{A})$.
But giving the latter isomorphism amounts to giving
a nondegenerate trace on the  $A_\infty$-algebra
in the sense of Tradler,
cf.  the discussion following Corollary \ref{crep_cor}.

Thus, we have proved the equivalence of the two definitions
of CY algebra.
\end{proof}

{\large\section{{\textbf{Noncommutative Hessian}}}\label{nchess}}

\subsection{A pairing.}
We recall that
the Hessian of a smooth function $f$  on a smooth  manifold $X$
is a symmetric bilinear form,
which is only well-defined, in general, on 
the restriction of the tangent sheaf $\T_X$ to the
critical set of $f$, i.e. on the zero locus of
the 1-form $df$. 

Below, we are going to study a
  noncommutative  analogue of this situation.
Specifically, 
 fix a $\C$-algebra $F$ and
 a cyclic 1-form $\al\in (\Om^1_RF)\br$. We have
the corresponding contraction, resp. reduced contraction, maps
\beq{ii}
\dder F \to F\o F,\,\th\mto i_\th\al,\quad\text{resp.},\en
\dder F\to F,\,\th\mto \bi_\th\al.
\eeq

Following \S\ref{cy3}, we consider the algebra $A:=F/\imm$
(we are now additionally assuming, for
simplicity of notation,  that $R=\C$). In the special case where $\al=d\Phi$, the algebra
$A$ plays the role of 
a  noncommutative  analogue of the critical set of a smooth function.
A  noncommutative  analogue of the Hessian of that function
is played by the  $\C$-bilinear pairing
\beq{pairingF}
\dder F \times \dder F\to F\o F\stackrel{\pr_A}\onto A\o A,
\quad \xi\times\th\mto
\pr_A\th(\bi_\xi\al)=\pr_Ai_\th(d\ccirc\bi_\xi\al),
\eeq
This analogy is supported by Proposition \ref{prop5}  below.

To state the Proposition,  introduce  the notation
$(u\o v)\y=v\o u$. Also, define a graded $A$-algebra
$\Om^\hdot(F|A):=T_A^\hdot(\Om^1(F|A))$ and put
$$\Om^2(F|A)\br=\Om^2(F|A)/[A,
\Om^2(F|A)]=\Om^1(F|A)\o_{A^e}\Om^1(F|A).$$

\begin{prop}\label{prop5} Assume that the algebra $F$ is smooth  and the 1-form $\al\in
 (\Om^1F)\br$ is such that $d\al=0$ holds in $\Om^2(F|A)\br.$ Then,

\vi The pairing  in \eqref{pairingF}  descends
 to a well-difined $\C$-bilinear map
\beq{pairing}
{\mathbb H}: \DPA \o \DPA\to A\o A,
\quad \xi\times\th\mto
{\mathbb H}(\xi,\th):=i_\th(d\ccirc\bi_\xi\al).
\eeq
 
\vii The  following {\em{Hessian symmetry property}} holds:
\beq{hessym}
{\mathbb
H}(\th,\xi)=({\mathbb H}(\xi,\th))\y,\qquad\forall\xi,\th\in\DPA.
\eeq
This property is also  equivalent to the condition that  the map
 $d\ccirc\bi$ in  \eqref{map} be {\sl selfadjoint}, i.e., that
one has $(d\ccirc\bi)^\vee=d\ccirc\bi.$ 
\end{prop}

The proof of the proposition is based on the following
\begin{prop}\label{symmetric} \vi For any algebra $F$ and 
  $\om\in(\Om^2F)\br$, in $F\o F,$
we have
$$i_\th\bi_\P\om=- (i_\P\bi_\th\om)\y,\quad\forall\th,\P\in\dder F.$$

\vii  For any  $\al\in (\Om^1F)\br$ such that   $d\al=0$ holds in $\Om^2(F|A)\br,$ we have
$$
i_\th(d\bi_\xi\al)
=(i_\xi(d\bi_\th\al))\y\en\oper{holds\en in}\en A\o A,\quad
\forall \xi,\th\in\DPA.
$$
\end{prop}
\begin{proof} To prove (i), we may assume $\om=\al\be$ for some $\al,\be\in \Om^1F$. We compute
$$i_\th\bi_\P(\al\be)=
\bi_\th\bigl(i''_\P\al\cd\be\cd i'_\P\al-i''_\P\be\cd\al\cd
i'_\P\be\bigl)=
i''_\P\al\cd i'_\th\be\,\o\, i''_\th\be\cd i'_\P\al-
i''_\P\be\cd i'_\th\al\,\o\, i''_\th\al\cd i'_\P\be.
$$
Exchanging the roles of $\th$ and $\P$ and comparing
the summands in the resulting expressions yields the required identity.

To prove part (ii), write $\bl_\xi=d\ccirc\bi_\xi+\bi_\xi\ccirc d$ for  {\em reduced Lie
derivative}
introduced in \cite{CBEG}. We use formula \cite[(A.6)]{VdB4}
which says that, for any $\th,\xi\in\dder
 F$, there exist certain
elements $\nu\o t\in(\dder
 F) \o F$ and $s\o \sigma\in F\o(\dder
 F)$ such that
one has
\beq{identity}
i_\th\bl_\xi\al-(L_\xi\bi_\th\al)\y=
\bi_\nu\al\o t+
s\o\bi_\sigma\al, \qquad\forall\al\in(\Om^1F)\br.
\eeq
(Van den Bergh used the notation
$\bll\th,\xi\brr'_l\o\bll\th,\xi\brr''_l$
for our $\nu\o t$, resp. $\bll\th,\xi\brr'_r\o\bll\th,\xi\brr''_r$
for our  $s\o \sigma$, where $\bll\th,\xi\brr$ stands for the
double bracket introduced in \cite{VdB4}).

Let $I=\imm,$ a two-sided ideal in $F$.
By definition, we have $\bi_\nu\al,\,\bi_\sigma\al\in
  I $.
Therefore, the left hand side of \eqref{identity} vanishes
modulo $ I \o F+F\o  I $.
Hence, 
 using Cartan's formula, in $A\o A=(F\o F)/( I \o F+F\o  I )$,
 we find
\begin{align}\label{schet}
i_\th(d\bi_\xi\al)-(i_\xi(d\bi_\th\al))\y
=\big(i_\th\bl_\xi\al-(L_\xi\bi_\th\al)\y\big)&-\big(i_\th\bi_\xi
 d\al-
(d i_\xi\bi_\th\al)\y\big)\nonumber\\
=&-i_\th\bi_\xi
 d\al+(d i_\xi\bi_\th\al)\y.
\end{align}

Since $\al$ is a 1-form, we have $i_\xi\bi_\th\al=0$,
so the second term in the last line of \eqref{schet} 
vanishes. Also, since $d\al=0$ in $\Om^2(F|A)\br,$
one can write $d\al$ as a sum of terms of the
form $f\,dx\,dy$, where $f\in  I $ and
 $x,y\in F$. Thus, the 1-form $\bi_\xi 
 d\al$ is a sum of the following expressions:
$$
\bi_\xi(f\,dx\,dy)=
\xi''(x)\,dy\cd f\cd \xi'(x)-\xi''(y)\cd f\,dx\cd\xi'(y)
\in (\Om^1F)\cd  I +  I \cd (\Om^1F).$$
We deduce that
 $i_\th\bi_\xi(f\,dx\,dy)$
  vanishes modulo $ I \o F+F\o  I $.
So, the  term 
 $i_\th\bi_\xi d\al$ in the last line of \eqref{schet} 
vanishes as well, and we are done.
\end{proof}

\begin{proof}[Proof of Proposition \ref{prop5}] First of all, equation
\eqref{hessym} follows from  Proposition \ref{symmetric}(ii).

Next, it is immediate from 
definitions that the pairing \eqref{pairingF} descends to a pairing
$\DPA\times \dder F \to A\o A$. Part (i) of Proposition \ref{prop5}
now follows by the symmetry, due to \eqref{hessym}.

To prove (ii), consider the morphism
$(d\bi)^\vee: \DPA\to\DPA^\vee,
\xi\mto (d\bi)^\vee(\xi)$.
From the definition of a dual morphism one finds that the element
$ (d\bi)^\vee(\xi)$ is given by the map
$$(d\bi)^\vee(\xi): \DPA\to A\o A,\quad\th\mto \HH(\xi,\th).$$
Therefore the equation $(d\bi)^\vee=d\bi$ is equivalent
to the symmetry property \eqref{hessym}. We remark that the flip involved in
\eqref{hessym} does not contradict the equation $(d\bi)^\vee=d\bi$.
The appearence of that flip is due to the fact that
the bimodule structure on double derivations is induced
by the {\em inner} bimodule structure on $A\o A$;
to compare the inner and outer structurs, one has to
use the flip  $A\ino A \iso A\out A,\, a\o b\mto b\o a,$
which is an isomorphism of $A$-bimodules.
\end{proof}

\begin{rem} The symmetry property \eqref{hessym} implies, in particular,
that
for any $\xi\in\DPA,$ the map
${\mathbb H}(\xi,-): \DPA\to A\o A,\,\th\mto{\mathbb H}(\xi,\th)$
is an $A$-bimodule map, i.e., one has
$${\mathbb H}(\xi,a\cd\th\cd b)=a\cd {\mathbb H}(\xi,\th)\cd b,\qquad
\forall\xi,\th\in\DPA,\,a,b\in A.$$
\end{rem}

\begin{rem}\label{FFA}
Similarly  to Lemma \ref{FF}(ii), one proves that, for a {\em smooth}
algebra $F$, 
the assignment $\be\mto\bi^\be$ 
yields a bijection
$$\Om^2(F|A)\br\iso
\{f\in\Hom_{\bimof{A}}\big(\dder(F|A),\,\Om^1(F|A)\big)\en\big|\en
f=f^\vee\}.
$$
\end{rem}

\subsection{}\label{cy3_pf} From now on, we fix a  {\em smooth} algebra
$F$, a 1-form $\al\in (\Om^1F)\br$, and put  $A:=F/\imm$.
The goal of this section is to prove
the following result.

\begin{thm}\label{HHthm} \vi Assume that in $(\Om^2(F|A))\br$ one has
$d\al=0$. Then,  there is a canonical exact sequence of $A$-bimodules
\beq{chain2}
0\to\dder A\stackrel{j}\too\DPA
\stackrel{d\ccirc\bi}\too
\Om^1(F|A)\stackrel{p}\too\Om^1A\to 0.
\eeq 

\vii Assume, in addition, that $\bi_{\De}\al\in(\imm)^2$.
Then \eqref{chain2} gives rise to
the following exact sequence of $A$-bimodules
$$0\to H^1(A, A\o A)\stackrel{\wb j}\too H^1(F,A\o A)\stackrel{\overline{d\ccirc \bi}}\too
H_1(F,A\o A)\stackrel{\wb p}\too
H_1(A,A\o A)\to 0.
$$
\end{thm}

\begin{proof} 
We recall that,
for a {\em smooth} algebra $F$,  both $\Der(F|A)$ and $\dder(F|A)$
are finitely generated projective $A$-bimodules.
Furthermore, there are canonical $A$-bimodule  isomorphisms
\begin{align}
 &\Der(F,A\o A)\cong A\o_F\dder F\o_F A=\dder(F|A)\quad\text{and}\label{DD1}\\
&
(\dder(F|A))^\vee\cong A\o_F(\dder F)^\vee\o_FA=
\Om^1(F|A).\label{DD2}
\end{align}

Next, we put $I:=\imm.$
One has an obvious identification
\beq{obv}\Der(A,A\o A)\cong\{\th\in \Der(F,A\o A)\mid \th|_I=0\}.
\eeq
Thus, using \eqref{DD1} and
\eqref{map}, we obtain a chain of  $A$-bimodule maps
$$
0\to\Der(A,A\o A)\stackrel{j}\to\Der(F,A\o A)=\dder(F|A)
\stackrel{d\ccirc\bi^\al}\too\Om^1(F|A)\stackrel{p_{F|A}}\too\Om^1A\to 0,
$$
where the map $j$ on the left is the natural imbedding.
We let
the chain above  be the sequence \eqref{chain2} in  part (i) of the Theorem.

\step{1.}
We prove that  the chain of maps in \eqref{chain2} is an exact sequence.

Injectivity of the map $j$ is clear from \eqref{obv}.
Further, the short exact  sequence \eqref{CQ} yields the
surjectivity of the last map $p$ and it also implies
that  the sequence \eqref{chain2} is exact
at the term $\Om^1(F|A).$

It remains to prove the exactness of \eqref{chain2} at the term
$\dder(F|A).$
To this end, observe first that an element $\xi\in\dder(F|A)$
 belongs to the image of the imbedding
$j$ in \eqref{chain2} if and only if we have  $\xi|_I=0$.  Hence, using
Proposition
\ref{prop5}(i) and the definition of the ideal $I$ we see that
 $\xi|_I=0$ holds  if and only if, for all
$\th\in\Der(F,A\o A)$, we have  $\xi(\bi_\th\al)=0$.
We can rewrite this, in terms of the
pairing \eqref{pairing}, as an equation $\HH(\xi,-)=0.$

Further, we have seen in the course of the proof of
Proposition
\ref{prop5}(ii), that the map $\xi\mto\HH(\xi,-)$
may be identified with the map $(d\ccirc\bi)^\vee$.
Moreover, since $(d\ccirc\bi)^\vee=d\ccirc\bi$, we conclude
$$\HH(\xi,-)=0 \quad\Leftrightarrow\quad(d\ccirc\bi)^\vee(\xi)=0
 \quad\Leftrightarrow\quad d\ccirc\bi(\xi)=0.
$$
Thus, we have proved that $\im(j)=\ker(d\ccirc\bi)$,
and  part (i) of the theorem follows.

\step{2.} 
Given an algebra $B$ 
write $M\br:=M/[B,M]=M\o_{B^e}B$ for the commutator quotient space of a $B$-bimodule $M$.
Further, let $\Om^1(B,M):=M\o_B\Om^1 B.$
This is a $B$-bimodule whose elements are sums
 of the form
$m\,db$, where $m\in M,b\in B.$
We form  the commutator quotient
$$\Om^1(B,M)\br:=(M\o_B\Om^1 B)\br=M\o_{B^e}\Om^1B.$$

For any $\th\in \dder B$, there is a natural reduced contraction map
$$\bi_\th^M:\ \Om^1(B,M)\br\too M,\quad 
m\,df\mto \th''(f)\cd m\cd \th'(f),\quad
\forall m\in M,f\in B.
$$
In particular, for  any 1-form $m\,df\in \Om^1(B,M)\br,$ we
get $\bi^M_\De(m\,df)=[m,f]\in M$.

According to \cite{CQ}, \cite{Gi3}, for any $B$-bimodule $M$,
there is a natural vector space
isomorphism
\beq{H1}
H_1(B,M)\cong\{\la\in\Om^1(B,M)\br\en\big|\en \bi_\De^M\la=0\}.
\eeq

Now, let an algebra $A$ be a quotient of $B$.
In the special case $M:=A\o A$,
there are natural bijections
$$\Om^1(B,A\o A)\br=
(A\o A)\o_{B^e}\Om^1 B=A\o_B\Om^1 B\o_B A=\Om^1(B|A).
$$
Furthermore, we have the following  equality of two contraction  maps
$$
\bi_\th^{A\o A}=(i_\th)\y: \;\Om^1(B|A)=A\o_B\Om^1\o_B A\too A\o A.
$$
Thus, we deduce a canonical
isomorphism
\beq{canB}
H_1(B,A\o A)\cong\{\la\in
\Om^1(B|A)\mid i_\De\la=0\}.
\eeq

\step{3.} We now return to our setting and let  $B$ be either the algebra $F$ or
its quotient
$A=F/\imm$.

We claim that  any element
contained in the image of the map
$d\ccirc\bi^\al$, in \eqref{chain2}, is annihilated by 
$i_\De$.
To see this, we use the symmetry property from Proposition
\ref{symmetric}(ii).
Thus, for any $\th\in\dder(F|A)$, we get
$i_\De d \bi_\th\al=(i_\th d\bi_\De\al)\y=0,$
since $\bi_\De\al=0$. 

Using \eqref{canB}, we conclude that the map
$d\bi: \Der(F,A\o A)\to \Om^1(F|A)$ may be interpreted as a map
$$d\bi:\  \Der(F,A\o A)\to \{\la\in
\Om^1(F|A)\mid i_\De\la=0\}=H_1(F,A\o A).
$$

Similarly, we  apply \eqref{canB} for $B:=A$ and use the following  commutativity  diagram
\beq{tri}
\xymatrix{
\Om^1(F|A)
\ar[rr]^<>(0.5){p}\ar[dr]^<>(0.5){\bi_\De}&&\Om^1A
\ar[dl]_<>(0.5){\bi_\De}\\
&A\o A&
}
\eeq
We see that the map $p$ in \eqref{chain2} descends to a map
$\wb p: H_1(F,A\o A)\to H_1(A,A\o A)=\{\la\in
\Om^1A\mid i_\De\la=0\}.$

\step{4.} The map $j$ in  \eqref{chain2} clearly
takes inner derivations to inner derivations,
specifically, we have $j(\Inn(A,A\o A))=\Inn(F,A\o A).$
Therefore, the imbedding $j$ descends to an {\em injective}
map $\wb j: H^1(A,A\o A)\into H^1(F,A\o A).$ 
Furthermore,
since $\Inn(F,A\o A)$ is contained in the image of $j$,
by Step 1 we conclude that inner derivations are killed
by the map $d\bi$. Thus, the map $d\bi$ descends to 
a map $\wb{d\bi}: H^1(F,A\o A)\to H_1(F,A\o A).$

Thus, the sequence of maps in \eqref{chain2}
gives rise to a  complex of $A$-bimodules,
as in the statement of part (ii) of the Theorem.

It remains to show that the constructed complex is an exact sequence.
Injectivity of the map
$\wb j$ has been already mentioned above.
Surjectivity of the map $\wb p$ in our complex follows from
the surjectivity of the map
$p$  in  \eqref{chain2} and commutativity of  diagram \eqref{tri}.
Finally, the exactness at the terms $H^1(F,A\o A)$ 
and $H_1(F,A\o A)$ is immediate,
in view of  isomorphism \eqref{canB}, from the exactness of  \eqref{chain2}
at the corresponding terms.
\end{proof}

{\large\section{{\textbf{Some homological algebra}}}\label{selfdual_sec}}
\subsection{Bimodule resolutions of CY algebras}\label{self_pf}
Throughout this section, we fix an $R$-algebra $A$ and write
$\Ext^\hdot=\Ext^\hdot_{\bimod{A}},\,\Hom=\Hom_{\bimod{A}},$ etc.
We begin with a standard 

\begin{lem}\label{st} Let $A$ be a coherent, finitely presented
$R$-algebra of finite Hochschild dimension $n+2$. Then, $A$ is friendly, and
we have

\vi Any  finitely presented
$A$-bimodule has a bounded resolution 
by finitely generated projective $A$-bimodules.

\vii  Assume that $A:=F/I,$ where $F$ is a smooth, finitely generated
$R$-algebra,
and  $I\sset F$ is a two-sided
finitely generated ideal. 

Then the algebra $A$, viewed as an $A$-bimodule,
has a resolution of the form
\beq{resa}
{\mathsf P}:\en 0\too P_{n+1}\too P_n\stackrel{d_n}\too\ldots
\stackrel{d_2}\too 
P_1\stackrel{d_1}\too P_0
\stackrel{d_0}\too A\o_RA\en(\;\stackrel{_{\oper{mult}}}\onto\; A),
\eeq
where  $\,P_0=\Om^1_R(F|A),$ and
$P_r, \,r=1,\ldots,n+1,$ are finitely generated projective $A$-bimodules.
\end{lem}
\begin{proof} 
Observe first that since $A$ is a  finitely
generated algebra,
$\Om^1_RA$ is a finitely generated $A$-bimodule.
Therefore, there exists  a cover $P_0\onto \Om^1_RA,$  by
a   finite rank
free $A$-bimodule $P_0$. Hence,
we see that the multiplication map provides 
a  projective presentation
$P_0\to A\o_RA\onto A.$ 
If $A=F/I$, as in part (ii) of the lemma,
then the map $\Om^1_R(F|A)\to \Om^1_RA$ is surjective
and, moreover, $\Om^1_R(F|A)$ is a finitely
generated projective $A$-bimodule,
so one can take $P_0:=\Om^1_R(F|A)$.

Further, since $A$ is coherent, 
the kernel of the map $P_0\to A\o_RA$ is a finitely
generated. Hence,
there exists a  finite rank
free $A$-bimodule $P_1$ and an exact
sequence $P_1\to P_0\to A\o_RA$. Continuing in this fashion,
one obtains a resolution of the form
\eqref{resa},
where $P_i$ is a finite rank
free $A$-bimodule for all $i=0,\ldots,n,$
and where $P_{n+1}:=\ker d_n,$
is a  finitely
generated $A$-bimodule.

Now, using  standard  long exact sequences of
Ext-groups, we deduce
$$\Ext^1(P_{n+1},N)=\Ext^{n+3}(A,N)=0,\quad
\forall N\in\bimod{A},
$$
where the equality on the right is due to 
the  bound on the Hochschild dimension of $A$.
It follows that $P_{n+1}$ is a projective $A$-bimodule. Thus,
 \eqref{resa} yields a length $n$ resolution of $A$
 by projective finitely
generated $A$-bimodules.

Since $A$ is free as either right or left $A$-module, the resolution
in \eqref{resa}
gives a {\em split} exact 
 sequence of  either right or left $A$-modules. Hence, tensoring with an
$A$-bimodule $M$ yields an exact sequence
$\ldots\to P_1\o_AM\to P_0\o_AM\to A\o_RM\onto M.$
Each term $P_i\o_AM$ is a direct summand of a direct
sum of free left $A$-modules of the form $A\o M$, hence it is
projective as a left $A$-module.
Similarly, we tensor the resulting exact sequence above
with our resolution of $A$,  viewed now  as a split
 exact sequence of  left $A$-modules. This way we obtain
a bicomplex with terms of the form
$P_i\o_AM\o_AP_j$. The total complex associated to this
bicomplex provides a bounded resolution of $M$ by  projective
$A$-bimodules. It follows that there exists an integer $m=m(M)\gg 0$
such that 
\beq{vanishX}
\Ext^{m+1}(M,N)=0,\quad
\forall N\in\bimod{A}.
\eeq

It remains to show that, if $M$ is 
a {\em finitely presented} $A$-bimodule, then one can find  a bounded
resolution of $M$ by {\em finitely generated} projective
$A$-bimodules. To do this, start with a presentation
$P_1\to P_0\onto M,$
by  finite rank
free $A$-bimodules. Then, use coherence
of $A$ and proceed as at the beginning of the proof
to obtain a resolution $Q_\idot$ of $M$ by finitely generated
 projective
$A$-bimodules.  Now, the vanishing
in \eqref{vanishX} implies that 
the kernel of the map $Q_m\to Q_{m-1},$
in that resolution, is 
projective, and we are done.
\end{proof}

\begin{lem}\label{imj} Let $A$ be a friendly  Calabi-Yau
algebra of dimension $n+2$. Then,

\vi A volume $\pi$ on $A$ induces an isomorphism
$\dis
\pi^R:\ A^R\iso\Ext^{n+2}(A,(A\o A)^R);
$
\vskip 1pt

\vii There exists a resolution of $A$ by finitely generated,
projective $A$-bimodules $
\SP:$
$$
\xymatrix{
0\to(A\o A)^R\ar[r]^<>(0.5){\bj}&
P_n\ar[r]^<>(0.5){d_n}&
P_{n-1}\ar[r]^<>(0.5){d_{n-1}}&\ldots\ar[r]^<>(0.5){d_2}&P_1\ar[r]^<>(0.5){d_1}&
P_0\ar[r]^<>(0.5){d_0}&
A\o_RA\ar@{->>}[r]^<>(0.5){\oper{mult}}&A
}
$$
such that the class of this exact sequence represents
the class $\pi^R(1)$.
\vskip 1pt

\viii  The image of $\dis\,\bj^\vee: P_n^\vee\to ((A\o A)^R)^\vee=A\o_R A$,
 the  dual of the imbedding $\bj$ in the above resolution, equals
 $\bj^\vee(P_n^\vee)= \Om^1_RA\sset A\o_R A$. 

\iv If $A=F/I$, as in Lemma \ref{st}(ii), then the resolution
 yields an exact sequence
$$\xymatrix{
0\to(A\o A)^R\ar[r]^<>(0.5){\bj}&
P_n\ar[r]^<>(0.5){d_n}&
P_{n-1}\ar[r]^<>(0.5){d_{n-1}}&\ldots\ar[r]^<>(0.5){d_2}&P_1\ar[r]^<>(0.5){\bp}&
I/I^2\ar[r]&0.
}
$$
\end{lem}
\proof Part (i) is clear. 
An $A$-bimodule isomorphism
 $A\iso\Ext^{n+2}_{\bimod{A}}(A,A\o A)$ restricts
to an isomorphism $A^R\iso\Ext^{n+2}_{\bimod{A}}(A,(A\o A)^R).$ 
This proves (i).

To prove (ii), we
start with a resolution of the form, cf. \eqref{resa},
\beq{KKK} 0\to K\too  P_n\stackrel{d_n}\too\ldots P_2\stackrel{d_2}\too P_1
\stackrel{d_1}\too A\o_RA\en(\stackrel{\eps}\onto A),
\eeq
where $P_i$ is a finite rank
free $A$-bimodule for all $i=1,\ldots,n-1,$
and where $K:=\Ker d_n.$ 

For any object $X$, using standard long exact sequences
of Ext-groups one finds that there are canonical isomorphisms
\begin{align}\label{ext_red}
\Ext^{i}(K,X)&\iso\Ext^{i+n}(\Om^1_RA,X)\iso
\Ext^{i+n+1}(A,X),\quad\forall i\geq 1,\en\text{and}\\
&\text{a surjection}\quad
\Hom(K,X)\onto\Ext^{n}(\Om^1_RA,X)=\Ext^{n+1}(A,X),\en\text{for $i=0$}.\nonumber
\end{align}

Thus,
we may identify the element
 $\pi^R(1)\in \Ext^{n+2}(A,(A\o A)^R)$,  with a class
in $\Ext^1(K,(A\o A)^R)$, that is, with an
extension
\beq{ANK}\wt{\pi}^R:\ 0\to(A\o A)^R\stackrel{\bj}\too N\stackrel{v}\too K\too 0.
\eeq

Splicing this extension with \eqref{KKK} yields an extension as
in   the displayed formula
in part (ii) of the Lemma that, by construction, represents the class of
 the image of $1\in A^R$ under the
 isomorphism $A^R\iso\Ext^{n+2}_{\bimod{A}}(A,(A\o A)^R).$ 
Furthermore, all the terms
$P_i$, except possibly $P_n$ are finitely generated projective
$A$-bimodules.

Let $\wt\SP$ be the push-out of $\SP$ via the
imbedding $(A\o A)^R\into A\o A$. Since $(A\o A)^R$ is a direct summand
in $A\o A$,
the latter imbedding is split. Therefore, $\SP$ is a direct
 summand of $\wt\SP$, and
the end of the exact sequence  $\wt\SP$ gives an extension $\wt\pi: A\o A\into \wt N\onto
\wt K$, such that \eqref{ANK} is its direct summand.

For any object $X$, we have the corresponding
 long exact sequence of Ext-groups
\beq{extx}
\Hom(\wt N,X)\stackrel{\wt\bj}\to\Hom(A\o A,X)\stackrel{\wt\pa}\to \Ext^1(\wt K,X)\stackrel{v}\to
\Ext^1(\wt N,X)\to\Ext^1(A\o A,X)=0,
\eeq
where $\wt\pa$ is the connecting homomorphism. We notice
that all the above maps are morphisms  of $A$-bimodules with respect to the 
bimodule structure induced by the inner 
bimodule structure on $A\o A$.

For  $X=A\o A$, the boundary map $\pa$ in \eqref{extx}
is a map $\wt\pa: \Hom(A\o A,A\o A)\to \Ext^1(K,A\o A).$
It is known that one has $\wt\pa(\Id)=\wt\pi$.
Therefore, by definition of $\wt \pi$,  the composite map
$$
m: \
A\o A\iso \Hom(A\o A,A\o A)\stackrel{\wt\pa}\too \Ext^1(\wt K,A\o A)\iso
\Ext^3(A,A\o A)\iso A,
$$
is such that  $m(1\o 1)=1$.
The map $m$ being an $A$-bimodule
morphism, for any $a,b\in A$, we deduce
$m(a\o b)=m(a\cdot(1\o 1)\cdot b)=  a\cdot m(1\o 1)\cdot b=a\cdot b$.
Thus, $m$ is nothing but the multiplication map
and we get $\Ker(\wt\pa)=\Ker(m)=\Om^1A$. Similarly,
restricting everything to  the direct summand $(A\o A)^R\sset A\o A$,
we obtain from $\wt\pa$ a similar
connecting homomorphism $\pa$, and we deduce similarly that $\im(\bj^\vee)=
\Ker(\pa)=\Om^1_RA.$

Next we observe that the map $m$ above is clearly surjective.
We deduce that
the connecting homomorphism $\wt\pa$ in \eqref{extx} is surjective
as well. Hence, it follows from that exact sequence
that $\Ext^1(\wt N,A\o A)=0$. This implies that for any projective
$A$-bimodule $P$ we have $\Ext^1(\wt N,P)=0$.

We know that $\wt N$ has a finite projective resolution, by Lemma \ref{st}.
Now, a routine induction on the length of projective
resolution shows that the vanishing of  $\Ext^1(\wt N,P)=0$
for all projective modules $P$ implies that $\wt N$
 is itself  projective. Therefore, $N=P_n$, being a direct summand of
$\wt N$, is projective as well.
This completes the proof of parts (ii), (iii).

Part (iv) follows from Lemma \ref{st}(ii),
and the cotangent
sequence \eqref{CQ}.\qed

\begin{thm}\label{ddpf} Let $A$ be a  Calabi-Yau
$R$-algebra of dimension $d\geq 2$ 
and let  $\pi\in
\Ext^d(A,A\o A)$ be a volume element.
Then, \vskip 2pt

\vi There exists a commutative diagram of the form
{\beq{PP}
\xymatrix{
\SP:\ar[d]^<>(0.5){f}\\
\SP^\vee:}
\xymatrix{
0\ar[r]& (A\o A)^R\ar[d]^<>(0.5){f_d}\ar[r]^<>(0.5){d_d}&
P_{d-1}\ar[r]^<>(0.5){d_{d-1}}\ar[d]^<>(0.5){f_{d-1}}&\ldots\ar[r]^<>(0.5){d_2}&
P_1\ar[r]^<>(0.5){d_1}\ar[d]^<>(0.5){f_1}&
P_0\ar[d]^<>(0.5){f_0}\ar@{->>}[r]&A\\
0\ar[r]& P_0^\vee\ar[r]^<>(0.5){d_1^\vee}&
 P_1^\vee\ar[r]^<>(0.5){d_2^\vee}& \ldots\ar[r]^<>(0.5){d_{d-1}^\vee}&
P_{d-1}^\vee\ar[r]^<>(0.5){d_d^\vee}&A\o_R A\ar@{->>}[r]&A,
}
\eeq}

such that  each of the rows  in \eqref{PP}  has the following
properties:
\vskip 2pt

\npb{it is a resolution of $A$ by  finitely generated projective
$A$-bimodules, and}

\npb{it represents the class $\pi$.}
\vskip 2pt

\vii If, in addition, $A$ is friendly then one may choose $\SP$ above so that
$P_0=A\o_R A$. Furthermore, one can construct a {\em selfdual}
diagram $f: \SP\to\SP^\vee$, i.e., such that we have
$f^\vee_i=f_{d-i},\,\forall i=0,1,\ldots,d;$ 
so, diagram \eqref{PP} reads
\beq{PPA}
\xymatrix{
0\ar[r]& (A\o A)^R\ar@{=}[d]^<>(0.5){\Id}\ar[r]^<>(0.5){d_d}&
P_{d-1}\ar[r]^<>(0.5){d_{d-1}}\ar[d]^<>(0.5){f_{d-1}}&\ldots\ar[r]^<>(0.5){d_2}&
P_1\ar[r]^<>(0.5){d_1}\ar[d]^<>(0.5){f_1}&
A\o_RA\ar@{=}[d]^<>(0.5){\Id}\ar@{->>}[r]^<>(0.5){\oper{mult}}&A\\
0\ar[r]& (A\o A)^R\ar[r]^<>(0.5){d_1^\vee}&
 P_1^\vee\ar[r]^<>(0.5){d_2^\vee}& \ldots\ar[r]^<>(0.5){d_{d-1}^\vee}&
P_{d-1}^\vee\ar[r]^<>(0.5){d_d^\vee}&A\o_R A\ar@{->>}[r]^<>(0.5){\oper{mult}}&A,
}
\eeq

\end{thm}

\begin{proof}[Proof of part $\mathsf{(i)}$] Choose
a resolution  by finitely generated
projective $A$-bimodules
$$\SP:\ 0\to P_{d+1}\to P_d\to\ldots\to P_1\to P_0\to A\to 0.
$$ 
 
Arguing as in the proof of 
Lemma \ref{st}, we may insure that $P_{d+1}=A\o A$.
(the argument does not depend on the assumption 
 that $A$ be coherent).
One may use the resolution $\SP$ 
to compute $\Ext^\hdot(A, A\o A).$ 
Thus, we form the dual sequence  $\SP^\vee$. Using that
$\Ext^i(A, A\o A)=0$ for all $i<d$ we deduce
the  exact sequence in the bottom row of diagram \eqref{PP}.
The argument in the proof  of
Lemma  \ref{imj}(iii) shows that the latter exact sequence
represents the class $\pi$ as well.

By abstract nonsense,  there exists a morphism $f:\SP\to\SP^\vee$
between the two resolutions above. We refer to
\cite{Bo} for more details.

The proof of part (ii) of the theorem in the special
case of CY algebras of dimension $d=3$ will be given in the next section. The proof
of the case $d>3$ is similar and will be omited.
\end{proof}

\subsection{Proof of Theorem \ref{ddpf}(ii) for CY algebras of dimension 3}\label{splice_sec}
 From now on, 
we let $A$ be a friendly algebra that satisfies condition
\eqref{cy_def} for $d=3$, and $\pi$ is a volume element.
 To simplify the notation, we assume
our ground ring to be $R=\C$.

Diagram \eqref{PPA} will be  constructed
by splicing together two shorter diagrams,
dual to each other,  along the  {\em self-dual} morphism $g=g\ve$, as
depicted below
\beq{splice}
\xymatrix{
0\ar[r]& A\o A\ar[r]^<>(0.5)\bj\ar@{=}[d]^<>(0.5){\Id}&
N\ar@{->>}[rr]^<>(0.5){v}\ar[d]^<>(0.5){f\ve}&&K\ar[d]^<>(0.5){g=g\ve}
\ar@{^{(}->}[rr]^<>(0.5){u}&&M\ar[d]^<>(0.5){f}
\ar@{->>}[r]^<>(0.5){\pp}&\Om^1A\ar[r]\ar@{=}[d]^<>(0.5){\Id}&0\\
0\ar[r]& A\o A\ar[r]^<>(0.5){\pp\ve}&M\ve\ar@{->>}[rr]^<>(0.5){u\ve}&&K\ve
\ar@{^{(}->}[rr]^<>(0.5){v\ve}&&N\ve\ar[r]^<>(0.5){\bj\ve}&\Om^1A\ar[r]&0.
}
\eeq
Here, $M$ and $N$ are
 finitely generated,
projective $A$-bimodules,
and  each row is composed by two
 short exact sequences.
Furthermore,
the vertical  morphism $g=g\ve$ in \eqref{splice}
 is {\em selfadjoint}, and splicing along this morphism
yields a diagram whose rows represent the class $\pi$, each.

To begin,  choose a finitely generated  projective  cover of $M\onto\Om^1A.$
Thus, there is an $A$-bimodule
 extension 
\beq{KMO}0\to K\stackrel{u}\too M\stackrel{\pp}\too \Om^1A\to0.
\eeq

Apply the functor $(-)^\vee$ to the  short exact sequence 
\eqref{ANK}. Using  that $\im(\bj^\vee)=\Om^1A\sset A\o A$
by Lemma \ref{imj}(iii), we obtain a short exact sequence
\beq{KNO}
0\to K^\vee\stackrel{v^\vee}\tooo N^\vee \stackrel{\bj^\vee}\tooo\Om^1A\too 0.
\eeq

\begin{lem}\label{gg_lem}
There exists a commutative diagram
$$
\xymatrix{
0\ar[r]&K\ar[d]^<>(0.5){g}
\ar@{^{(}->}[rr]^<>(0.5){u}&&M\ar[d]^<>(0.5){f}
\ar@{->>}[r]^<>(0.5){\pp}&\Om^1A\ar[r]\ar@{=}[d]^<>(0.5){\Id}&0\\
0\ar[r]& K\ve
\ar@{^{(}->}[rr]^<>(0.5){v\ve}&&N\ve\ar@{->>}[r]^<>(0.5){\bj\ve}&\Om^1A\ar[r]&0.
}
$$
such that the map $g$ is selfdual, i.e., we have $g\ve=g.$ 
\end{lem}

\begin{proof} 
Similarly to what we have done above the statement of the Lemma,
 we apply the functor $(-)^\vee$ to
 \eqref{KMO}.
We get an
exact sequence
$$
0\to (\Om^1A)^\vee\stackrel{\pp^\vee}\too M^\vee\stackrel{u^\vee}\too K^\vee\too
\Ext^1(\Om^1A, A\o A)=\Ext^2(A,A\o A)=0.
$$
Observe  that since $\Ext^0(A,A\o A)=\Ext^1(A,A\o A)=0$,
we have a canonical isomorphism
$\dis A\o A\stackrel{\ad}\iso \dder A=(\Om^1A)^\vee.$
Therefore, the  sequence above reads
\beq{AMK}
0\too A\o A\stackrel{\pp^\vee\ccirc\ad}\tooo
 M^\vee\stackrel{u^\vee}\tooo K^\vee\too 0.
\eeq

Next, we splice together short exact sequences
\eqref{ANK}, \eqref{KMO}, and the fundamental
short exact sequence $\Om^1A\into A\o A\onto A$. We obtain
the following complex 
\beq{resolution}
0\to A\o A\stackrel{\bj}\too N\stackrel{d}\too M\stackrel{\pp}\too A\o A\to 0,
\eeq
where $d=v\ccirc u$. This complex  provides a
 projective $A$-bimodule resolution of $A$.

Dually, we splice  together short exact sequences
\eqref{AMK} and \eqref{KNO} to obtain another  projective
 resolution of $A$, which is dual to the first one
$$
0\to A\o A\stackrel{\pp\ve}\too M\ve\stackrel{d\ve}\too N\ve\stackrel{\bj\ve}\too 
A\o A\to 0.
$$
Thus, there exists a  morphism between the two resolutions,
that is one has the following commutative diagram of long 
exact sequences
\beq{big}
\xymatrix{
0\ar[r]& A\o A\ar[r]^<>(0.5)\bj\ar[d]^<>(0.5){h}&
N\ar@{->>}[rr]^<>(0.5){v}\ar[d]^<>(0.5){f_1}&&K\ar[d]^<>(0.5){g}
\ar@{^{(}->}[rr]^<>(0.5){u}&&M\ar[d]^<>(0.5){f}
\ar@{->>}[r]^<>(0.5){\pp}&\Om^1A\ar[r]\ar@{=}[d]^<>(0.5){\Id}&0\\
0\ar[r]& A\o A\ar[r]^<>(0.5){\pp\ve}&M\ve\ar@{->>}[rr]^<>(0.5){u\ve}&&K\ve
\ar@{^{(}->}[rr]^<>(0.5){v\ve}&&N\ve\ar@{->>}[r]^<>(0.5){\bj\ve}&\Om^1A\ar[r]&0.
}
\eeq
It is immeadiate from commutativity of the diagram that the map
$f$ takes $K$ to $K\ve$, hence, restricts to
a map $g: K\to K\ve.$
Observe also that the vertical morphism $h$ on the left 
of the diagram is given by multiplication,
via the inner bimodule structure on $A\o A,$
by some element $z=\sum z'\o z''\in A\o A.$

Clearly,
one may compute $\Ext^3(A,A\o A)$
by applying  $\Hom(-,A\o A)$ to
 any of the two resolutions in the diagram above.
The resulting self-map on
 $\Ext^3(A,A\o A)\cong A$ induced by the vertical morphism
 between the resolutions is 
a map $A\to A$ given by multiplication by 
the element $\sum z'\cdot z''\in A$.
On the other hand, this endomorphism of
 $\Ext^3(A,A\o A)$ must be equal to the identity.
Thus, we conclude that $\sum z'\cdot z''=1$.

Next, consider the  connecting homomorphism $\pa$
\beq{connect}
\Hom(K,M\ve)\stackrel{u\ve}\too\Hom(K,K\ve)
\stackrel{\pa}\too \Ext^1(K, A\o A)=A,
\eeq 
arising from the 
the short exact sequence
$A\o A\into M\ve\stackrel{\pp\ve}\onto K\ve$.
We would like to compute the elements $\pa(g)$ and $\pa(g\ve).$

To compute $\pa(g)$, we use the diagram
$$
\xymatrix{
0\ar[r]& A\o A\ar[r]^<>(0.5)\bj\ar[d]^<>(0.5){h}& N\ar[r]^<>(0.5){v}\ar[d]^<>(0.5){f_1}
&K\ar[r]\ar[d]^<>(0.5){g}& 0\\
0\ar[r]& A\o A\ar[r]^<>(0.5){\pp\ve}&M\ve\ar[r]^<>(0.5){u\ve}& K\ve\ar[r]&
0.
}
$$

The top row of this diagram is a projective resolution of $K$.
We deduce  ${\Ext^1(K, A\o A)}$ $=\Hom(A\o A,A\o A)/\Hom(N,A\o A)=(A\o A)/N^\vee.$
With this identification,
the isomorphism $\Ext^1(K, A\o A)\cong A$ goes
to the isomorphism $(A\o A)/N^\vee\iso A$ induced by the map $\bj\ve,$ see
\eqref{big}. Furthermore, the
class $\pa(g)\in \Ext^1(K, A\o A)$ corresponds, under the identification,
to the projection to $(A\o A)/N^\vee$ of the 
element $h\in \Hom(A\o A,A\o A).$ Thus, in 
 $(A\o A)/N^\vee\cong A$, we have
$\pa(g)=h=\sum z'\cdot z''=1.$

Next, we compute  $\pa(g\ve).$ To this end, we
dualize the right part of diagram \eqref{big} and obtain
a commutative diagram
$$
\xymatrix{
0\ar[r]& A\o A\ar[r]^<>(0.5){\bj}\ar@{=}[d]^<>(0.5){\Id}&
N\ar[r]^<>(0.5){v}\ar[d]^<>(0.5){f\ve}&K\ar[d]^<>(0.5){g\ve}\ar[r]&0\\
0\ar[r]& A\o A\ar[r]^<>(0.5){\pp\ve}&M\ve\ar[r]^<>(0.5){u\ve}& K\ar[r]&0.
}
$$
We deduce as before that, in 
 $(A\o A)/N^\vee\cong A$, we have
$\pa(g\ve)=\Id=1.$ 

Thus, we have shown that $\pa(g)=1=\pa(g\ve).$
Hence, we deduce from \eqref{connect} that there exists
a morphism $r: K\to M\ve$ such that
 $u\ve\ccirc r=g\ve-g$. We have
$(u\ve\ccirc r)\ve=(g\ve-g)=-(g\ve-g)=-(u\ve\ccirc r).$
Therefore,  we find
\begin{align*}
(g+\frac{1}{2}u\ve\ccirc r)\ve-
(g+\frac{1}{2}u\ve\ccirc r)&=g\ve-g+\frac{1}{2}(u\ve\ccirc r)\ve-
\frac{1}{2}(u\ve\ccirc r)\\
&=g\ve-g-(u\ve\ccirc r)=g\ve-g-(g\ve-g)=0.
\end{align*}

Replacing $g$ by $g'=g+\frac{1}{2}u\ve\ccirc r$ does not affect commutativity of 
the diagram in the statement of the Lemma, and the result follows.
\end{proof}

\subsection{Proof of Theorem \ref{BV}}\label{bv_pf}
To prove the equation stated in part (i) of Theorem \ref{BV},
choose a resolution $\SP$, as in Proposition \ref{ddpf}(i).
The Hochschild cohomology algebra,
 $H^\hdot(A,A)$, may be computed as the cohomology
of the DG algebra $(\End \SP, \ad d)$,
where we put $\End \SP:=\bigoplus_{r\in\Z}
\Hom_{\bimod{A}}(\SP,\SP[r])$ and where $d: \SP\to\SP[1]$ denotes the
differential in $\SP$.
 The  cup product
on  $H^\hdot(A,A)$ is known to be induced by the
natural  composition-product
 $ \End \SP \o   \End \SP
\to \End \SP$ in the DG algebra $\End \SP$.

Similarly,  $H_\idot(A,A)$  may be 
computed by means of the complex $(\SP\o_A   \SP, d)$, and
the contraction-action
of the algebra $H^\hdot(A,A)$ on $H_\idot(A,A)$ is known to be induced by 
the natural pairing
\beq{ppp}(\SP \o_A  \SP)\o   \End \SP\to(\SP \o_A  \SP),\quad
(p'\o p'', f)\mto p'\o f(p'').
\eeq

The morphism $f: \SP\to \SP^\vee$, see \eqref{PP}, yields 
quasi-isomorphisms
of complexes
$$\SP \o_A  \SP\cong \SP^\vee \o_A  \SP\cong \End \SP.$$
Using this, one
may reinterpret \eqref{ppp} as a   pairing
$\End \SP \o   \End \SP
\to \End \SP$. The latter map is nothing but the composition
map, and the statement of Theorem \ref{BV}(i)  follows.

To prove part (ii) of the theorem, we recall that for any
$\eta\in H^p(A,A)$ one can define
a {\em  Lie derivative} operation
$L_\eta: H_j(A,A)\to H_{j-p+1}(A,A),$ cf. eg. \cite{TT}.
The  Lie derivative and contraction operators
 satisfy all the
standard identities well-known from differential geometry.
In particular, for any $\xi, \eta\in H^\hdot(A,A),$
one has $i_{[\xi,\eta]}=[L_\xi, i_\eta],$
where $[-,-]$ always stands for the {\em super}-commutator.
Further, according to \cite{Re}, one has 
the following analogue of the Cartan identity
$L_\eta:= [{\mathsf{B}},i_\eta]$.
Using this identity and the fact that $i_\xi i_\eta=i_{\xi\cup\eta}$,
 one can rewrite the equation above in the 
form
\begin{align}\label{stand}
i_{[\xi,\eta]}&=L_\xi i_\eta-i_\eta  L_\xi &
 \big(\forall\xi, \eta\in H^\hdot(A,A)\big)\nonumber
\\
&=i_\xi {\mathsf{B}}  i_\eta+ {\mathsf{B}} i_\xi i_\eta- i_\eta {\mathsf{B}} i_\xi- {\mathsf{B}} i_\eta i_\xi&\\
&=i_\xi {\mathsf{B}} i_\eta- i_\eta {\mathsf{B}} i_\xi+{\mathsf{B}} i_{\xi\cup\eta}-
i_{\xi\cup\eta} {\mathsf{B}}.&\nonumber
\end{align}

Now, let $c\in H_d(A,A)$ be the image of
$1\in H^0(A,A)$ under the isomorphism
${\mathbb D}\inv: H^0(A,A)\iso H_d(A,A)$, cf. \eqref{vdb}.
We apply the operation on each side of \eqref{stand} to 
the element $c$.
Further,
applying the isomorphism ${\mathbb D}$
to the  resulting equation, we get
\beq{result}\SH(i_{[\xi,\eta]}c)=\SH(i_\xi {\mathsf{B}} i_\eta c)
-\SH( i_\eta {\mathsf{B}} i_\xi c)+\SH({\mathsf{B}} i_{\xi\cup\eta}c)-
\SH(i_{\xi\cup\eta} {\mathsf{B}} c).
\eeq

Next, we  use the definition of the
operator $\Delta$ and the statement of part (i) of the theorem
to rewrite \eqref{result} as follows
$$[\xi,\eta]\cup {\mathbb D}(c)=\xi\cup\Delta(\eta\cup{\mathbb D}(c))
-\eta
\cup\Delta(\xi\cup\SH(c))+\Delta(\xi\cup\eta\cup{\mathbb D}(c))-
\xi\cup\eta\cup\Delta({\mathbb D}(c)).
$$

But we have $\SH(c)=1$ and, clearly,
$\Delta(1)=0$. Thus, the rightmost summand in this formula vanishes
and we see that the above equation becomes 
the  identity from part (ii) of the theorem.
\qed

\begin{rem} The proof above applies without change to yield
the standard BV identity in $\La^\hdot\T(X)$
for a Calabi-Yau manifold $X$, cf.  Example \ref{cyx}.
In that case, one can of course apply
the statement of Theorem \ref{BV}  directly,
as well.
\end{rem}

\subsection{Proof of Proposition \ref{kunneth}}\label{prop_pf}
 Let $R$ be a finite dimensional semisimple algebra and
$(D=\bigoplus_{r\geq 0} 
D_r,\d)$
 a DG $R$-algebra.

\begin{defn}\label{sf} We call a left \DG $D$-module {\em d-projective}
if it  can be obtained from $D$, viewed as a rank 1 free left  \DG
$D$-module,
by repeated application of the following operations

\npb{Taking arbitrary direct sums of \DG modules;}

\npb{Taking \DG direct summands;}

\npb{Degree shifts $M \mto M[k],$ for any $k\in\Z$.}
\end{defn}

The following result is clear

\begin{lem}\label{claim2}  Let 
$M$ be a d-projective left
\DG $D$-module.
Then, 

\vi For any quasi-isomorphism $f:N\stackrel{\oper{qis}}\too N'$ of
right \DG $D$-modules,
the 
induced map  $f\o\Id_M: N\o_D M\stackrel{\oper{qis}}\too
N'\o_D M$ is  a quasi-isomorphism. 

\vii Any morphism $g: M'\to M,$ of  left
\DG $D$-modules, that induces an isomorphism on cohomology
has a quasi-inverse $h:M\to M',$ a morphism of \DG $D$-modules
such that the maps $g\ccirc h$ and $h\ccirc g$ are both homotopic
to the identity.\qed
\end{lem}

Let  $B$ be an
$R$-algebra. Given $0\to P_n\to P_{n-1}\ldots\to P_0\to B\to 0$, a bounded
complex of 
$B$-bimodules, let  $\SP=\bigoplus_{r\geq 0}P_r$.
This is a DG bimodule with an obvious grading and differential.
The tensor algebra, $D:=T_B\SP$, 
acquires a DG algebra structure, with the grading 
$D=\bigoplus_{r\geq0}D_i$ and differential
$\d: D_\idot\to D_{\idot-1}$, both being
induced from  $\SP$.
Let $D_+:=\bigoplus_{r>0} 
D_r$ denote the augmentation ideal of $D$, and put $I=\d(D_1)$,
a  two-sided in $D_0=B$. The sum
$I+D_+$ is clearly a graded  $\d$-stable two-sided ideal in the algebra $D$.

Observe also that any power of a $\d$-stable  two-sided ideal in $D$
is again  a $\d$-stable  two-sided ideal.

\begin{lem}\label{ideal} Assume that  $B$ is a smooth
algebra and that  $P_r$ 
is a finitely generated
 projective $B$-bimodule for any $r=0,\ldots,n$.

Then, the two-sided ideal
$(I+D_+)^m$ is a d-projective left  \DG  $D$-module,
for  any $m\geq 1.$
\end{lem}

The proof of the lemma is based on the fact that, for a smooth algebra
$B$, any $B$-submodule of a projective left $B$-module is
projective. In particular, $I^m$ is a  projective left $B$-module,
for any  $m\geq 1.$ 

Now, fix $m\geq 1$ and put ${\mathbf J}:=(I+D_+)^m$.
One proves that ${\mathbf J}/D_+\cdot{\mathbf J}$ is a projective graded
$B$-bimodule. Furthermore, one constructs
an isomorphism $D\o_B({\mathbf J}/D_+\cdot{\mathbf J})\iso {\mathbf J}$,
of left DG $D$-modules. We leave details to the reader.\hfill$\lozenge$

\begin{proof}[Proof of  Proposition \ref{kunneth}]
To prove the implication
 (i)$\en\Rightarrow\en$(ii), we
consider a natural commutative diagram
\beq{pda}
\xymatrix{
&&\Om^1_R\FD\ar[dll]_<>(0.5){p\o\Id\o p}\ar[drr]^<>(0.5){p\o p}&&\\
A\o_R\Om^1_R\FD\o_RA\ar@{=}[rr]&&\Om^1_R(\FD|A)\ar[rr]^<>(0.4){p_{\FD|A}}&&\Om^1_RA.
}
\eeq
Here, we view $\Om^1_RA$ as a DG vector space with zero
differential, so the maps in  \eqref{pda} 
are morphisms  of DG vector space.
We claim that all these morphisms are  quasi-isomorphisms.

In \eqref{pda}, the map $p\o p$ acts as $u\, dv\mto p(u)\, dp(v)$, where we identify
the space of 1-forms with a tensor product via the canonical 
isomorphisms
$\FD\o_R (\FD/R)\iso \Om^1_R\FD,\, u\o v\mto u\, dv,$ resp. $A\o_R  (A/R)\iso \Om^1_RA$.
We know that $\d(R)=0$ and  that $\d(\FD)\cap R=0,$ by our assumptions.
Since $R$ is a semisimple
finite dimensional algebra, we may use the K\"unneth formula to conclude
that the map $p\o p$ is a  quasi-isomorphism.

To prove that the map $p\o\Id\o p$ in diagram  \eqref{pda}
 is a  quasi-isomorphism,
we apply Lemma \ref{claim2}(i) to the algebra $D:=\FD^e,$ cf. Notation \ref{e}.
We view $M:=\Om^1_R\FD$ as a left DG $D$-module.
This  DG module is d-projective since  $\FD$ is a  smooth DG algebra,
by  assumption.
Applying Lemma \ref{claim2} to
the quasi-isomorphism $f=p\o p^\opp:\FD^e\qis A^e,$
yields a quasi-isomorphism
$$
\Om^1_R\FD=D\o_D\Om^1_R\FD=\FD^e\o_D\Om^1_R\FD
\qis A^e\o_D\Om^1_R\FD=
\Om^1_R(\FD|A).
$$
The composite above is  the map  $p\o\Id\o p$ in \eqref{pda}.

Thus, we have shown that the maps
$p\o p$ and $p\o\Id\o p$ are both quasi-isomorphisms.
We conclude by commutativity of  \eqref{pda} that
 the third
map, $p_{\FD|A}$, is a quasi-isomorphism as well.
We may express this by saying that the DG module
$\cone(p_{\FD|A})$,
the cone of the DG morphism $\Om^1_R(\FD|A)\to \Om^1_RA,$ is {\em acyclic}.

We now prove the implication (ii)$\en\Rightarrow\en$(i).
To this end, let $\FI:=I + \bigoplus_{r>0}\FD_r$.
Thus, $\FI\sset \FD$ is a graded $\d$-stable two-sided ideal,
and $A=\FD/\FI$. Therefore, from the cotangent 
exact sequence \eqref{CQ} for the pair $\FI\sset\FD$, using the quasi-isomorphism
$\Om^1(\FD|A)\qis\Om^1A,$ we deduce that the
complex $(\FI/\FI^2, \d)$ is acyclic.

Let $\FD\supset\FI\supset \FI^2\supset\ldots,$ be the  $\FI$-adic
filtration on $\FD.$ This is a  filtration by $\d$-stable two-sided
ideals. Hence, one can form $\bigoplus_q \FI^q/\FI^{q+1},$
an associated graded DG algebra, and there is a standard
 spectral sequence
$$E^2_{p,q}=H_p(\FI^q/\FI^{q+1})\quad\Longrightarrow\quad
E^\infty_{p,q}=\FI^q\cd H_p(\FD)/\FI^{q+1}\cd H_p(\FD).
$$

Recall next that there are two cases in the statement of 
Proposition \ref{kunneth}. In the first case  the algebra
$\FD$ is assumed to have an additional grading.
In the second case
each homogeneous component, $\FD_p, \,p=0,1,\ldots,$ of 
the algebra $\FD$, is assumed to be complete in 
$I$-adic topology. Either of the two assumptions insures
that the  spectral sequence above
converges. Hence, the vanishing of
the $E^2_{p,q}$-terms with $p+q>0$
would imply  that $H_p(\FD)=0$ for all $p> 0.$

To prove the $E^2$-vanishing we will exploit
the idea of the proof of \cite{CQ}, Proposition 5.2.
Specifically,  it follows from Lemma \ref{ideal}  that
$\FI^q$ is a d-projective left DG $\FD$-module,
for each $q=1,2,\ldots.$ In particular, $M=\FI^q$ is a
flat  left  $\FD$-module. Following \cite{CQ}, we observe that for any
flat  $\FD$-module $M$, the canonical map
$(\FI/\FI^2)\o_\FD M\to (\FI\cdot M)/(\FI^2\cdot M)$ is an isomorphism.
Therefore, in our case, we deduce that multiplication in the
algebra $\FD$ gives  an isomorphism
$(\FI/\FI^2)\o_\FD \FI^q\to \FI^{q+1}/\FI^{q+2}$, for any   $q\geq 0$.
Hence, applying Lemma \ref{claim2} in the case where $M=\FI^q$ and where
$f: \FI/\FI^2\qis0$ is the zero map,
we deduce that, for any $q\geq 1$, one has $H_p(\FI^{q+1}/\FI^{q+2})=
H_p\big((\FI/\FI^2)\o_\FD \FI^q\big)=0.$
It follows that $E^2_{p,q}=0$ for all $p+q>0,$
and we are done.
\end{proof}

\medskip

\noindent
Department of Mathematics, University of Chicago,
Chicago, IL 60637, USA\\
{\tt ginzburg@math.uchicago.edu}
\end{document}